\documentclass[12pt]{article}

\usepackage{amsfonts,psfig}

\newcommand{\text}[1]{\mbox{#1}}

\newcommand{\limsuppn}{\limsup_{n\to+\infty}{1\over n}}

\newcounter{main}

\newtheorem{theorem}{Theorem}[section]
\newtheorem{proposition}[theorem]{Proposition}
\newtheorem{lemma}[theorem]{Lemma}

\newtheorem{property}{Property}[section]
\newtheorem{remark}{Remark}[section]
\newtheorem{definition}{Definition}[section]
\newtheorem{maintheorem}{Theorem}
\newtheorem{claim}{Claim}[section]

\newcommand{\blanksquare}{\,\,\,$\sqcup\!\!\!\!\sqcap$}
\newenvironment{proof}{{\flushleft {\bf Proof: }}}{\blanksquare}
\newenvironment{rproof}[1]{{\flushleft{\bf Proof of #1: }}}{\blanksquare}
\newcounter{example}
\newenvironment{example}%
{{\stepcounter{example}}{\flushleft {\bf Example \arabic{example}:}}}%
{\par}

\newcommand{\al} {\alpha}       
\newcommand{\be} {\beta}        
\newcommand{\ga} {\gamma}    
\newcommand{\de} {\delta}       \newcommand{\De}{\Delta}

\newcommand{\ep} {\epsilon}     

\newcommand{\ze} {\zeta}
\newcommand{\te} {\theta}

\newcommand{\la} {\lambda}      \newcommand{\La}{\Lambda}

\newcommand{\si} {\sigma}       
\newcommand{\ro}{\rho}

\newcommand{\ups}{\upsilon}     
\newcommand{\vfi}{\varphi}
\newcommand{\om} {\omega}       \newcommand{\Om}{\Omega}

\newcommand{\vecto}[2]{(#1_1,\ldots,#1_{#2})}

\newcommand{\deter}{{\rm det\,}}
\newcommand{\diame}[1]{{\rm diam\,}(#1)}
\newcommand{\diff}{{\rm Diff\,}}
\newcommand{\dimns}{{\rm dim\,}}

\newcommand{\inter}[1]{{\rm int\,}(#1)}

\newcommand{\refer}[1]{(\ref{#1})}

\newcommand{\slope}[1]{{\rm slope\,}(#1)}

\newcommand{\supp}{{\rm supp\,}}

\newcommand{\natur}{{\mathbb N}}
\newcommand{\real} {{\mathbb R}}
\newcommand{\difsim}{{\bigtriangleup}}

\newcommand{\relativ}{{\mathbb Z}}
\newcommand{\complex}{{\mathbb C}}
\newcommand{\torus}{{\mathbb T}}

\newcommand{\un}{\underline}

\newcommand{\pa}{\partial}
\newcommand{\til}{\tilde}

\newcommand{\mapto}{\longrightarrow}
\newcommand{\impl}{\Longrightarrow}

\newcommand{\equi}{\Longleftrightarrow}

\newcommand{\SA}{{\cal A}}
\newcommand{\SB}{{\cal B}}
\newcommand{\SC}{{\cal C}}
\newcommand{\SD}{{\cal D}}

\newcommand{\SF}{{\cal F}}
\newcommand{\SG}{{\cal G}}
\newcommand{\SH}{{\cal H}}

\newcommand{\SM}{{\cal M}}

\newcommand{\SO}{{\cal O}}

\newcommand{\SQ}{{\cal Q}}
\newcommand{\SR}{{\cal R}}
\newcommand{\Ss}{{\cal S}}

\newcommand{\SU}{{\cal U}}
\newcommand{\SV}{{\cal V}}

\newcommand{\subl}{\underline}
\newcommand{\tra}{\overline}

\newcounter{contador}
{\begin{list}{#1}{\usecounter{contador}}}%
{\end{list}}

{\begin{list}{#1}{\usecounter{#2}}}%
{\end{list}}

\newenvironment{nliste}[1]%
{\begin{list}{#1}%
{\usecounter{contador}\setlength{\leftmargin}%
{0.5cm}\setlength{\labelsep}{0.1cm}\setlength{\itemsep}%
{0cm}\setlength{\parsep}{0.1cm}}}%
{\end{list}}

\newcommand{\Parentalgarismo}{{(\arabic{contador})}}

\newcommand{\ai}{{\'{\i}}}

\newcommand{\nuinf}{{\nu^\infty}}

\newcommand{\domcomplinv}{( \SU_0,\ldots,\SU_{r-1} )}
\newcommand{\domcomplinvv}{( \SU_0',\ldots,\SU_{r'-1}' )}
\newcommand{\domcomplinvl}{( \SU_0',\ldots,\SU_{r-1}' )}
\newcommand{\domcomplinvtil}{( \til{\SU}_0,\ldots,\til{\SU}_{r-1})}

\newcommand{\dominvu}{ \SU_0, \ldots, \SU_{r-1} }

\title{Attractors and Time Averages
for Random Maps}

\author{V\'{\i}tor Ara\'ujo}

\date{November 1998}

\begin{document}
\maketitle

\begin{abstract}

Considering random noise in finite dimensional
parameterized  families
of diffeomorphisms of a compact finite dimensional
boundaryless manifold $M$, we show the existence of
time averages for almost every orbit of each point of $M$,
imposing mild conditions on the families; see
section~\ref{teoremaprincipal}.
Moreover these averages are given by a finite
number of physical absolutely continuous
stationary probability measures.

We use this result to deduce that situations with
infinitely many sinks and H\'enon-like
attractors are not stable under
random perturbations, e.g., Newhouse's and Colli's
phenomena in the generic unfolding of a quadratic homoclinic
tangency by a one-parameter family of diffeomorphisms.

\smallskip

\noindent
Key-words: random perturbations, time averages,
physical probabilities, homoclinic bifurcations.

\bigskip

\begin{center}
{\bf Resum\'e}
\end{center}

\medskip

On consid\'ere bruit al\'eatoire dans
des familles param\`etriz\'es de
dimension finite de diff\'eomorphismes
d'une vari\'et\'e compacte $M$, de dimension
finite et sans bord, et montre l'existence
de moyennes temporelles assymptotiques pour
presque toute orbite de chaque point de $M$,
en imposant des conditions pas tr\'es forts
sur les familles; v. section~\ref{teoremaprincipal}.

Cettes moyennes sont d\'efinies par un nombre
finite de mesures de probabilit\'e stationnaires
physiques et absolument continues.

On utilise ce resultat pour d\'eduire que
les situations de coexistence d'une infinit\'e
de puis et d'attracteures de type H\'enon
ne sont pas stables par
perturbations al\'eatoires, e.g.,
les ph\'enom\`enes de Newhouse et Colli
dans le d\'edoublement g\'en\'erique
d'une tangence homoclinique quadratique
par une famille
de diff\'eomor\-phismes \`a un param\`etre.

\smallskip
\noindent
Mots-cl\'es: perturbations al\'eatoires,
moyennes temporelles, probabilit\'es physiques,
bifurcations homocliniques.

\end{abstract}

\pagebreak

\tableofcontents

\pagebreak

\section{Introduction}
\label{intro}

Newhouse proved in~\cite{N1,N2,N3} that many
surface diffeomorphisms have infinitely many
attracting periodic orbits (sinks),
a serious blow to early hopes that generic
systems might have only finitely many attractors.
Indeed, see~\cite{N3} and also~\cite{PT},
arbitrarily close to
any $C^2$ diffeomorphism
on a surface $M$ with a homoclinic tangency
there exist open subsets of ${\rm Diff}^2(M)$
whose generic elements have infinitely many
sinks or sources.

This result was extended to arbitrary dimensions
by Palis-Viana in~\cite{PV},
see also~\cite{Ro} and~\cite{GST}.
Diffeomorphisms with infinitely many coexisting
hyperbolic attractors were constructed by
Gambaudo-Tresser in~\cite{GT}. Colli
showed in~\cite{C}
that diffeomorphisms displaying infinitely
many H\'enon-like strange attractors
are dense in some open subsets of
${\rm Diff}^\infty(M)$, if $\dimns{M}=2$.
Even more recently, Bonatti-D{\ai}az in~\cite{BD}
showed that coexistence of infinitely many sinks
or sources is generic in some open subsets of
${\rm Diff}^1(M)$, if $\dimns{M}\ge3$.

However, apart from these existence results,
diffeomorphisms with infinitely many  attractors
or repellers are still a mystery.
Results of~\cite{Ma}, \cite{DPU}, \cite{BDP}
show that maps which cannot be approximated
by others with infinitely many sinks or sources
have properties of partial hyperbolicity.
In this case the dynamics of these maps can
be understood to some degree,
see e.g.~\cite{BP},\cite{PS},\cite{GPS},\cite{BV},\cite{ABV}.
It would be nice to know that systems
with infinitely many sinks or sources are
negligible from the measure theoretical
point of view. Indeed,
it has been conjectured that such systems
correspond to zero Lebesgue measure
in parameter space for generic families
(finite number of parameters) of maps,
see~\cite{TY} and~\cite{PT}.
Nevertheless this is not yet know.

Here we show that this phenomenon of
coexistence of infinitely many sinks or
sources can indeed be discarded in the
setting of maps endowed with random noise.
We prove that (Theorem 1)
{\em every diffeomorphism of a compact
finite dimensional boundaryless manifold $M$
under absolutely continuous random perturbations
along a parameterized family
has only finitely many physical measures
whose basins cover Lebesgue-a.e. point of $M$\/}.

\medskip

In the context of the generic unfolding
of quadratic homoclinic tangencies by uniparametric
arcs of surface diffeomorphisms,
where the coexistence phenomenon of
infinitely many attractors
was first shown to occur,
we prove (Theorem 2) a result similar to
the previous one concerning points whose
perturbed orbits visit a neighborhood of
the tangency infinitely often with positive
probability, which we call {\em recurrent points\/}.

This result is a corollary of the former since
we show the random parametric perturbations 
applied on the recurrent points to be
{\em absolutely continuous\/} as well.
For an uniparametric arc to satisfy this
property in a surface
a quadratic homoclinic tangency is used:
the mixture of expanding and contracting
directions near a homoclinic tangency point,
in a neighborhood of it in the manifold for
every diffeomorphism close to the one
exhibiting the tangency,
is what permits us to get absolute
continuity even when only a single
parameter is at hand.

We conclude (section~\ref{sinks})
that {\em there cannot be infinitely
many attractors (or physical measures) whose
orbits (resp. supports) pass near
a quadratic homoclinic tangency point or
its generic unfolding under random parametric
perturbations (i.e. random errors in the parameters)\/}
---
{\em in this sense, diffeomorphisms with infinitely
many attractors are not stable under
random perturbations\/}.

\medskip

These results can be seen from the
perspective of a
broad program proposed by J.Palis in~\cite{Pa}.
In particular, he conjectured that systems
with finitely many attractors are dense in
the space of all systems. Moreover,
these attractors should have nice statistical
properties, including existence of physical
measures supported on them, and stochastic
stability under small random noise ---
see e.g.~\cite{V2}.

\medskip

Fornaess and Sibony in~\cite{FS} have shown a
result similar to Theorem~\ref{teorema1} to hold
in the context of random perturbations of rational
functions.
The precise form of the statement of this theorem and
of some definitions was inspired on theorem 1.1 of
theirs.

\medskip

Relevant setting and all definitions
are in sections~\ref{notacoes}
and~\ref{invariantdomains}
along with the
precise statement of the result, including
the kind of noise to be used and
some examples.
A summary of the steps of the proof is
given in section~\ref{esquema},
where we also sketch the contents
of sections 5 through 9.
In section~\ref{Bowen} we apply our results
to perturbations of an example of Bowen.
This provides a good insight into the
meaning of these results.

\medskip

Relevant settings, definitions and the statement
of Theorem~\ref{teorema2} are
in section~\ref{bifurcations}.
Its proof in sections~\ref{fisica:prova}
and~\ref{regular:prova}.

\medskip

Several questions arise in this context of
systems with random noise and the simple methods
used in this work to derive theorems 1 and 2
should be generalized and extended.
Some of those questions are presented in the last
section (section~\ref{conjecturas}) of this paper.

\medskip

{\bf Aknowledgements:\/}
I thank {\em Instituto de Matem\'atica Pura e Aplicada\/}
(IMPA) at Rio de Janeiro, where this paper has been written,
for its excellent research atmosphere and facilities.
Many thanks to Marcelo Viana for uncountable discussions
and suggestions, to Eul\'alia Vares for clarifying lemma
\ref{lema0-1} to me and to JNICT/FCT~\footnote{Subprograma
Ci\^encia e Tecnologia do $2^{\subl{o}}$
Quadro Comunit\'ario de Apoio -- PRAXIS XXI/BD/3446/94}
(Portugal) and IMPA/CNPq (Brasil)
for partial financial support,
which made this work possible.

\section{Some Notations, Definitions\\and the Main Theorem}
\label{notacoes}

Throughout this paper 
$M$ will signify a compact boundaryless manifold with finite dimension,
$m$ will be some normalized ($m(M)=1$) Riemannian volume form
on $M$ and $d_M:M\times M\mapto \real$ a distance given by
some Riemannian structure on M,
fixed once and for all.
When not otherwise mentioned, absolute continuity will
be taken with respect to the probability $m$.

The random perturbations to be considered will act
on the dynamics of diffeomorphisms of a parameterized
family given by the $C^1$ function
$f:M\times B^n \mapto M$, where
$B^n=\{ y\in\real^n: \| y \|_2<1 \}$ is the
unit ball of $\real^n$, $n\ge1$,
$\| \cdot \|_2$ is the Euclidean norm and
the map $f_t:M\mapto M,\; x\in M \mapsto f(x,t)$ is
a diffeomorphism for every $t\in B^n$.

\subsection{Perturbations around a Parameter}
\label{pertaroundpar}

Let us fix $a\in B^n$ and take $\ep>0$ such that
the closed $\ep$-neighborhood of $a$ be contained
in $B^n$, $\tra{B}^n(a,\ep)\subset B^n$.
We define the {\em perturbation space around $a$ of
size $\ep$\/} to be
\[
\De=\De_{\ep}(a)={\tra{B}^n(a,\ep)}^\natur=
\left\{
\subl{t}=(t_j)_{j=1}^\infty : \| t_j-a \|_2 \le \ep, j\ge1
\right\}
\]
with the product topology, which is equivalent to the
topology induced by the metric
$
d( \subl{t},\subl{s} ) =
\sum_{j=1}^\infty 2^{-n}\cdot\| t_j-s_j \|_2,\;\;
\subl{t},\subl{s}\in\De,
$
and the measure $\nuinf$ given by the product of the
normalized Lebesgue volume measure $\nu$ over each $\tra{B}^n(a,\ep)$.
For sets $A_1,\dots,A_k$ of the Borel family in $\tra{B}^n(a,\ep)$ we
have
$
\nuinf( A_1\times\dots\times A_k\times {\tra{B}^n(a,\ep)}^\natur)=
\nu(A_1)\dots\nu(A_k)
$
and if $A\subset \tra{B}^n(a,\ep)$ then
$\nu(A)=|\tra{B}^n(a,\ep)|^{-1}\cdot |A|$,
where $|A|$ will mean
the Lebesgue volume measure of $A$.

Now we define the {\em perturbed iterates\/} of $f$ by
\[
f_{\subl{t}}^k(z) = 
	f^k(z,\subl{t}) = f_{t_k}\circ\dots\circ f_{t_1}(z),\;\;
z\in M,\, \subl{t}\in\De
\]
and state the useful convention that $f^0(z,\subl{t})=z$ and
\[
f_{V}^k(U) =
	f^k (U, V) = 
\{ f_{\subl{t}}^k (z) : \subl{t}\in V,\, z\in U \},\;\;
U\subset M,\, V\subset \De
\]
for every $k\ge1$.
We emphasize a very often used property in what
follows.
\begin{property}\label{prop0}
For every fixed $k\ge1$ it holds that
\begin{enumerate}
\item $(z,t_1,\ldots,t_k)\in M\times
\tra{B}^n(a,\ep)\times
\stackrel{k}{\ldots}
\times \tra{B}^n(a,\ep)
\mapsto f^k(z,t_1,\ldots,t_k)=
f_{t_k}\circ\dots\circ f_{t_1}(z) \in M$
is differentiable;

\item $(z,\subl{t}) \in M\times\De \mapsto
f^k (z,\subl{t})\in M$ is continuous (with the
product topology);

\item $z\in M\mapsto f^k(z,t_1,\ldots,t_k)\in M$
is a diffeomorphism for every $t_1,\ldots,t_k\in \tra{B}^n(a,\ep)$.
\end{enumerate}
\end{property}

Given $\subl{t}\in\De$ and $z\in M$
we will call $\{f_{\subl{t}}^k(z)\}_{n=1}^\infty$ the
{\em $\subl{t}$-orbit\/} of $z$ and many times write
$\SO(z,\subl{t})$.

\medskip

In this way, perturbations are implemented by a random
choice of parameters of a parameterized family
of diffeomorphisms
at each iteration,
the choice being made in a $\ep$-neighborhood
of a fixed parameter  according to a uniform probability.
Such choices are represented by a vector $\subl{t}$ in $\De$,
an infinite product of intervals, and the greater or
lesser importance of the set of perturbations taken into
account will be evaluated by the measure $\nuinf$.

\medskip

This kind of random iteration will be referred to
as {\em parametric noise\/}.
With the settings given above, the family of
diffeomorphisms acting on $M$ with parametric
noise of level $\ep$ around $f_a$ will be
written $\SF_{a,\ep} = \{ f_t: t\in\tra{B}^n(a,\ep) \}$.
To simplify writing the factors of $\De$ we set
$T=\tra{B}^n(a,\ep)$ from now on, so that
$\De=T^\natur$.

\subsection{Stationary Probabilities}
\label{probestacsec}

We can define a shift operator
$
S:M\times\De \mapto M\times\De, \;
(z,\subl{t}) \mapsto ( f_{t_1}(z) , \si(\subl{t}) ),
$
where $\si$ is the left shift on sequences of $\De$:
$
\si(\subl{t})=\subl{s} \quad \text{with} \quad
\subl{t}=(t_1,t_2,t_3,\ldots) \;\; \text{and} \;\;
\subl{s}=(t_2,t_3,t_4,\ldots).
$
By the definition of $S$ and property \ref{prop0}(2)
we deduce that $S$ is continuous.

\medskip

A probability measure $\mu$ in $M$ is said
a {\em stationary probability\/} if the measure
$\mu\times\nuinf$ is $S$-invariant:

\begin{equation}\label{probabilestacionaria}
\mu\times\nuinf \; (S^{-1} A) = \mu\times\nuinf
\;(A),
\;\; \text{for every Borel subset}\;
A \; \text{of} \;\; M\times\De.
\end{equation}

\smallskip

This is equivalent to say that $\mu$ satisfies
the following identity
\begin{equation}\label{probestac}
\int \int \vfi(f(z,t)) \, d\mu(z) d\nu(t) =
\int \vfi (z) \, d\mu(z), \; \forall \vfi \in C^0(M).
\end{equation}

\medskip

In fact, writing \refer{probabilestacionaria}
for $A=U\times\De$, where $U$ is a Borel subset of $M$,
we have
\begin{eqnarray}
\mu\times\nuinf \; (S^{-1} (U\times\De))
&=&
\mu\times\nuinf \;
\left(
\bigcup_{s\in T} f_s^{-1} (U) \times \{s\} \times \De
\right) \nonumber
\\
&=&
\mu\times\nu \;
\left(
\bigcup_{s\in T} f_s^{-1} (U) \times \{s\}
\right)
\times \nuinf (\De) \nonumber
\\
&=&
\int \int 1_U( f(x,s) ) \, d\mu(x) \, d\nu(s) \label{indicadora}
\end{eqnarray}
which is equal to
$\mu\times\nuinf \: (U\times\De) = \mu(U)$,
that is,
\[
\int \int 1_U(f(x,s)) \, d\mu(x) \, d\nu(s) = \mu(U)
= \int 1_U(x) \, d\mu(x),
\]
where $1_U$ is such that $1_U(x)=1$ if $x\in U$ and
$1_U(x)=0$ otherwise.
Then \refer{probestac} holds for every
$\vfi\in L^1_\mu(M,\real) \supset C(M,\real)$,
because simple functions are dense in $L^1_\mu$
and the relation \refer{probestac} is linear.

Conversely, if \refer{probestac} holds for every
$\vfi\in C^0(M,\real)$, then it holds for every element
of $L^1_\mu(M,\real)$ because $\mu$ and $\nu$ are
Borel measures and $f:M\times B\mapto M$ is
continuous (so that the left hand side of
\refer{probestac} gives a regular measure
over $M$).
In particular, it holds
for $\vfi=1_U$, and \refer{indicadora} is equal to
$
\int 1_U(x)\, d\mu(x) = \mu(U) = \mu\times\nuinf \;
(U\times\De)
$
proving that \refer{probestac} implies
$\mu\times\nuinf (S^{-1} (U\times\De)) =
\mu\times\nuinf (U\times\De)$. Now we see that,
if $V\subset\De$ is also a Borel subset,
\begin{eqnarray*}
\mu\times\nuinf \; (S^{-1} (U\times V))
&=&       
\mu\times\nuinf \;
\left(    
\bigcup_{s\in T} f_s^{-1} (U) \times \{s\} \times V
\right)
\\        
&=&       
\mu\times\nu \;
\left(    
\bigcup_{s\in T} f_s^{-1} (U) \times \{s\}
\right)   
\times \nuinf (V) 
\\        
&=&       
\int \int 1_U( f(x,s) ) \, d\mu(x) \, d\nu(s)
\times \nuinf (V)
\\
&=&
\int 1_U(x)\, d\mu(x) \times \nuinf (V) =
\mu\times\nuinf( U\times V)
\end{eqnarray*}
proving the equivalence between \refer{probestac}
and \refer{probabilestacionaria}.

\subsection{Ergodicity, Generic Points, Ergodic Basin}
\label{baciaergodica}

In the same way we have defined a stationary probability,
by utilizing the shift $S$, we will say that $\mu$
is a {\em stationary ergodic\/} probability measure
if $\mu\times\nuinf$ is $S$-ergodic.

\medskip

In this situation, Birkhoff's ergodic theorem ensures that
$
\lim_{n\to\infty} \frac1n \sum_{j=0}^{n-1} \psi(S^j(x,\subl{t}))
= \int \psi \, d(\mu\times\nuinf)
$
for $\mu\times\nuinf$-a.e.$(x,\subl{t})\in M\times\De$
and for every $\psi\in C^0(M\times\De,\real)$.
In particular, putting $\psi=\vfi\circ\pi$,
with $\vfi\in C^0(M,\real)$ and $\pi:M\times\De\mapto M$
the projection on the first factor, we obtain
$\psi( S^j(x,\subl{t}) ) = \vfi( f^j(x,\subl{t}) ),
\; j=0,1,2,\ldots$ and
$\int \psi \, d(\mu\times\nuinf) = \int \vfi \, d\mu$,
thus for every continuous $\vfi:M\mapto\real$
\begin{equation}\label{mediaorbital}
\lim_{n\to\infty} \frac1n \sum_{j=0}^{n-1}
\vfi( f^j(x,\subl{t}) ) = \int \vfi \, d\mu, \;\;
\text{for} \; \mu\times\nuinf-\text{a.e.} \;
(x,\subl{t})\in M\times\De.
\end{equation}

We now remark that, because $\mu\times\nuinf$
is a product measure, we have the following property.
Let $X$ be the set of $(x,\subl{t})$ that satisfy
\refer{mediaorbital} for every continuous function
$\vfi:M\mapto\real$.
If $\{\vfi_n\}_{n=1}^\infty$ is a denumerable and
dense sequence in $C^0(M,\real)$ and $X_n$ the
set of those points $(x,\subl{t})\in M\times\De$
that satisfy \refer{mediaorbital} for
$\vfi_n$, $n\ge1$, then it is easy to see
(cf.~\cite[Capt. II.6]{Mn}) that
$X=\cap_{n\ge1} X_n$ is a set of $\mu\times\nuinf$-measure 1.
Let us consider now
$X(x)=\{ \subl{t}\in\De: (x,\subl{t})\in X \}$,
the section of $X$ through $x\in M$.
Then we have $\nuinf (X(x)) = 1$
for $\mu$-a.e. $x\in M$.
Indeed, by Fubini's theorem,
$
\mu\times\nuinf (X) = \int \nuinf( X(x) )
\, d\mu(x) =1
$
with $0\le\nuinf (X(x)) \le 1$
for every $x\in M$. Hence, the last identity
implies the statement, because
$\mu$ is a probability measure.

\medskip

The points $x$ that satisfy $\nuinf( X(x)) = 1$,
that is, for which the limit in \refer{mediaorbital}
exists and equals $\int \vfi\,d\mu$ for
$\nuinf$-a.e. $\subl{t}\in\De$ and every continuous
$\vfi:M\mapto\real$,
will be called {\em $\mu$-generic points\/}.
The set of $\mu$-generic points, when $\mu$ is
stationary and ergodic, will be the
{\em ergodic basin\/} of $\mu$ and will
be written $E(\mu)$.

\medskip

To complete this setting of terms and symbols,
those ergodic stationary probability measures $\mu$
whose basin has positive volume,
$m( E(\mu))>0$, will be called
{\em physical measures\/} of the perturbed system.
We also convention to write $f^k(x,\nuinf)$
for the push-forward of $\nuinf$ by $f^k(x,\cdot)$,
that is
$
f^k(x,\nuinf) \vfi = \int \vfi( f^k(x,\subl{t}) )
\, d\nuinf(\subl{t})
$
for every
$
k\ge1, \; x\in M
$
and
$
\vfi\in C^0(M,\real)
$
by definition.

\subsection{Statement of the Results}
\label{teoremaprincipal}

\begin{maintheorem}\label{teorema1}
Let $f:M\mapto M$ be a diffeomorphism of
class $C^r$, $r\ge1$, of a compact connected
boundaryless manifold $M$ of finite dimension.
If $f=f_a$ is a member of a parametric
family under parametric noise of level $\ep>0$,
as in subsection \ref{pertaroundpar},
that satisfies the hypothesis:
there are $K\in\natur$ and $\xi_0>0$
such that, for all $k\ge K$ and $x\in M$
\begin{enumerate}
\item[A)] $f^k( x,\De) \supset B( f^k(x) , \xi_0 )$;
\item[B)] $f^k(x,\nuinf) \ll m$;
\end{enumerate}
then there is a finite number of probability
measures $\mu_1,\ldots,\mu_l$ in $M$
with the properties
\begin{enumerate}
\item $\mu_1,\ldots,\mu_l$ are physical absolutely continuous
probability measures;
\item $\supp{\mu_i}\cap\supp{\mu_j}=\emptyset$
for all $1\le i < j \le l$;
\item for all $x\in M$ there are open sets
$V_1=V_1(x),\ldots,V_l=V_l(x)\subset\De$
such that
\begin{enumerate}
\item $V_i \cap V_j = \emptyset,\;\; 1\le i < j \le l$;
\item $\nuinf( \De \setminus (V_1\cup\ldots\cup V_l)) =0$;
\item for all $1\le i\le l$ and
$\nuinf$-a.e. $\subl{t}\in V_i$
we have
\[
\lim_{n\to\infty} \frac1n \sum_{j=0}^{n-1}
\vfi( f^j(x,\subl{t}) ) = \int \vfi \, d\mu_i, \;\;
\text{for every} \; \vfi\in C(M,\real).
\]
\end{enumerate}
Moreover, the sets
$V_1(x),\ldots,V_l(x)$ depend continuously
on $x\in M$ with respect to the distance
$d_{\nu}(A,B)=\nuinf(A\difsim B)$ between
$\nuinf-\bmod 0$ subsets of $\De$.
\end{enumerate}
\end{maintheorem}

The theorem assures the existence of a finite number
of physical probability
measures with respect to the perturbed system
$\SF_{a,\ep}$, as defined in the previous subsections,
which describe the asymptotics of the Birkhoff
averages of almost every perturbed orbit
of every point of $M$.
Section~\ref{Bowen} gives perhaps a clearer
meaning for this result.

\medskip

The conditions on the noise are about
``how much spread'' suffer the orbits
under perturbation when compared with those without
perturbation. They demand that the perturbations
``scatter'' the orbits in an ``uniform'' way around
the nonperturbed ones, at least from some iterates
onward, and ask for negligible perturbations
(of $\nuinf$ measure zero) to produce
negligible effects: the result of such
perturbations should only be a set of
$m$ measure zero.

These hypothesis try to translate the
intuitive idea of random perturbations
not having ``privileged direction or size'',
causing deviations from the ideal orbit
that will ``fill'' a full neighborhood of
that orbit and ``ignoring'' sets of perturbations
of zero probability.
In the light of this, parametric noise satisfying
conditions $A)$ and $B)$ may aptly be referred
to as {\em physical parametric noise}.

\medskip

\begin{example}\label{exemplotoro}
Let $M=\torus^n$ be the $n$-torus, $n\ge1$,
and $f_0:\torus^n\mapto\torus^n$ a
$C^r$-diffeomorphism, $r\ge1$.
Since $\torus^n$ is parallelizable,
$T\torus^n \cong \torus^n \times \real^n$,
we can find $n$ globally orthonormal
(hence nonvanishing) vector fields in
${\cal X}^r(M)$.
For instance, through the identification
$\torus^n \cong \real^n / \relativ^n$
via the natural projection,
we may take
$X_1(x)= e_1 = (1,0,\ldots,0),
X_2(x) = e_2 = (0,1,\ldots,0), \ldots,
X_n(x) = e_n = (0,0,\ldots,1)$
for all $x\in\torus^n$.

We construct a family of differentiable
maps defining
$f:\torus^n \times \real^n \mapto \torus^n$
by
\[
(x,t)\in \torus^n\times\real^n \mapsto
f_0(x) + t_1 X_1(f_0(x)) + \ldots
+ t_n X_n( f_0(x) ) \bmod \relativ^n,
\]
or equivalently by
$
f_t(x) = f(x,t_1,\ldots,t_n) = f_0(x) +
( t_1, \ldots, t_n) \bmod \relativ^n.
$

We note that since $\| t \|_2 < \ep$
implies $\| f_t - f_0 \|_{C^r} < \ep$
for every $\ep>0$ and $\diff^r(\torus^n)$
is open in $C^r(\torus^n,\torus^n)$
(cf. \cite[Capt. I]{PM}), there is $\ep_0>0$
such that the restriction
$
f_{|}:\torus^n\times B^n(0,\ep_0) \mapto \torus^n
$
is a $C^r$-family of $C^r$-diffeomorphisms
of $\torus^n$.

It is not difficult to see that $f$
satisfies hypothesis $A)$ and $B)$ of
theorem \ref{teorema1} for $K=1$ and for
every family
$\SF_{a,\ep} = \left\{
f_t : \| t-a \|_2 < \ep \right\}$
such that 
$\tra{B}^n(a,\ep) \subset B^n(0,\ep_0)$.
We may say, in the light of this, that
this specific kind of random parametric
perturbation is an
{\em absolutely continuous random
perturbation\/}.

Theorem \ref{teorema1} follows and we see
that
{\em any random absolutely continuous
perturbation of a diffeomorphism of
the torus (or of any parallelizable manifold)
is such that Birkhoff
averages exist for almost every
orbit of every point of the torus.
Moreover, their values are defined
by a finite number of absolutely
continuous physical stationary
probability measures\/}.
\end{example}

\begin{remark}
\label{hafamilia}
Example 1 shows that
given any diffeomorphism $f$ of a
parallelizable manifold we may
easily embed $f$ in a suitable
parameterized family of diffeomorphisms
satisfying hypothesis $A)$ and $B)$.
\end{remark}

\begin{example}
\label{constroifamilia}
We now construct an absolutely continuous random
perturbation around any given diffeomorphism $f\in\diff^r(M)$,
$r\ge1$,
of every compact finite dimensional boundaryless
manifold $M$, assuming $M$ to be endowed with
some Riemannian metric.
It is most likely that this kind of construction
can be carried out with $n=\dimns (M)$ or $n+1$
parameters.

\medskip

We start by taking a finite number of coordinate
charts $\{ \psi_i: B(0,3) \mapto M\}_{i=1}^l$
such that $\{ \psi_i( B(0,3)) \}_{i=1}^l$
is an open cover of $M$ and
$\{ \psi_i( B(0,1) ) \}_{i=1}^l$ also
(this is a standard construction,
cf.~\cite[Sec. 1.2]{PM}).
In each of those charts we define
$n=\dimns(M)$ orthonormal vector fields
$\til{X}_{i1},\ldots,\til{X}_{in}: B(0,3) \mapto T_{\psi_i(B(0,3))}M$
and extend them to
the whole of $M$ with the help of bump functions.
This may be done in such a way that the extensions
$X_{ij}$ are null outside $\psi_i(B(0,2))$
and coincide with $\til{X}_{ij}$ in
$\psi_i(\tra{B(0,1)})$,
$i=1,\ldots,l$; $j=1,\ldots,n$.
We then see that
\begin{itemize}
\item At every $x\in M$ there is some $1\le i\le l$
such that $X_{i1}(x),\ldots,X_{in}(x)$
is an orthonormal basis for $T_x M$ --- and likewise
for $X_{i1},\ldots,X_{in}$ --- because
$\{\psi_i(B(0,1))\}_{i=1}^l$ was an open cover of $M$.
\end{itemize}

Finally we define the following parameterized family
$$
F:(\real^n)^l \mapto C^r(M,M), \qquad
F
\left( (u_{ij})_{i=1,\ldots,l \atop j=1,\ldots,n}
\right) (x)
= \Phi\left(
f(x), \sum_{i=1}^l\sum_{j=1}^n u_{ij}\cdot X_{ij} ,1
\right)
$$
where $\Phi:TM\times\real\mapto M$ is the geodesic flow
associated to the given Riemannian metric.
Then for some $\ep_0>0$ we get a finite dimensional
parameterized family of
diffeomorphisms
$F_|:B^{n\cdot l}(0,\ep_0) \mapto \diff^r(M)$
satisfying conditions $A)$ and $B)$
of Theorem~\ref{teorema1} for $K=1$
and some $\xi_0>0$, and for every family
$\SF_{a,\ep}=\{ F_t: \| t-a \|_2 <\ep \}$
where $a\in\tra{B}^{n\cdot l} (a,\ep)\subset
B^{n\cdot l}(0,\ep_0)$.
\end{example}


\begin{example}
\label{ex-FS}
In the context of \emph{random perturbation
of rational functions}, as in~\cite{FS},
hypothesis $A)$ and $B)$ are immediate.

Indeed, let 
$R:\tra{\complex}\times W\mapto \tra{\complex}$
be analytic, where $W\subset\complex$
is open an connected, $z\mapsto R(z,c)$
is rational for all $c\in W$ and
$c\in W \mapsto R(z,c)$ is nonconstant
for every $z\in\tra{\complex}$
(i.e., $R$ is a \emph{generic
family of rational functions}).
Then it is easy to get a $\xi=\xi(c_0,\ep)>0$
such that
$R(z,B(c_0,\ep))\supset B( R(z,c_0), \xi)$
for all $z\in\tra{\complex}$,
whenever $B(c_0,\ep)\subset W$,
by compactness of $\tra{\complex}$
and because analytic nonconstant
functions are open.
Moreover, if $\la$ is Lebesgue measure
normalized and restricted to $B(c_0,\ep)$,
then $R(z,\la)\ll$ Lebesgue on $\complex$.
Hence we get $A)$ and $B)$ with $K=1$.

Theorem~\ref{teorema1} then proves something
more than Theorem 0.1 of~\cite{FS}:
we get physical measures
whose support contains neighborhoods
of the attracting cycles of $R_{c_0}$
and which give the time averages of
almost every orbit of each point of
the Riemann sphere.
\end{example}

\begin{example}\label{exemplotransitivo}
Let $f:M\times T\mapto M$ be a parameterized
family of diffeomorphisms as in section \ref{notacoes}
such that for some $a\in T$ the diffeomorphism
$f_a$ is transitive.
Let us suppose further that for some $\ep>0$
the parametric noise of level $\ep$ around $f_a$,
$\SF_{a,\ep}$, satisfies hypothesis $A)$ and $B)$.
Hence theorem \ref{teorema1} holds and let $\mu_i$
be one of the physical absolutely continuous
probabilities given by the theorem.

Since $f_a$ is transitive, there is a residual set
$\SR$ in $M$ whose points $x_0\in\SR$ give dense
$f_a$-orbits: $\tra{ \{ f_a^k(x_0) \}_{k=0}^\infty }=M$.
Moreover, the c-invariance of $\supp{\mu_i}$
(v. section~\ref{invariantdomains},
definition~\ref{c-inv-set}) and
hypothesis $A)$ imply that
$\inter{\supp{\mu_i}}\neq\emptyset$, and
thus there is
$x_0\in \left( \SR \cap \inter{\supp{\mu_i}} \right)$.

We deduce that
$\supp{\mu_i} \supset \tra{ \{ f^k(x_0,\De) \}_{k=1}^\infty }
\supset \tra{ \{ f_a^k(x_0) \}_{k=0}^\infty }=M$
and so there is only one physical absolutely continuous
probability in $M$, whose support is the whole of $M$.

In particular, {\em every diffeomorphism of
the torus $\torus^n\,(n\ge1)$ with a dense
orbit, under absolutely continuous noise of
arbitrary level $\ep>0$, has a single
physical absolutely continuous probability
whose support is $M$
(and likewise if $M$ is any parallelizable
compact boundaryless manifold).}

\end{example}

\medskip

In section \ref{bifurcations}
we shall see that certain
arcs (uniparametric families) of
diffeomorphisms of class $C^r$ ($r\ge3$)
generically unfolding a quadratic
homoclinic tangency satisfy
both conditions of theorem \ref{teorema1},
restricted to a neighborhood of the point
of homoclinic tangency.
For more specifics, check the abovementioned
section.
We will then have

\medskip

\noindent
{\bf Theorem 2} {\em
There are open sets of arcs (in the $C^3$ topology)
$\{ f_t \}_{t\in]-1,1[}$ of diffeomorphisms of
class $C^3$ of a compact boundaryless surface
generically unfolding
a quadratic homoclinic tangency at $f_0$ such that,
in a neighborhood $\SQ$ of a point of
homoclinic tangency and for all $f_{t_0}$
sufficiently near $f_0$ under parametric
noise of sufficiently small level $0<\ep<\ep_0$,
there are a finite number of probability measures
$\mu_1,\ldots,\mu_h$ in $\SQ$ that satisfy
the conditions $1)$ and $2)$ and also $3)$
of theorem \ref{teorema1},
for points $x\in M$ whose orbits
$\SO(x,\subl{t})$ have an infinite number
of iterates in $\SQ$ with respect to a $\nuinf$
positive measure set of perturbations.
}

\medskip

This result, combined with Newhouse's phenomenon,
shows that the infinity of periodic hyperbolic
attractors (sinks) that coexist in a neighborhood
of a point of homoclinic tangency, for ``many''
parameter values near the bifurcation parameter,
cannot ``survive'' the random parametric perturbation.
Moreover it must subsist, at most, a finite number
of {\em analytic continuations under random
perturbation\/} of a sink. Section \ref{sinks}
will specify this conclusions and extend the
result in a simple manner to Colli's phenomenon,
where the infinity of hyperbolic periodic attractors
is replaced by an infinity of H\'enon-like strange
attractors.

Now we will concentrate on the proof of
theorem \ref{teorema1}.

\section{Invariant Domains}
\label{invariantdomains}

Let $\mu$ be a stationary probability measure with respect to
a parametric perturbation of noise level $\ep>0$
around $f_a$. Then $\supp{\mu}$ is
$S$-invariant:
$S( \supp ( \mu\times\nuinf ) )
\subset \supp ( \mu\times\nuinf )$.

Let us observe that since
$\supp(\mu)=\supp(\mu)\times\De$
we have for all
$(x,\subl{t})\in\supp(\mu)\times\De$
that $f^k(x,\subl{t})\in \supp\mu$,
for all $k\ge1$.
That is, $\supp\mu$ is completely
invariant  according to


\begin{definition}\label{c-inv-set}
A part $C$ of $M$ is said
{\em completely invariant\/} or {\em c-invariant\/}
if
$f^k(x,\subl{t})\in C$
for all $x\in C$, $\subl{t}\in\De$
and $k\ge1$.
\end{definition}

With the purpose of showing the existence
of the kind of stationary probability measures
stated in theorem \ref{teorema1} and to
better understand the dynamics of the points
in their support as well, we make a series
of definitions.

\begin{definition}\label{defdominv}
An {\em invariant domain under an $\ep$-perturbation
with respect to the family $f$
around the parameter $a\in I$\/}
will be a finite collection $\dominvu$
of pairwise {\em separated open sets\/}, that is,
$i\neq j\impl \tra{\SU_i}\cap\tra{\SU_j} = \emptyset$,
such that
$f^k (\SU_0,\De)\subseteq \SU_{k \bmod r}$
for all $k\ge1$,
and it will be written $D=(\SU_0,\ldots,\SU_{r-1})$.
The number $r\in\natur$
above will be referred to as the {\em period\/}
of the invariant domain.
\end{definition}

Let us observe that the open set $\SU_0$ has a
privileged role in the above definitions.

\begin{definition}\label{defdomsimetr}
An invariant domain that also satisfies
\begin{equation}\label{s-invariant}
f^k (\SU_i,\De) \subseteq \SU_{ (k+i) \bmod r }, \;\;\forall k\ge1
\end{equation}
whatever $i\in\{ 0,\dots,r-1 \}$ will be a
{\em symmetrically invariant\/} domain or
{\em s-invariant\/} domain.
\end{definition}

This kind of domains will be at the heart of the
arguments within next sections
and the proof of their existence and finite
number is the key to every other result in this paper.

\begin{remark}\label{closure-s-inv}
Since the $f_t$ are diffeomorphisms for all $t\in T$,
we see that if the collection
$D=(\SU_0,\ldots,\SU_{r-1})$ is
s-invariant, then
$\tra{\SD}=(\tra{\SU}_0,\ldots,\tra{\SU}_{r-1})$
also satisfies \refer{s-invariant} and conversely:
if the closure
$\tra{\SD}=(\tra{\SU}_0,\ldots,\tra{\SU}_{r-1})$
satisfies \refer{s-invariant} with
$\SU_0,\ldots,\SU_{r-1}$ pairwise disjoint open
sets, then $D=(\SU_0,\ldots,\SU_{r-1})$ is
an s-invariant domain.
\end{remark}

\subsection{Partial Order and Minimality}
\label{defordparc}

Let $\SD$ be the family of s-invariant domains.
We define the following partial order relation between
its elements.

Let $D=\domcomplinv$ and $D'=\domcomplinvv$ be elements
of $\SD$.

\medskip

First, $D=D'$ if there are $i,i'\in\natur$ such that
$\SU_{(i+k) \bmod r} = \SU_{(i'+k) \bmod r'}', \;
\forall k\ge1$
which implies $r=r'$, because the open sets that form
each invariant domain are pairwise disjoint.

\begin{figure}[h]
 
\centerline{\psfig{figure=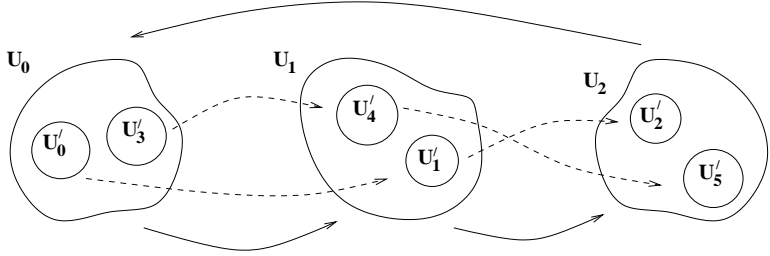,width=4.0in}}
 
\caption{\label{ordenados} Domains $D$,$D'$ with $D'\prec D$} 
\end{figure}

We say $D\prec D'$ if there are $i,i'\in\natur$
such that
$\SU_{i\bmod r} \subseteq \SU_{i'\bmod r'}' $
but
$\SU_{i\bmod r} \neq \SU_{i'\bmod r'}' $,
and
$\SU_{(i+k) \bmod r} \subseteq \SU_{(i'+k) \bmod r'}' $
for all $k\ge1$
(see figure \ref{ordenados} for an example with
$r=3$ and $r'=6$).

\medskip

We write $D\preceq D'$ if, and only if, $D=D'$ or $D\prec D'$.

Clearly $(\SD,\preceq)$ is now a partially ordered set.

\begin{definition}\label{defminimal}
A {\em\bf minimal invariant domain\/} is a domain $D\in\SD$
which is minimal with respect to the partial order $\preceq$
just defined.
\end{definition}

{\em Minimal domains will be represented by the letter $\SM$
throughout this text.}

\section{A Tour of the Proof}
\label{esquema}

With the notions given in previous sections we
can now divide the proof of theorem \ref{teorema1}
in the following steps:

\begin{nliste}{\Parentalgarismo}
\item To show that $\SD$ has some minimal invariant domain
and that any invariant domain contains some minimal one
(section \ref{minimalexist}).
\item To show that minimal invariant domains are
pairwise disjoint (section \ref{minimais-disjuntos}).
\end{nliste}

By now we can already deduce the number of minimals is finite.
In fact, a minimal invariant domain $\SM$
is completely invariant
and by hypothesis $A)$ of theorem \ref{teorema1} we see
that every open set of the finite collection forming $\SM$
contains a ball of radius $\ge\xi_0>0$.
The compactness of $M$ and step $2$ above ensure
there can only be a finite number of such open sets and
thus a finite number of minimals also.

\begin{nliste}{\Parentalgarismo}
\item[(3)] Every minimal domain is
{\em randomly transitive\/} or {\em r-transitive\/},
this notion will be specified in
subsection \ref{transitivos}.
\item[(4)] The orbits of every point $z\in\SM$ under
noise generate a stationary probability measure
$\mu$ which is absolutely continuous
(section \ref{abscont}).
\end{nliste}

From $3$ and $4$ we deduce that
there exists an absolutely continuous stationary
probability $\mu$ in the closure of each minimal $\SM$
(since $\SM$ contains every orbit of $z\in\SM$) whose
support is the closure of $\SM$
(by the c-invariance of the support and item 3):
$\supp{\mu}=\tra{\SM}$.

\begin{nliste}{\Parentalgarismo}
\item[(5)] Every stationary absolutely continuous
probability measure $\mu$ supported on
a minimal domain $\SM$
is ergodic and its ergodic basin
$E(\mu)$ contains the whole of $\SM$:
$E(\mu)\supset\tra{\SM}$
(section \ref{ergodicidade}).
\end{nliste}

Being ergodic, absolutely continuous and
supported on the whole of $\SM$, this probability
$\mu$ is physical, since the minimal invariant domain
is a collection of open sets.
Consequently since for every such measure
$E(\mu)\supset\tra{\SM}$ holds,
this is the only
stationary ergodic absolutely continuous probability
measure supported on $\SM$.
It will be referred to as the
{\em characteristic probability\/} of the minimal
$\SM$.

\begin{nliste}{\Parentalgarismo}
\item[(6)] Every stationary probability measure is
supported on some s-invariant domain
(section \ref{decomposition}).
\end{nliste}

This crucial step gives the converse of
the property deduced from step $5$.
Moreover, combining with the results
of the previous steps we will
deduce from step $6$ that

\begin{nliste}{\Parentalgarismo}
\item[(7)] Every stationary probability measure
is a finite convex linear combination of
characteristic probabilities
(section \ref{decomposition}).

\item[(8)] Finally, in section \ref{conclusion},
we will use items 4 and 7 to deduce that
$\nuinf$-a.e. perturbation $\subl{t}\in\De$ is
such that $\SO(z,\subl{t})$ eventually falls
into some minimal $\SM$. The perturbations sending
$z$ into different minimals form the partition
of item 3 of Theorem \ref{teorema1}.
Since $\SM$ supports
a characteristic measure which is physical,
we further derive that Birkhoff averages exist
for $\SO(z,\subl{t})$ and satisfy~\refer{mediaorbital}.
\end{nliste}

\section{Fundamental lemmas}
\label{fundamental}

The measure theoretical lemma that follows will be used frequently
within the arguments of this and next sections.

\begin{lemma}\label{lema0-1}
Given $V\subset\De$ with $\nu^\infty (V)>0$, we define for fixed
$\subl{\te}\in\De$ and $k\ge 1$
\[
V({\subl{\te}},k) = \{ \subl{\om}\in V: \om_1=\te_1,\ldots,\om_k=\te_k \}
\]
the $k$-section of $\;V$ along $\subl{\te}$. Then we have
\[
\nu^\infty ( \si^k V({\subl{\te}},k) ) \mapto 1 \;\;
	\text{ when } \;\; k\mapto\infty
\]
for $\nu^\infty$-a.e. $\subl{\te}\in V$,
where $\si:\De\mapto\De$ is the left shift on sequences:
$\si(\subl{\psi}) = \subl{\vfi}$
with $\vfi_n=\psi_{n+1}, \; n=1,2,3\ldots$
\end{lemma}

{\bf Note:} From now on we will say that a vector $\subl{\te}$
satisfying the above limit with respect to a set $V\subset\De$
is {\em $V$-generic\/}.

\begin{figure}[h]
                 
\centerline{\psfig{figure=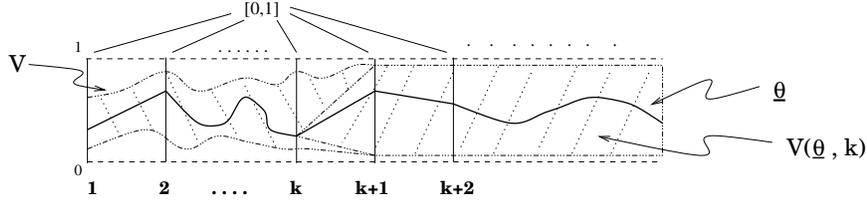,width=4.5in}}
 
\caption{{\em
\label{seccoes} Representation of the infinite product
of the interval $[0,1]$, a vector $\subl{\te}$ and the sets
$V$ and $V(\subl{\te},k)$ }}
\end{figure}

\begin{proof}
We may assume, for definiteness,
that $\De=[0,1]^{\natur}$ with
$\nu$ the Lebesgue measure in $[0,1]$ so that
$\nu^\infty$ is a probability in $\De$. Let
$V\subset\De$ be such that $\nu^\infty ( V ) > 0$.

If $\SB$ is the Borel $\si$-algebra in $[0,1]$ and
$\SB_k=
\overbrace{\SB\times\ldots\times\SB}^k
\times[0,1]^\natur, \; k \ge 1$, then
$\SA=\si( \cup_{k=1}^\infty \SB_k )$ is the $\si$-algebra
of $\De$ over which $\nu^\infty$ is defined, the $\si$-algebra
generated by all $\SB_k$.
For every $f\in L^1(\De,\SA,\nu^\infty)$ and each $k\ge1$ the
map $A\in\SB_k\mapsto\int_A f \, d\nu^\infty$
defines a finite measure on $(\De,\SB_k,\nu^\infty)$, which clearly
is absolutely continuous with respect to the measure
$A\in\SB_k\mapsto\nu^\infty(A)$ (the restriction
of $\nu^\infty$ to $\SB_k$).
By Radon-Nikodym's theorem there is
$E(f|\SB_k)\in L^1(\De,\SB_k,\nu^\infty)$, the
{\em conditional expectation of $f$ with respect to the
$\si$-algebra $\SB_k$\/}, such that
\begin{equation}\label{espcondk}
\int_A E(f|\SB_k) \, d\nu^\infty = \int_A f \, d\nu^\infty, \;
\forall A\in\SB_k
\end{equation}
and this function is unique with this property in
$L^1(\De,\SB_k,\nu^\infty)$.

\smallskip

Let $X_k=E(f|\SB_k), \; k=1,2,\ldots$ We are going to see that
$\{ X_k \}_{k=1}^\infty$ is a martingale with respect to the
sequence $\{ \SB_k  \}_{k=1}^\infty$ of $\si$-algebras.

Indeed, because $\SB_k\subset\SB_{k+1}$ we have
$
\int_A E(f|\SB_{k+1}) \, d\nu^\infty = \int_A f \, d\nu^\infty
$
for all $A\in\SB_k$
and by \refer{espcondk} and uniqueness of conditional expectation
\[
E(X_{k+1}|\SB_k) = E( E(f|\SB_{k+1}) |\SB_k  )=E(f|\SB_k)=X_k \;\;
\nu^\infty-\text{a.e.}
\]

By the martingale convergence theorem (cf. \cite{P} for
simple definitions and proofs), the sequence
$\{X_k\}_{k=1}^\infty$ has a $\nu^\infty$-a.e. limit that
we shall write $X\in L^0(\De,\SA)$.

By \refer{espcondk} and because
$f\in L^1(\De,\SA,\nu^\infty)$ we have, assuming
$f\ge 0$, that $X_k\ge 0 \;\; \nu^\infty$-a.e., $k\ge 1$, and
consequently $X\ge 0 \;\; \nu^\infty$-a.e.. Moreover
$
\int | X_k | \, d\nu^\infty = \int X_k \, d\nu^\infty =
	\int f \, d\nu^\infty = \int | f | \, d\nu^\infty
$
and so $X\in L^1(\De,\SA,\nu^\infty)$
by dominated convergence and 
$\int | X | \, d\nu^\infty = \int X \, d\nu^\infty$
gives
$\int f \, d\nu^\infty$.
Furthermore, if $A\in\SB_k$ then
$
\int_A X_j \, d\nu^\infty = \int_A f \, d\nu^\infty
$
for all $j\ge k$
and from this we get
$\int_A X\, d\nu^\infty = \int_A f \, d\nu^\infty$
for all $A\in\SB_k$ and $k\ge 1$.

\medskip

By the absolute continuity of the integral of a $L^1$-function and
by definition of $\SA$, for every
$\ep>0$ and $A\in\SA$ there are $\de>0,\, k\ge 1$ and $B\in\SB_k$
such that
$
\nuinf (A\difsim B)<\de
$,
$
\int_{A\difsim B} |X| \, d\nuinf < \ep
$
and
$
\int_{A\difsim B} |f| \, d\nuinf < \ep.
$
Now we have, in succession
\begin{eqnarray*}
\left| \int_A f \, d\nuinf - \int_B X_k \, d\nuinf \right| &=&
\left| \int_A f \, d\nuinf - \int_B f \, d\nuinf \right| \le
\int_{A\difsim B} |f| \, d\nuinf < \ep ;
\\
\left|\int_A X \, d\nuinf - \int_B X_k \, d\nuinf \right| &=&
\left| \int_A X \, d\nuinf - \int_B X \, d\nuinf \right| \le
\int_{A\difsim B} |X| \, d\nuinf < \ep ;
\end{eqnarray*}
and from this we get
$
\left| \int_A X \, d\nuinf - \int_A f \, d\nuinf \right| \le 2\ep
$
with $\ep>0$ arbitrary.

We conclude that $\int_A X \, d\nuinf = \int_A f \, d\nuinf,\;
\forall A\in\SA\,$ and so $X=f\;\nuinf$-a.e..

\medskip

In particular if $f=1_V$ we have
$X_k \mapto 1_V \; \nuinf$-a.e.
and
$
\int_B E(1_V|\SB_k) \, d\nuinf =
\int_B 1_V \, d\nuinf
$
equals $\nuinf (V\cap B)$
by definition of conditional expectation.
But
$\nuinf (V\cap B)= \int_B 1_V \, d\nuinf$
also equals
$
\int_B \nu_k^\infty ( V(\subl{\te},k) ) \, d\nu^k ( \subl{\te} )
$
for every $B\in\SB_k$ and $k\ge 1$
by Fubini's theorem, where
$\nu_k^\infty( A ) = \nuinf( \si^k A )$ and
$\nu^k ( A ) = \nu^k ( \pi_k ( A ) ), \, A\in\SA$ with
$\pi_k:\De\mapto [0,1]^k$ the natural projection
$\subl{\te}=( \te_i )_{i=1}^\infty
\mapsto \vecto{\te}{k}$.
That is
$
\nu_k^\infty ( V(\subl{\te},k) ) = E(1_V|\SB_k) = X_k
$
$\nuinf$-a.e.  $\subl{t}\in\De$,
and the proof is complete.

\end{proof}

This lemma will be utilized essentially in the
following way. Let $V,W$ be subsets of $\De$
with $\nuinf$-positive measure and $\subl{t}$
a $V$-generic vector. Then there is $k_0\in\natur$
such that
\begin{equation}\label{usolema0-1}
k\ge k_0 \impl
\nuinf \left( W \cap \si^k V(\subl{t},k) \right) >0.
\end{equation}
Since $\nuinf( V(\subl{t},k) )=0$ for all
$\subl{t}\in\De$ and $k\ge1$ we may wonder
whether we may use \refer{usolema0-1} in arguments
proving some $\nuinf$-a.e. result.
The answer is in the following

\begin{lemma}\label{lema0-1>0}
Let $V,W \subset \De$ be such that
$\nuinf(V),\nuinf(W) > 0$.
Then for $\nuinf$-a.e. $\subl{t}\in V$
there is a $k_0\in\natur$ such that
for all $k\ge k_0$ and every $\eta>0$
\[
\nuinf \left\{
\subl{s}\in V : d(\subl{s},\subl{t}) < \eta
\;\;\text{and}\;\;
\si^k\subl{s} \in W
\right\} >0.
\]
\end{lemma}

Hence we may not have \refer{usolema0-1} but we
know we can choose with positive probability
a vector in $V$ arbitrarily close to $\subl{t}$
whose $k$th shift is in $W$.
This will be enough for our purposes.

\begin{proof}
Let $V\subset\De$ be such that $\nuinf(V)>0$.
For every $n\ge1$ and $j\ge1$ let $K_{n,j}$
be a compact set inside $V$ such that
$\nuinf(V\setminus K_{n,j}) < (n\cdot 2^j)^{-1}$
and $E(1_V | \SB_j )_{|K_{n,j}}$ is continuous
-- we are using Luzin's theorem
(v. \cite[Capt. IV, Sec. 21]{M}).
Then $C_n=\cap_{j\ge1} K_{n,j}$ is
a compact subset of $V$,
$\nuinf(V\setminus C_n) \le n^{-1}$
and $E(1_V | \SB_j )_{|C_n}$ is continuous
for every $j,n\ge1$.

We have $V=\cup_{n\ge1} C_n$, $\nuinf\bmod 0$
and so $\nuinf$-a.e. $\subl{t}\in V$ is in some
$C_n$, $n\ge1$.
Moreover $\nuinf$-a.e. $\subl{t}\in V$ is
a $\nuinf$-density point of some $C_n$ and we
may suppose $\nuinf(C_n)>0$ for all $n\ge1$
(otherwise we consider only $n\ge n_0$ for some
big $n_o\in\natur$).

\medskip

>From now on we suppose $\subl{t}$ is $V$-generic
and a $\nuinf$-density point of some $C_n$ with
$\nuinf( C_n)>0$. We let $W\subset\De$ be such
that $\nuinf(W)>0$, set
$\de=\frac14\nuinf(W) >0$ and let $k_0\in\natur$
be such that
$
\nuinf( \si^k V(\subl{t},k) ) \ge 1-\de,
$
for every $k\ge k_0$
by lemma \ref{lema0-1}.
By the choice of $\subl{t}$ and $C_n$ we have
$\nuinf( B(\subl{t},\eta) \cap C_n) >0$
for all $\eta>0$ and for some $\eta_0>0$
we have further that, fixing $k\ge k_0$,
\[
d(\subl{s},\subl{t}) < \eta_0, \subl{s}\in C_n
\impl
\nuinf( \si^k V(\subl{s},k) ) \ge 1-2\de
\]
by the continuity of $E(1_V|\SB_k)_{|C_n}$
at $\subl{t}$.
Therefore we deduce that
\[
d(\subl{s},\subl{t}) < \eta_0, \subl{s}\in C_n
\impl \nuinf \left( W \cap \si^k V(\subl{s},k)
\right) \ge 2\de >0
\]
and so, for any $\eta>0$, we have
\begin{eqnarray*}
&\nuinf&
\left\{ \subl{s}\in V :
d(\subl{s},\subl{t}) < \eta \;\;\text{and}\;\;
\si^k\subl{s}\in W \right\} \ge
\\
&\ge&
\nuinf \left\{ \subl{s}\in V :
d(\subl{s},\subl{t}) < \eta_1=\min\{\eta_0,\eta\}
\;\;\text{and}\;\; \si^k\subl{s}\in W
\right\} \ge
\\
&\ge&
\int_{B(\subl{t},\eta_1)\cap C_n}
1_W( \si^k\subl{s} ) \, d\nuinf(\subl{s})
= \int_{B(\subl{t},\eta_1)\cap C_n}
\int_{\si^k V(\subl{s},k)}
1_W ( \subl{u} ) \, d\nuinf(\subl{u})
\, d\nu^k(\subl{s})
\\
&=&
\int_{B(\subl{t},\eta_1)\cap C_n}
\nuinf\left( W \cap \si^k V(\subl{s},k) \right)
\, d\nu^k(\subl{s})
\ge 2\de \cdot \nu^k( B(\subl{t},\eta_1)\cap C_n)
\\
&\ge&
2\de \cdot \nuinf( B(\subl{t},\eta_1)\cap C_n ) >0
\end{eqnarray*}
where we have used Fubini's theorem and
$\nu^k$ is as before in lemma \ref{lema0-1}.
\end{proof}

\medskip

In section~\ref{regular:prova}
a slight generalization of lemma \ref{lema0-1}
will be needed.

\begin{definition}\label{seccaodupla}
Given $V\subset\De$ and $\subl{t},\subl{s}\in\De$
we define a {\em double section\/} through
$\subl{t}$ and $\subl{s}$ at $k\ge1$ by
$
V(\subl{t},k,\subl{s})=
\left\{ \subl{\vfi}\in V :
\vfi_1=t_1,\ldots,\vfi_k=t_k \;\;\text{and}\;\;
\vfi_{k+2}=s_1,\vfi_{k+3}=s_2,\ldots
\right\} .
$
\end{definition}

%
%

\begin{lemma}\label{lema0-1-2}
Let $V\subset\De$ be such that $\nuinf(V)>0$.
Then for $\nuinf$-a.e. $\subl{t}\in V$ and
for every $0<\ga,\de<1$ there exists $k_0\in\natur$
such that for all $k\ge k_0$ there is a set
$W_k\subset V$ with the properties
\begin{enumerate}
\item $\subl{t}\in W_k$;
\item $\nuinf( W_k ) >0$;
\item $\nu\left(
p_{k+1} W_k( \subl{t},k,\subl{s}) \right)
\ge 1-\de$ for $\nuinf$-a.e. $\subl{t}\in W_k$
and $\subl{s}$ in a subset of $\De$ with
$\nuinf$-measure $\ge 1-\ga$;
\end{enumerate}
where $p_k:\De\mapto B$ is the projection on the
$k$th coordinate.
\end{lemma}

\begin{proof}
(An application of lemma \ref{lema0-1} and
Fubini's theorem)

Defining
$
V_n=\left\{ \subl{t}\in V:
\nuinf( \si^k V(\subl{t},k) ) \ge 1-\de \cdot (1-\ga), \;
\forall k\ge n \right\}
$
we have $V_n\subset V_{n+1}$ and lemma
\ref{lema0-1} says $V=\cup_{n\ge1} V_n,\;\nuinf\bmod0$.
We set $k_0\in\natur$ such that
$\nuinf(V_k) \ge \frac45\nuinf(V)$ for every $k\ge k_0$.
Definition \ref{seccaodupla} and Fubini's theorem
imply
\[
1-\de \cdot (1-\ga) \le \nuinf( \si^k V(\subl{t},k) )
=\int \nu [ p_{k+1} V(\subl{t},k,\subl{s}) ]
\, d\nuinf(\subl{s})
\]
for every $\subl{t}\in V_{k_0}$ and $k\ge k_0$.
We define now for each $\subl{t}\in V_{k_0}$ and
$k\ge k_0$ the set
\[
W_k(\subl{t}) = \left\{
\subl{s}\in\De:
\nu[ p_{k+1} V(\subl{t},k,\subl{s})] \ge 1-\de
\right\}
\]
and by the last inequality we see that
$\nuinf( W_k(\subl{t}) ) \ge 1-\ga$.
Then defining for $k\ge k_0$
\[
W_k= \bigcup \left\{
V(\subl{t},k,\subl{s}) :
\subl{t}\in V_{k_0} \;\;\text{and}\;\;
\subl{s}\in W_k(\subl{t})
\right\},
\]
we get
\begin{eqnarray*}
\nuinf( W_k ) 
&=&
\int_{V_{k_0}}
\int_{W_k(\subl{t})}
\nu[ p_{k+1} V(\subl{t},k,\subl{s})]
\, d\nuinf(\subl{s}) \, d\nu^k(\subl{t})
\\
&\ge&
\int_{V_{k_0}} (1-\de)\cdot
\nuinf ( W_k(\subl{t}) ) \, d\nu^k(\subl{t})
\ge
(1-\de)\cdot (1-\ga) \cdot \nu^k(V_{k_0})
\\
&\ge& (1-\de)(1-\ga) \cdot \nuinf(V_{k_0})
\ge
(1-\de)(1-\ga) \cdot \frac45 \cdot \nuinf(V) >0.
\end{eqnarray*}

We finally note that $\nuinf$-a.e. $\subl{t}\in V$
is in every $V_{k}$ for sufficiently big $k$.
\end{proof}

\medskip

The following notions will be extremely useful.
They are mere adaptations of the usual notions
of $\om$-limit to the context of random
parametric perturbations.

\begin{definition}\label{defswlimite}
We take $z$ to be some point in $M$, $U$ some subset of $M$,
$\subl{t}$ some vector in $\De$ and define
\begin{eqnarray*}
	\om( z , \subl{t} ) &=& \left\{ w \in M :
	\exists n_1<n_2<\ldots \;\text{in}\; \natur \;
	\text{such that}\;
	f^{n_j}_{\subl{t}} (z) \mapto w \;\text{when}\;
	j\mapto\infty
	\right\}
	\\
	&&( \text{the usual definition of} \;\;
	\om \text{-limit for the orbit} \;
	\SO(z,\subl{t}) );
	\\
	\om( U, \subl{t} ) &=& \left\{ w\in M :
	\exists \{ u_j \}_{j=1}^\infty \subset M \,
	\exists n_1<n_2<\ldots \;\text{in}\; \natur \;
	\text{such that} \;
	f^{n_j}_{\subl{t}} (u_j) \mapto w
	\right.
	\\
	&& \left. \text{when} \; j\mapto\infty
	\right\} \;
	( \text{the} \; \om \text{-limit of a set under a
	perturbation vector} \; \subl{t} );
	\\
	\om( z, \De ) &=& \left\{ w\in M :
	\exists \{ \subl{\te}^{(j)} \}_{j=1}^\infty \subset \De \,
	\exists n_1<n_2<\ldots \;\text{in}\; \natur \;
        \text{such that} \;
	f_{ \subl{\te}^{(j)} }^{n_j} (z) \mapto w
	\right.
	\\
	&& \left. \text{when} \; j\mapto\infty
	\right\} \;
	( \text{the}\; \om \text{-limit of a point
	under every perturbation} );
	\\
	\om( U, \De ) &=& \left\{ w\in M :
	\exists \{ u_j \}_{j=1}^\infty \subset M \;
	\exists \{ \subl{\te}^{(j)} \}_{j=1}^\infty \subset \De \,
	\exists n_1<n_2<\ldots \;\text{in}\; \natur \;
        \text{such that} \;
	\right.
	\\
	&& \left. f_{ \subl{\te}^{(j)} }^{n_j} (u_j) \mapto w \;
	\text{when} \; j\mapto\infty
	\right\} \;
	( \text{the same as before with respect to a set} ).
\end{eqnarray*}

\end{definition}

\begin{lemma}\label{wlimiteV}
Let us suppose $U$ to be a subset of $M$ whose orbits,
under a positive $\nuinf$-measure set
$V\subset\De$ of perturbations,
go through a finite family
of pairwise separated open sets
$A_0,\ldots,A_{l-1}$ in a cyclic way, that is
\begin{equation}\label{wlV-0}
f_V^k ( U ) \subset A_{ k\bmod l },\, \forall k\ge 1
\end{equation}
(example: the set $\SU_0$ of an invariant domain
$D\in\SD$ with respect to $\SU_0,\ldots,\SU_{r-1}$)

Then the set $\om( U , \subl{\te} )$ of accumulation points
of the orbit of $\;U$ under a $V$-generic perturbation
$\subl{\te}\in V$ is such that
$\om( U,\subl{\te} ) \subset \tra{A}_0 \cup \ldots \cup \tra{A}_{l-1}$
and if $z\in\om( U,\subl{\te} )\cap \tra{A}_i$ with $0\le i \le l-1$
and $\subl{\psi}\in\De$, then
$f_{\subl{\psi}}^k(z)\in\tra{A}_{(i+k)\bmod l}$
for all $k\ge 1$.
\end{lemma}

\begin{lemma}\label{wlimiteDelta}
If in the last lemma we had $V=\De$ then the set
$\om( U,\De )$, besides having orbits that go
in a cyclic way through the
$\tra{A}_i,\, i=0,\ldots,l-1$, under {\em any\/} perturbation,
would also be invariant under {\em every\/} perturbation:
$f_{\subl{\psi}}^k (z) \in \om ( U,\De )$
for all $k \ge 1$,
for all $z\in\om( U,\De )$
and for all $\subl{\psi}\in\De$.
\end{lemma}

These lemmas essentially state that whenever we look at
limits of {\em generic\/} perturbations we find a point
whose perturbed orbit does {\em not\/} depend on the
perturbation chosen, in the sense that it is carried cyclically
through some specified family of sets. This property
is the key idea behind the construction of
s-invariant domains in lemma \ref{construir}

\begin{rproof}{\ref{wlimiteV}}
Let us fix $\subl{\psi}\in\De$ and $z\in\om(U,\subl{\te})$ with
$\nuinf(\si^j V(\subl{\te},j) ) \mapto 1$ when $j\mapto\infty$.

Then there are sequences
$\{ u_j \}_{j=1}^\infty \subset U$ and
$\{ n_j \}_{j=1}^\infty\subset\natur$ with
$n_1<n_2<\ldots$ such that
$
z_j = f^{n_j} ( u_j, \subl{\te} ) \mapto z
$
when $j\mapto\infty$.
It is clear that $z\in \tra{A}_0 \cup \ldots \cup \tra{A}_{l-1}$.

Let us now fix $k\in\natur$ and assume $z\in\tra{A_i}$
for some $i\in\{0,\ldots,l-1\}$. We want to show that
$
f^k (z, \subl{\psi}) \in \tra{A}_{(i+k)\bmod l}.
$

Once $k$ is fixed,
property \ref{prop0}
implies that, for given
$\de>0$, there are  $\ga,\ups>0$ such that
\begin{eqnarray}
d(\subl{\psi},\subl{\vfi}) < \ga,\, \subl{\vfi}\in\De
&\impl&
d_M(f^k (z,\subl{\psi}), f^k (z,\subl{\vfi}) ) < \de \nonumber ;
\\
d_M(z_1,z_2)<\ups,\, z_1,z_2\in M,\, \subl{\vfi}\in\De
&\impl&
d_M( f^k (z_1,\subl{\vfi}), f^k (z_2,\subl{\vfi}) ) < \de.
\label{wlV-3}
\end{eqnarray}

By lemma \ref{lema0-1>0} and the convergence of
$\{ z_j \}_{j=1}^\infty$,
making $W=B(\subl{\psi},\ga/2)$ we may
choose a sufficiently big $j\in\natur$ such that
$
d_M( f^{n_j} (u_j, \subl{\te}) , z ) < \ups/2
$
and a sufficiently small $\eta>0$ such that,
with positive probability, there can be found
$\subl{\vfi}\in V$ with
\begin{equation}\label{wlV-aprox}
d(\subl{\vfi}, \subl{\te} ) < \eta, \;
d( \si^{n_j} \subl{\vfi}, \subl{\psi} ) < \ga/2
\;\; \text{and also} \;\;
d_M( f^{n_j}( u_j, \subl{\vfi}), z) < \ups.
\end{equation}

Hence,
by the choice of $\ga$ and $\ups$ we will have that:
\begin{eqnarray*}
d_M( f^k (z,\subl{\psi}),f^k (z_j,\si^{n_j} \subl{\vfi}) )
& \le &
d_M( f^k (z,\subl{\psi}), f^k (z,\si^{n_j} \subl{\vfi}) ) +
\\
& &
+ d_M( f^k (z,\si^{n_j} \subl{\vfi}),
	f^k (z_j,\si^{n_j} \subl{\vfi}) )
\\
& \le &
\de + \de = 2\de.
\end{eqnarray*}

But we can take $\ups>0$ so small that
besides \refer{wlV-aprox} and we get
\begin{equation}\label{wlV-6}
d_M(w,z)<\ups,\, w\in A_0\cup\ldots A_{l-1} \impl z \in A_i.
\end{equation}
With this we have $z_j\in A_i$ and also
$f^k (z_j,\si^{n_j} \subl{\vfi})\in A_{(i+k)\bmod l}$
by the hypothesis \refer{wlV-0}, with $\de>0$ arbitrary,
and the lemma follows immediately.
\end{rproof}

\begin{rproof}{\ref{wlimiteDelta}}
Let us take $z\in\om( U,\De)$ and suppose $z\in\tra{A_i}$ for
some $i\in\{0,\ldots,l-1\}$.
We fix $k\ge1$ and $\subl{\psi}\in\De$.

Then there are $\{ \subl{\te}^{(j)} \}_{j=1}^\infty \subset \De$,
$\{ u_j \}_{j=1}^\infty \subset U$ and
$\{ n_j \}_{j=1}^\infty \subset \natur$ with
$n_1< n_2 < \ldots$ in such a way that
$
z_j = f^{n_j} ( u_j,\subl{\te}^{(j)} ) \mapto z
$
when
$j \mapto \infty$.

For $\de>0$ let us take $\ups>0$ as in \refer{wlV-3},
and $\ups$ so small  that \refer{wlV-6} holds.
Moreover, let $j_0\in\natur$ be such that
$
j\ge j_0 \impl
d_M( f^{n_j}( u_j, \subl{\te}^{(j)} ) , z ) < \ups.
$

We now have
$d_M( f^k (z,\subl{\psi}) , f^k (z_j,\subl{\psi}) ) =
d_M( f^k (z,\subl{\psi}) , f^{k+n_j} (u_j,\subl{\til{\te}}^{(j)} ) )
\le \de$
for $j\ge j_0$,
where
$\subl{\til{\te}}^{(j)}=
(\te_1^{(j)},\ldots,\te_{n_j}^{(j)},\psi_1,\ldots,
\psi_k,\psi_{k+1},\ldots) \in \De$.
But $\de>0$ is arbitrary, thus we get that
$
f^{k+n_j} ( u_j, \subl{\til{\te}}^{(j)} )
\mapto f^k_{\subl{\psi}} (z)
$
when
$j \mapto \infty$.

Now we see that for $f^k (z,\subl{\psi})$ there exist
$\{ \subl{\til{\te}}^{(j)} \}_{j=1}^\infty\subset\De$,
$\{ u_j \}_{j=1}^\infty \subset U$ and
$\{k+n_j\}_{j+1}^\infty \subset\natur$ with
$k+n_1 < k+n_2 < \ldots$ in such a way that
$f^k ( z, \subl{\psi} ) \in\om( U,\De )$.
\end{rproof}

\medskip

We state the following lemma
(which should be a corollary
of the previous two)
with a slight abuse of language:
we say an invariant domain 
$D=\domcomplinv$
contains (is contained by) a set $C$
if
$\SU_0\cup\ldots\SU_{r-1}\supset C$
(respectively
$C\supset\SU_0\cup\ldots\SU_{r-1}$).

\medskip

\begin{lemma}\label{construir}
If $C$ is a c-invariant set contained
in some domain $D=\domcomplinv$
invariant
with respect to a system $\SF_{a,\ep}$
under parametric noise
satisfying hypothesis $A)$ of theorem \ref{teorema1},
then it contains
some s-invariant domain.
\end{lemma}

\begin{proof}
Let $C$ and $D$ be as stated and
let us consider $X=\om(C,\De)$
(cf. definition \ref{defswlimite}).

By lemma \ref{wlimiteDelta} we know
that
$X \subseteq \tra{C} \subset
\tra{\SU}_0\cup\ldots\tra{\SU}_{r-1}$
is a c-invariant set whose points are
carried cyclically through the $\tra{\SU}_i$,
$i=0,1,\ldots,r-1$.

By hypothesis $A)$ of theorem \ref{teorema1}
it holds that $\inter{X}\neq\emptyset$.
Thus the collection
$
\til{D}=( {\SU}_0\cap\inter{X} , \ldots,
{\SU}_{r-1}\cap\inter{X} )
$
is a member of $\SD$, an s-invariant domain.

Indeed, since the $f_t$ are diffeomorphisms
for all $t\in B$, the interior of $X$ must
be sent into the interior of $X$. But, by
lemma \ref{wlimiteDelta}, the orbits of
points of $X$ must respect the cyclic order of
the $\SU_i$, $i=0,\ldots,r-1$.

We conclude that $X$ contains an s-invariant
domain in its interior (the open sets
forming $\til{D}$ are pairwise separated
by construction).
Since $X\subset\tra{C}$,
we have the same for $C$.

\end{proof}

\begin{definition}\label{defGH}
Let $D=\domcomplinv$ be an
s-invariant domain ($D\in\SD$)
and $z\in M$. We define
$
G(z)=G_{D}(z) =
	\{ \subl{t} \in\De : 
		\exists n\in\natur \;\;\text{such that}\;\;
		f^n_{\subl{t}}(z)\in D
	\}
$
and $H(z)=H_D (z) = \De \setminus G(z)$,
the perturbation vectors that will
send $z$ into $D$ and those that
never do so, respectively.
\end{definition}

\begin{lemma}\label{H=Delta}
Let us suppose that $z\in M$ is such that
$\nuinf ( H_D (z) ) >0 $ for some
$D\in\SD$ and $\subl{t}$ is a
$H$-generic vector ($H=H(z)=H_D (z)$).

Then $H(w)=H_D (w) = \De$ for every
$w \in \om (z,\subl{t})$.
\end{lemma}

This lemma assures that those points whose
perturbed orbits never fall in some invariant domain $D$
for {\em many\/} ($\nuinf$-positive measure)
perturbations have $\om$-limit points
(under generic perturbations) which are {\em never\/}
sent into the same domain $D$ by every perturbation.
This is another ``independence of perturbation'' property
for the orbits of $\om$-limit points.

\begin{proof}
Let us fix a $H$-generic perturbation vector
$\subl{t}$ and $w\in\om(z,\subl{t})$.

By contradiction, let us suppose there are
$\subl{s}\in\De$ and $n\in\natur$
such that
$f^n_{\subl{s}} (w) \in D$.
Then there must be a neighborhood $U_w$
of $w$ in $M$ and  a neighborhood $V_{\subl{s}}$
of $\subl{s}$ in $\De$ such that
$
f^n ( U_w \times V_{\subl{s}} ) \subseteq D
$
by the continuity of
$f^n: M\times \De \mapto M$
(by property \ref{prop0}).

But $w\in\om(z,\subl{t})$ and $\subl{t}$
is $H$-generic, thus there are
$k\in\natur$ and $\subl{\te}\in H$
very close to $\subl{t}$,
with positive probability,
such that
$f^k (z,\subl{t}) \in U_w$
and
$\si^k\subl{\te} \in V_{\subl{s}}$
by lemma \ref{lema0-1>0}, since $\nuinf(H)>0$.
Therefore
$f^{k+n} (z,\subl{\te}) \in D$
contradicting
$\subl{\te}\in H$.
\end{proof}

\begin{lemma}\label{finalmentecai}
Let $z$ be a point of $M$ and $V$ a subset of $\De$,
with $\nuinf(V)>0$, such that 
for $\nuinf$-a.e. vector
$\subl{t}\in V$ and every $w\in \om (z,\subl{t} )$
there is $\subl{s}\in\De$
($\subl{s}=\subl{s}(\subl{t},w$))
such that the orbit $\SO(w,\subl{s})$
eventually falls in some minimal invariant
domain:
\[
\exists \SM=\SM(\subl{s}) \;\text{minimal}\;\;
\exists n=n(\subl{s})\in\natur:
f^n_{\subl{s}} (w) \in \SM.
\]

Then we will have
a $\preceq$-minimal domain $\SM$,
a set $W\subseteq V$, with
$\nuinf( W ) > 0$, and a $m\in\natur$
such that
$
f^m_{\subl{\te}} (z) \in \SM
$
for every
$\subl{\te}$ in $W$.
\end{lemma}

Let us observe that the hypothesis does not
prevent the point from being sent into different
invariant domains by different perturbations, but
the lemma ensures there will be a positive
measure set of perturbation vectors sending
the point {\em into the same invariant domain!\/}
In other words, the system under parametric noise
cannot be unstable to the extent of sending
a given point into completely different places
by nearby perturbations.

\begin{proof}
As in the proof of lemma \ref{lema0-1>0}
let us fix $\de >0$ and a compact $C$
contained in $V$ such that
$\nuinf( V\setminus C ) < \de$
and $E(1_V | \SB_j )_{|C}$
is continuous for every $j\ge1$.
We may assume $\nuinf( C ) >0$.

\medskip

Now we take $\subl{t}\in C$ such that $\subl{t}$ is
both $V$-generic and a $\nuinf$-density point of $C$.

\medskip

Let $w$ be a point in $\om(z,\subl{t})$ and
$\{ n_j \}_{j=1}^\infty\subseteq\natur$ a
sequence $n_1<n_2<\ldots$ such that
$f^{n_j}_{\subl{t}} (z) \to w$ when $j\to\infty$.
We will fix, from the hypothesis,
a minimal domain $\SM$,
an integer $k\in\natur$ and
a perturbation vector
$\subl{\te}\in\De$ such that
$f^k (w,\subl{\te}) \in \SM$.

\medskip

Since $\SM$ is open and $f^k:M\times\De\mapto M$
is continuous (see property \ref{prop0}), there are
neighborhoods $U_w$ of $w$ in $M$ and
$U_{\subl{\te}}$ of $\subl{\te}$ in $\De$
such that $f^k( U_w\times U_{\subl{\te}} ) \subseteq \SM$.

\medskip

By the choice of $w$ and $\subl{t}$ there is $m\in\natur$
with the property
\[
j\ge m \impl \left(
	f^{n_j}_{\subl{t}} (z) \in U_w
	\quad \text{ and } \quad
	\be =
	\nuinf( \si^{n_j} V(\subl{t},n_j) \cap U_{\subl{\te}} )
	> 0
	\right).
\]

Because $E(1_V | \SB_j )_{|C}$
is continuous, there is $\ro>0$ such that
\[
\subl{s}\in B(\subl{t},\ro)\cap C \impl
\left|
E(1_V|\SB_{n_m}) (\subl{s}) - E(1_V|\SB_{n_m}) (\subl{t})
\right| < \frac{\be}2
\]
and
$
f^{n_m}( z, B(\subl{t},\ro) ) \subseteq U_w
$ 
by the continuity of $f^{n_m}:M\times\De\mapto M$.

Then we have
$\nuinf( \si^{n_m} V(\subl{s},n_m) \cap U_{\subl{\te}} )
\ge \be/2 > 0$
for every
$\subl{s}\in B(\subl{t},\ro) \cap C$
because
\begin{eqnarray*}
\lefteqn{
\left|
\nuinf( \si^{n_m} V(\subl{t},n_m) \cap U_{\subl{\te}} )
\right.
-
\left.
\nuinf( \si^{n_m} V(\subl{s},n_m) \cap U_{\subl{\te}} )
\right| =
} \hspace{2cm}
\\
&=&
\nuinf \left[
\left(
U_{\subl{\te}} \cap \si^{n_m} V(\subl{t},n_m)
\right)
\difsim
\left(
U_{\subl{\te}} \cap \si^{n_m} V(\subl{s},n_m)
\right)
\right]
\\
&=&
\nuinf \left[
U_{\subl{\te}} \cap
\left( \;
\si^{n_m} V(\subl{t},n_m)
\difsim
\si^{n_m} V(\subl{s},n_m) \;
\right)
\right]
\\
&\le&
\nuinf \left( \;
\si^{n_m} V(\subl{t},n_m)
\difsim
\si^{n_m} V(\subl{s},n_m) \;
\right)
\\
&=&
\left|
\nuinf( \si^{n_m} V(\subl{t},n_m) )
-
\nuinf( \si^{n_m} V(\subl{s},n_m) )
\right| < \frac{\be}2
\end{eqnarray*}

It follows that
$
W = \bigcup\limits_{\subl{s}\in B(\subl{t},\ro)\cap C}
\left\{
(s_1,\ldots,s_{n_m},u_1,u_2,\ldots) :
\subl{u}\in
\left(
\si^{n_m} V(\subl{s},n_m) \cap U_{\subl{\te}}
\right)
\right\}
$
is a subset of $V$ such that
\begin{eqnarray*}
\nuinf(W) &=&
\int_{ B(\subl{t},\ro)\cap C }
\nuinf \left(
\si^{n_m} V(\subl{s},n_m) \cap U_{\subl{\te}}
\right) \;
d\nu^{n_m} ( \subl{s} )
\\
&\ge& \frac{\be}2 \cdot
\nuinf \left(
B(\subl{t},\ro)\cap C
\right)
> 0
\end{eqnarray*}
because $\subl{t}$ is a $\nuinf$-density
point of $C$.
Moreover
\[
f^{n_m}( z,W ) \subseteq
f^{n_m} \left(
z, ( B(\subl{t}, \ro) \cap C )
\right) \subseteq
f^{n_m} ( z, B(\subl{t}, \ro) )
\subseteq U_w
\] 
and
$f^k( U_w\times U_{\subl{\te}} ) \subseteq \SM$
with
$\si^{n_m} W \subseteq U_{\subl{\te}}$.
Thus
$
f^{n_m+k} ( z,W ) \subseteq \SM,
$
completing the proof of the lemma.
\end{proof}

\section{Finite Number of Minimal Invariant Domains}
\label{numerofinitodeminimais}

Two basic properties of the members of $\SD$
are the following direct consequences of
hypothesis $A)$ of theorem \ref{teorema1}
and definitions \ref{defdominv} and \ref{defdomsimetr}.

\begin{property}\label{prop1}
Any s-invariant domain $D=\domcomplinv$ is such
that every open set $\SU_i$ contains some ball of radius
$\xi_0>0$, $i=0,\ldots,r-1$.
Consequently, each open set has a volume ($m$ measure)
greater than some constant $l_0>0$.
\end{property}

\begin{property}\label{prop2}
The period of any invariant domain
$D\in\SD$ is bounded from above
by a constant $T_{p}\in\natur$ dependent on
$l_0$ ($T_{p} \le 1/l_0$).
\end{property}

\subsection{Minimals Exist}\label{minimalexist}

We start by showing that Zorn's lemma can be applied
to the partially ordered set $(\SD,\preceq)$ of
completely and symmetrically invariant domains of $M$.
Having established this, we conclude that there are
minimal invariant domains in $M$.

\medskip

Let $\SC$ be a $\preceq$-chain in $(\SD,\preceq)$,
that is, if $D,D'\in\SC$ then either $D\preceq D'$
or $D'\preceq D$. By property \ref{prop2},
the domains of $\SC$ have a finite number of
distinct periods.
So if $\rho:\SC\mapto\natur$ is the map that associates
to each $D\in\SC$ its period $\rho(D)\in\natur$, then
$\rho(\SC)=\{ r_1, \ldots, r_l \}$ and
$\SC = \bigcup\limits_{j=1}^l \rho^{-1} \{ r_j \}$.
We need to find a lower bound for $\SC$ in $(\SD,\preceq)$.
We can suppose that $\SC$ does not have a minimum,
otherwise we would have nothing to prove.
Now we establish

\begin{claim}
There is a $j_0\in\{1,\ldots,l\}$ such that the subchain
$\Ss = \rho^{-1} \{ r_{j_0} \}$ does not have a
lower bound in $\SC$. Moreover $\Ss$ {\em precedes\/}
every element of $\SC$: for all
$D\in\SC$ there is a $ D'\in\Ss$ such that $D'\preceq D$.
\end{claim}

Indeed, if every subchain of constant period
$\Ss_j=\rho^{-1} \{ r_j \}$ had a lower bound
$D_j\in\SC$ for $j=1,\ldots,l$,
then the minimum of the subchain
$\Ss'=\{ D_1,\ldots,D_l \} \subseteq \SC$
(which always exists because $\Ss'$ is finite)
would be a minimum for $\SC$,
in contradiction to the supposition
we started with. So there is some $\Ss=\Ss_{j_0}$
without a lower bound in $\SC$.

\medskip

Now for the second part of the claim. Let us suppose,
by contradiction, that there is a $\til{D}\in\SC$ such that
$D \not\preceq \til{D}$ for every $D\in\Ss$. But we are
within a chain, thus $\til{D} \prec D$ for all $D\in\Ss$,
that is, $\til{D}$ would be a lower bound for $\Ss$
in $\SC$, and this contradiction proves the claim.

\medskip

Now we just need to show that $\Ss$ has some
lower bound in $(\SD,\preceq)$ in order to get a
lower bound for $\SC$.

\begin{figure}[h]
                 
\centerline{\psfig{figure=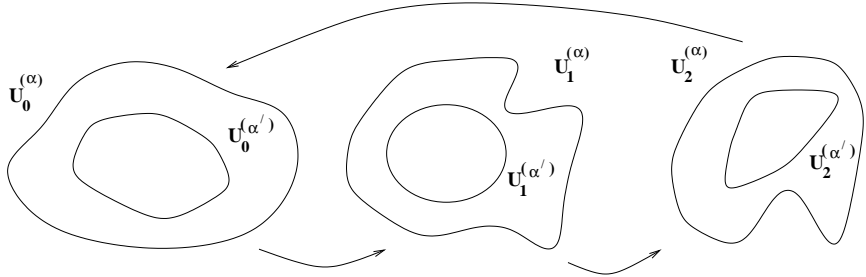,width=4.5in}}
                 
\caption{{\em
\label{encaixados} $D_{\al}\prec D_{\al'}$ with
$D_{\al}$,$D_{\al'}$ in a subchain of period three
after suitable arrangement of indexes}}
\end{figure} 

To do that, let us first observe that $\Ss$ is made
by nested invariant domains of equal period, all
symmetrically invariant. Thus we can always write
$D\in\Ss$ as $D=\domcomplinv$ and, for any other
$D'=\domcomplinvl$, we can never have two different
$\SU_i',\SU'_j$ intersect the same $\SU_k$,
$i,j,k\in\{ 0,\ldots,r-1 \}$ and $i\neq j$
(see figure \ref{encaixados} for a representation
of $\Ss$ with period three).

\medskip

Hence we can rearrange the lower indexes of the open sets
that form the domains of $\Ss$ in order to
obtain
$\Ss=\{ D_\al \}_{\al\in\SA}$
with
$D_\al=( \SU_0^{(\al)}, \ldots, \SU_{r-1}^{(\al)} )$
for $\al\in\SA$,  $\SA$ some set of indexes,
and satisfying the following property
\[
D_\al \preceq D_{\al'} \equi
	\left(
	\SU_i^{(\al)} \subseteq \SU_i^{(\al')}, \;
	i=0,\ldots,r-1
	\right) ;
\]
for all $\al,\al' \in \SA$.

We can now consider the intersections
$\til{\SU}_i = \bigcap\limits_{\al\in\SA}
\SU_i^{(\al)},\; i=0,\ldots,r-1$,
and observe that, because each $D_\al$ is
s-invariant, the family $\domcomplinvtil$ 
satisfies
\begin{equation}\label{simetria}
f^k (\til{\SU}_i,\De) \subseteq \til{\SU}_{ (k+i) \bmod r },
\;\;\forall k\ge1 \;\; \forall i=0,\ldots,r-1
\end{equation}
and since fixing $\al_0\in\SA$ we have
$\til{\SU}_i\subset\SU_i^{(\al_0)}$ for
$i=0,\ldots,r-1$, the
$\til{\SU}_0,\ldots,\til{\SU}_{r-1}$ are
pairwise separated, because the
$\SU_0^{(\al_0)},\ldots,\SU_{r-1}^{(\al_0)}$
already were pairwise separated.

Finally, hypothesis $A)$ of theorem \ref{teorema1}
and \refer{simetria} ensure that
every $\til{\SU}_i$ has nonempty interior
($i=0,\ldots,r-1$).
Since the $f_t$ are diffeomorphisms for $t\in B$,
hence open maps, we deduce that
$
\til{D} = ( \inter{\til{\SU}_0},\ldots,\inter{\til{\SU}_{r-1}} )
$
is an s-invariant domain of
$\SD$ which clearly is a lower bound for the subchain
$\Ss$. Consequently we got a lower bound for the chain
$\SC$ we started with and proved that Zorn's lemma
can be applied to $(\SD,\preceq)$.

\medskip

Moreover,
it is easy to see that
{\em each member of $\SD$ contains
a minimal domain\/}.

\medskip

In fact, let us now fix $D_0\in\SD$ and consider the
partially ordered set $(\SD_{D_0},\preceq)$, where
$
\SD_{D_0} = \{ D\in\SD : D\preceq D_0 \}.
$
Since it can be shown that each chain of
$(\SD_{D_0},\preceq)$ has a lower bound in
$\SD_{D_0}$,
in the same way we did before, there must be
some minimal domain in $(\SD_{D_0},\preceq)$ which,
by the definition of $\SD_{D_0}$, is also a minimal
domain of $(\SD,\preceq)$.

\medskip

We conclude that each domain in $\SD$ contains a
minimal domain of $(\SD,\preceq)$.

\subsection{Minimals are Pairwise Disjoint}
\label{minimais-disjuntos}

Let us now observe that, because each open set of the
collection that forms an invariant domain
has a volume (Riemannian measure $m$ on $M$) of
at least $l_0>0$ by property \ref{prop1},
to prove there is a finite
number of $\preceq$-minimals we need only show they are
pairwise disjoint.

\medskip

Let $D=\domcomplinv$ and $D'=\domcomplinvv$ be
two minimals of $(\SD,\preceq)$ whose open sets
have some intersection,
$\SU_i \cap \SU_j' \neq \emptyset$ say,
for some $i\in\{ 0,\ldots,r-1 \}$ and
$j\in\{ 0,\ldots,r'-1 \}$.

Because both $D$ and $D'$ are
s-invariant, we have for all $k\ge1$
\[
f^k_{\De} ( \SU_i \cap \SU_j' )  \subset 
f^k_{\De} ( \SU_i ) \subset \SU_{ (i+k) \bmod r } \;
\text{ and } \;
f^k_{\De} ( \SU_i \cap \SU_j' )  \subset 
f^k_{\De} ( \SU_j' ) \subset \SU_{ (j+k) \bmod r'}'
\]
and thus
$
f^k_{\De} ( \SU_i \cap \SU_j' )  \subset
\SU_{ (i+k) \bmod r } \cap \SU_{ (j+k) \bmod r'}'  .
$
Therefore if we define
\[
\hat{D} =
( \SU_i \cap \SU_j',
\SU_{ (i+1) \bmod r } \cap \SU_{ (j+1) \bmod r'}',
\ldots,
\SU_{ (i+[r,r']-1) \bmod r } \cap \SU_{ (j+[r,r']-1) \bmod r' }' )
\]
we will get $\hat{D} \in \SD$
(here $[r,r']$ is the least common multiple of
$r$ and $r'$).

The invariance property is clear.
Let us check that the open sets forming $\hat{D}$
are pairwise {\em separated\/}. Indeed, if
we had
$
\tra{\SU_{ (i+k_1) \bmod r } \cap \SU_{ (j+k_1) \bmod r' }' }
\cap
\tra{\SU_{ (i+k_2) \bmod r } \cap \SU_{ (j+k_2) \bmod r' }' }
\neq \emptyset
$
with $0 \le k_1 < k_2 \le [r,r']-1$
then, in particular,
\[
\tra{\SU}_{ (i+k_1) \bmod r }
\cap
\tra{\SU}_{ (j+k_2) \bmod r }
\neq \emptyset \;
\text{ and } \;
\tra{\SU}_{ (i+k_1) \bmod r' }'
\cap
\tra{\SU}_{ (j+k_2) \bmod r' }'
\neq \emptyset.
\]

However by definitions \ref{defdominv}
and \ref{defdomsimetr}
we conclude that
$
k_1 \equiv k_2 \pmod{r}
$
and
$
k_1 \equiv k_2 \pmod{r'}
$
with $0 \le k_1 < k_2 \le [r,r']-1$,
contradicting the Chinese Remainder Theorem.

\medskip

We have now $\hat{D}\preceq D$ and $\hat{D}\preceq D'$,
so the minimality of both $D$ and $D'$ implies
$D=\hat{D}=D'$. We have shown that if two $\preceq$-minimals
intersect then they are equal. Consequently, we have that
they are pairwise disjoint and, as mentioned above,
we conclude there is a finite number of minimals in
$(\SD,\preceq)$.

\subsection{Minimals are Transitive}
\label{transitivos}

The following is an expression of the
dynamical indivisibility of minimal
invariant domains.

\begin{lemma}\label{r-transitiveness}
Every minimal invariant domain $\SM=\domcomplinv$
is transitive in the following sense.
For every $z\in\SM$
(meaning $z\in \SU_0\cup\ldots\cup\SU_{r-1}$)
the sequence $\{ f^n_{\De} (z)  \}_{n=1}^\infty$
is dense in $\SM$.
\end{lemma}

We will say that minimal invariant domains are
{\em randomly transitive\/} or
{\em r-transitive\/} when referring to
this kind of transitiveness.

\begin{proof}
In fact, let $\SM=\domcomplinv\in\SD$
be a minimal and let us take some point
$z\in\SU_i$ with $i\in \{0,\ldots,r-1\}$
and $X = \om ( z,\De )$
(cf. definition \ref{defswlimite}).

By lemma \ref{wlimiteDelta},
we have
$X\subseteq\tra{\SM}=\tra{\SU}_0
\cup\ldots\cup\tra{\SU}_{r-1}$,
$X$ is c-invariant and
goes cyclically through the
$\tra{\SU}_0,\ldots,\tra{\SU}_{r-1}$,
under every perturbation vector of $\De$.
Besides, by lemma \ref{construir} there
is $D\in\SD$ such that $D\subset X$.
So $D\preceq\SM$,
in contradiction with the $\preceq$-minimality
of $\SM$.    

Hence it must be that $\SM=D$ and
then $\{ f^n_{\De} (z)  \}_{n=1}^\infty$ is dense
in $\SU_0\cup\ldots\SU_{r-1}$,
as stated.   
\end{proof}

\medskip

Given a minimal $\SM=\domcomplinv$, since it is
s-invariant, it will also be
invariant with respect to $f_t$ for every
$t\in T$, because the vector
$(t,t,t,\ldots)$ is in $\De$.

This means we have
$
f^k_t ( \SU_i ) \subset \SU_{ (i+k) \bmod r }
$
for all $k\ge1$ and $i=0,\ldots,r-1$.

\medskip

However, we cannot state any kind of indivisibility
for this domain with respect to $f_t$ because
the domain was originally a minimal domain,
but with noise. The perturbations around the
system $f_a$ may have mixed, in a single collection
of open sets, several attractors indivisible
with respect to $f_t$, but that under random choices
of parameters were indistinguishable. We cannot
proceed further in this because we made no
hypothesis about the dynamics of the $f_t$
without noise.

\section{Stationary Probability Measures}
\label{ergodica}
 
\subsection{Existence and Absolute Continuity}
\label{abscont}

Let $z$ be a point of $M$.
The formalization of the dynamics under noise
by means of the operator $S$ enables us to
naturally associate a probability measure to
the orbits of the system:
the push-forward of $\nuinf$ from $\De$
to $M$ via the map $f^k$ given by
$f^k(z,\nuinf)$, $k\ge1$. We have defined
this as the probability which integrates
continuous functions $\vfi:M\mapto\real$
as
$
f^k(z,\nuinf) \vfi=
\left[
f^k(z,\cdot)_\ast \nuinf
\right]
\vfi = \int \vfi
\left(
f^k(z,\subl{t})
\right)  \,
d\nuinf (\subl{t}), \;\; k\ge1.
$

These probabilities are not stationary
in general, but if we consider their averages
\begin{equation}\label{average}
\mu_n = \frac1n \sum\limits_{j=1}^{n}
f^j(z, \nuinf)  ,\;\; n= 1, 2,\ldots,
\end{equation}
we obtain a sequence of probability measures
in $M$ which,
by compactness of the space
${\Bbb P}(M)$ of probabilities measures over $M$
with the weak topology,
has some limit point
$
\mu_\infty = \lim\limits_{i\to\infty}
\frac{1}{n_i} \sum\limits_{j=1}^{n_i}
f^j(z,\nuinf) .
$
This means the integral of a continuous
$\vfi:M\mapto\real$ with
respect to $\mu_\infty$ is given by
$
\mu_\infty ( \vfi ) = \lim\limits_{i\to\infty}
\frac{1}{n_i} \sum\limits_{j=1}^{n_i}
\int \vfi
\left(
f^j(z,\subl{t})
\right)  \,
d\nuinf (\subl{t}) .
$

This accumulation point is a stationary probability.
In fact,
\[
\int \int \vfi\circ f(w,s) \,
d\mu_\infty (w) \, d\nu(s)
=
\int
\left[
\lim\limits_{i\to\infty}
\frac{1}{n_i}
\sum\limits_{j=1}^{n_i}
\int ( \vfi\circ f_s )
\left( f^j(z,\subl{t}) \right) \,
d\nuinf (\subl{t})
\right] \,
d \nu (s)
\]
and
\begin{eqnarray*}
\lefteqn{
\int
\left[
\frac{1}{n_i}
\sum\limits_{j=1}^{n_i}
\int ( \vfi\circ f_s )
\left( f^j(z,\subl{t}) \right) \,
d\nuinf (\subl{t})
\right] \,
d \nu (s) =
} \hspace{-2cm}
\\
& & =
\frac{1}{n_i} \sum\limits_{j=1}^{n_i}
\int \int \vfi
\left(
f_s\circ f_{t_j}\circ \ldots \circ f_{t_1} (z)
\right) \,
d\nuinf (\subl{t}) \, d \nu (s)
=
\frac{1}{n_i} \sum\limits_{j=1}^{n_i}
\int \vfi
\left(
f^{j+1}_{\subl{t}} (z)
\right)
\, d\nuinf (\subl{t})
\\
& & =
\frac{1}{n_i} \sum\limits_{j=1}^{n_i}
\int \vfi
\left(
f^{j}_{\subl{t}} (z)
\right)
\, d\nuinf (\subl{t})
+
\frac{1}{n_i}
\left[
\int \vfi
\left(
f^{n_i+1}_{\subl{t}} (z)
\right)
\, d\nuinf (\subl{t})
-
\int \vfi
\left(
f_{\subl{t}} (z)
\right)
\, d\nuinf (\subl{t})
\right]
\end{eqnarray*}
for $i\ge1$.
Since 
$\sup\limits_{w\in M} |\vfi(w)| = \| \vfi \|$
is finite, the second term of the last expression
converges to zero when $i\to\infty$, while
the first term gives the integral of $\vfi$ with
respect to $\mu_\infty$,
that is
\begin{eqnarray*}
\int \vfi \, d\mu_\infty
&=&
\lim\limits_{i\to\infty}
\int
\left[
\frac{1}{n_i}
\sum\limits_{j=1}^{n_i}
\int ( \vfi\circ f_s )
\left( f^j(z,\subl{t}) \right) \,
d\nuinf (\subl{t})
\right] \,
d \nu (s)
\\
&=&
\int
\left[
\lim\limits_{i\to\infty}
\frac{1}{n_i}
\sum\limits_{j=1}^{n_i}
\int ( \vfi\circ f_s )
\left( f^j(z,\subl{t}) \right) \,
d\nuinf (\subl{t})
\right] \,
d \nu (s)
\\
&=&
\int \int \vfi\circ f_s(w) \,
d\mu_\infty (w) \, d\nu(s),
\end{eqnarray*}
where we have used the dominated convergence
theorem to exchange the limit and the
integral signs.
In addition,
because $C^0(M,\real)$ is dense
in $L^1(M,\mu_\infty)$ with
the $L^1$-norm, we see the last
identity holds for every
$\mu_\infty$-integrable
$\vfi:M\mapto\real$.

\medskip

Moreover, if $E$ is any Borel subset of $M$
we can write
\begin{eqnarray*}
\mu_\infty (E)
&=&
\int 1_E \, d\mu_\infty
=
\int \int 1_E( f_t (x) ) \, d\mu_\infty(x) \, d\nu(t)
\\
&=&
\int \int \int 1_E
\left(
f_{t_2}\circ f_{t_1} (x)
\right) \, d\mu_\infty(x) \, d\nu(t_1) \, d\nu(t_2)
\\
&\vdots&
\\
&=&
\int \int 1_E( f^k_{\subl{t}}(x) ) \, d\mu_\infty(x)
\, d\nuinf(\subl{t})
\\
&=&
\int \int 1_E( f^k_{\subl{t}}(x) ) \, d\nuinf(\subl{t})
\, d\mu_\infty(x)
\\
&=&
\int f^k(x,\nuinf)(E) \, d\mu_\infty(x)
\end{eqnarray*}
for $k=1,2,3 \ldots$
Hypothesis $B)$ of theorem \ref{teorema1}
guarantees that $f^k(x,\nuinf)\ll m$ for $k\ge K$.
Thus $\mu_\infty (E) = 0$ whenever $m(E)=0$.
We have just proved

\begin{lemma}\label{stationary+abscont}
Given $z\in M$, any accumulation point of the
averages \refer{average} is a stationary
absolutely continuous probability measure over $M$.
\end{lemma}

Let us remark that $\mu_\infty=\mu_\infty(z)$
depends on $z\in M$ and the accumulation point
of the averages \refer{average} may not
be unique.

\subsection{Ergodicity and Characteristic Probabilities}
\label{ergodicidade}

Let us suppose $z\in D$ for some $D\in\SD$. Then it
is clear that $\supp{\mu_\infty}\subset\tra{D}$,
whatever accumulation point of the averages
\refer{average} we choose.
Moreover, by remark \ref{closure-s-inv} we have
that $\tra{\SD}=(\tra{\SU}_0,\ldots,\tra{\SU}_{r-1})$
satisfies \refer{s-invariant} also. Thus
if $w\in\supp{\mu_\infty}$, we get by hypothesis $A)$
that
$
f^k(w,\De)\supset B(f^k_{t_0}(w),\xi_0)
$
for all $k\ge K$,
and by the invariance of the support we conclude
$\supp{\mu}\cap D\neq\emptyset$ because
$\inter{\pa D}=\inter{\pa\SU_0\cup\ldots\cup\pa\SU_{r-1}}
=\emptyset$.
In addition, if $z$ belonged to a minimal $\SM\in\SD$,
then the invariance of the support, the fact
that $\supp{\mu}\cap D\neq\emptyset$ and
the r-transitiveness of $\SM$ (given by
lemma~\ref{r-transitiveness}) together imply
$\supp{\mu}=\tra{\SM}$.

\begin{lemma}\label{ergodiclaim}
If $\SM\in\SD$ is a minimal invariant domain
and $\mu$ a stationary absolutely continuous
probability measure with 
$\supp{\mu}=\tra{\SM}$, then
\[
\vfi(x)=\int \vfi \left( f_t(x)  \right) d\,\nu(t),\,
\mu-\text{a.e.}\;x \impl
\vfi \;\; \text{is} \;\;
\mu-\text{a.e. constant}
\]
for every bounded measurable function
$\vfi:M\mapto\real$.
\end{lemma}

\begin{proof}
What we want to prove
is equivalent to the following
for every Borel set $E$:
\begin{equation}\label{ec2}
1_E(x)=f(x,\nu) 1_E \,
\mu-\text{a.e.}\;x \impl
\mu (E)=0 \;\; \text{or} \;\; 1 .
\end{equation}

Let $E$ be a Borel set
that satisfies the left
hand side of \refer{ec2} and let us suppose
that $\mu(E)>0$.
Let $\delta>0$ be given and take an open set $A\supset E$
such that $\mu(A\setminus E)<\delta$. Then the following
holds $\mu$-a.e. $x$
\[
1_E(x) = f(x,\nu) 1_E = f(x,\nu) f(x,\nu) 1_E = f^2(x,\nu^2) 1_E.
\]
We can iterate any number of times obtaining
$1_E(x)=f^k(x,\nu^\infty) 1_E$ for all $k\ge1$ and
$\mu$-a.e. $x$.

In particular, if $x\in E$, this means that
$f^k(x,\un t)\in E \subset A$ for $\nu^\infty$-a.e. $\subl{t}$
and every given $k\ge1$. By Property~\ref{prop0}(2) and because
a full $\nu$-measure set is dense in $\Delta$ we have that
\[
f^k(x,\subl{t}) \in \tra{A} 
\quad\mbox{for all}\quad k\ge1 \quad\mbox{and every}\quad
\subl{t}\in\Delta
\]
whenever $x\in E$.
Then, if we define
${\SU}'_i=\SU_i\cap \cup_{x\in E,k\ge K}
\inter{f^k(x,\De)}$,
$i=0,\ldots,r-1$,
we see that the ${\SU}'_i\subset A$ are
open, nonempty (by hypothesis (A) and because
$\supp{\mu}=\tra{\SM}$ and
$\inter{\tra{\SM}} = \SM = \SU_0\cup\ldots\cup\SU_{r-1}$)
and so $D=\domcomplinvl$ is an s-invariant domain.

\medskip

In fact, fixing $y\in{\SU}'_i$ for some
$0\le i \le r-1$, $\subl{s}\in\De$ and
$n\ge1$, there are $k\ge K$ and $\de>0$
such that
$B(y,\de)\subset f^k(x,\De)$ and
$f^n_{\subl{s}}( B(y,\de) )
\subset f^{k+n}(x,\De) \subset \SU_{i+n \bmod r}$
by definition of ${\SU}'_i$.
Hence
$f^n_{\subl{s}}(y)\in \inter{ f^{k+n}(x,\De) }
\cap \SU_{i+n \bmod r}$
after property~\ref{prop0}(3).

\medskip

We have
built an s-invariant domain $D\in\SD$
such that $D\preceq\SM$.
The minimality of $\SM$ gives $D=\SM$ and hence
$\mu(A)=1$, that is,
$\mu(E)\ge \mu(A)-\delta = 1-\delta$.
Since $\delta>0$ was arbitrary, the proof is complete.
\end{proof}

\medskip
\medskip

Lemma~\ref{ergodiclaim} implies that $\mu_\infty$ is
ergodic, that is, $\mu_\infty\times\nuinf$ is
$S$-ergodic. (For ease of writing we make
$\mu=\mu_\infty$ in the following discussion.)

\medskip

Indeed, let us assume that
$\psi:M\times\De \mapto \real$
is an $S$-invariant bounded
measurable function:
$\psi\left( S(z,\subl{t}) \right) =
\psi( z,\subl{t} )$,
$\mu\times\nuinf$-a.e
$(z,\subl{t})\in M\times\De$.
 
For each $k\ge0$ we define
\[
\psi_k (x,t_1,\ldots,t_k) =
\int
\psi (x,t_1,t_2,\ldots,t_k,t_{k+1},\ldots)
\, d\nu(t_{k+1}) \, d\nu(t_{k+2}) \ldots
\]
and we have, by the invariance of $\psi$,
\begin{eqnarray*}
\psi_0 (x) &=&
\int \psi (x,t_1,t_2,\ldots) \,
d\nu(t_1) \, d\nu(t_2) \ldots
\\
&=&
\int
\psi ( f_{t_1}(x), t_2, t_3, \ldots ) \,
d\nu(t_2) \, d\nu(t_3) \ldots d\nu(t_1)
\\
&=&
\int \psi_0 \left( f_{t_1} (x) \right) \, d\nu(t_1),
\quad \mu-\text{a.e.} \;\; x\in M.
\end{eqnarray*}
 
Therefore, by lemma~\ref{ergodiclaim}, we conclude
that $\psi_0$ is $\mu$-a.e. constant.
In general, for $k\ge1$,
\begin{eqnarray*}
\psi_k (x,t_1,\ldots,t_k)
&=&
\int
\psi (x,t_1,t_2,\ldots,t_k,t_{k+1},\ldots)
\, d\nu(t_{k+1}) \, d\nu(t_{k+2}) \ldots
\\
&=&
\int
\psi ( f_{t_1}(x), t_2,\ldots, t_k, t_{k+1}, \ldots ) \,
d\nu(t_{k+1}) \, d\nu(t_{k+2}) \ldots
\\
&=&
\psi_{k-1}
\left( f_{t_1}(x), t_2, \ldots, t_k \right),
\;\; \mu\times\nu^k-\text{a.e.} \;\;
(x,t_1,\ldots,t_k).
\end{eqnarray*}
 
We then have $\psi_1\equiv\psi_0,\;\;\mu\times\nu$-a.e.;
$\psi_2\equiv\psi_1,\;\;\mu\times\nu^2$-a.e.;$\ldots$
and so, by induction
\[
\psi_k\equiv\psi_0\equiv\text{constant},\;\;
\mu\times\nu^k-\text{a.e., for every } k\ge1.
\]
 
However
if we identify
$\psi_k (x,\subl{t})$
with $\psi_k(x,t_1,\ldots,t_k)$,
then $\psi_k$ coincides with
$E(\psi| \SB_k)$,  $\mu\times\nuinf$-a.e.
and we have seen in lemma \ref{lema0-1}
that $E(\psi| \SB_k) \mapto \psi,
\;\; \mu\times\nuinf$-a.e.,
when $k\to\infty$.
Hence we have also $\psi\equiv\text{constant},
\;\; \mu\times\nuinf$-a.e., and
conclude that $\mu\times\nuinf$ is $S$-ergodic.

\medskip

Ergodicity, Birkhoff's ergodic theorem and the 
absolute continuity imply that $\mu=\mu_\infty$
is physical. Indeed $\mu(E(\mu))=1$ by Birkhoff's
ergodic theorem, so $m(E(\mu)) > 0$ by the
absolute continuity of $\mu$ with respect to $m$.

\medskip

We easily deduce that any two physical
probability measures $\mu_1, \mu_2$
whose support is $\tra{\SM}$
must be equal.

We notice first that for any given $x\in\tra{\SM}$
the union $\cup_{k\ge K} f^k(x,\Delta)$ contains $\SM$.
This is clear since this union is easily seen to be
completely invariant, which implies that it must
contain an invariant domain, so it contains $\SM$
by the minimality property.

Let $x\in E(\mu_1)\cap\SM$ and let $y\in E(\mu_2)\cap\SM$ be
a density point of $E(\mu_2)$ with respect to $m$
(it exists since $m(E(\mu_2)\cap\SM)>0$). By the argument
in the previous paragraph there are $k\ge K$ and $\subl{t}$ such that
$f^k(x,\subl{t})=y$.

Property~\ref{prop0}(1) ensures that $m\ll f^k(x,\nu^\infty)$.
The choice of $y$ implies that
for a given $\eta>0$ and a sufficiently small $\delta$ we have
$m(B(y,\delta)\cap E(\mu_2))\ge (1-\eta) m(B(y,\delta))>0$.
Thus $f^k(x,\nu^\infty)(E(\mu_2)\cap B(y,\delta))>0$.

According to the definition of ergodic basin, this implies
that $\mu_1=\mu_2$.

\medskip

The above arguments prove the existence of
a {\em characteristic measure} for each
minimal invariant domain. One extra very usefull property
can also be deduced.

\begin{lemma}\label{lema:suppnabacia}
For every stationary ergodic $\mu$ we have $\supp(\mu)\subset B(\mu)$.
\end{lemma}

\begin{proof}
We start observing that Hypothesis (A) and Property~\ref{prop0}(1)
imply that there is a $\xi$-ball Lebesgue-a.e. inside the basin of $\mu$.
Indeed, taking $x\in B(\mu)$ we know that $f^k(x,\subl{t})\in B(\mu)$
for all $k\ge1$ and $\nuinf$-a.e. $\subl{t}$. For $k\ge K$ let
$B=B(f^k(x),\xi)$. Then $B\subset\inter{\supp(\mu)}$. By contradiction,
if there existed some $E\subset B\setminus B(\mu)$ with $m(E)>0$,
then because $m\ll f^k(x,\nu^\infty)$ we would have $f^k(x,\nu^\infty)(E)>0$.
So we conclude that  $m(B\setminus B(\mu))=0$.

For future use set $\vfi\in C^0(M,\real)$ such that $0\le \vfi\le1,
\vfi\mid(M\setminus B)\equiv 0$ and $\vfi\mid B(f^k(x),\xi/2)\equiv 1$,
so there is $c>0$ such that $\mu(\vfi)>c$.

Now we know from the arguments prior to the statement of the lemma
that every $x\in\supp(\mu)$ admits $k\ge1$ and a subset $V_{x,k}$ with
$\nuinf(V_{x,k})>0$ such that $f^k(x,V_{x,k})\subset B(\mu)$.
Fixing $x\in\supp(\mu)$, what we want to show is that
\[
W=\{ \subl{t}: f^k(x,\subl{t})\in E(\mu) \quad \mbox{for some}\quad k\ge1\}
\]
has full $\nu^\infty$-measure.
Let us suppose that $F=\Delta\setminus W$ is satisfies $\nu^\infty(F)>0$ 
and take $\subl{t}\in F$. 
Then for $k\ge K$ the point $y=f^k(x,\subl{t})$ is in $\supp(\mu)$ 
and we may assume that $\nuinf(\si^k F(\subl{t},k))\ge 1-c/10$
because $\nuinf(F)>0$ and by Lemma~\ref{lema0-1}.
Since $\mu(\vfi)>c$ there is $\ell\ge1$
such that $f^\ell(y,\nuinf)\vfi\ge c/2$.
Hence 
\[
\nuinf(\{\subl{s}:f^k(y,\subl{s})\in B\}\cap \si^k F(\subl{t},k))>0
\]
and $\nuinf$ almost every vector in this set sends $y$ in $B(\mu)$
which contradicts the definition of $F$.
\end{proof}

Along this subsection we have proved

\begin{proposition}\label{physical}
Given a minimal $\SM\in\SD$ there is only
one physical absolutely continuous
probability measure whose
support is contained in $\tra{\SM}$.
Moreover, every $x\in\tra{\SM}$ is in the
ergodic basin of this 
characteristic measure.
\end{proposition}

\section{Decomposition of Stationary Probabilities}
\label{decomposition}

Let $\mu$ be a stationary probability.
Then $\supp\mu$ is a c-invariant set.
By hypothesis $A)$ of theorem \ref{teorema1}
we deduce that $\inter{\supp\mu}\neq\emptyset$.

Let $C_1, C_2,\ldots$ be the connected components
of $\inter{\supp\mu}$: it is an at most countable family
of connected sets and $\inter{\supp\mu}=\cup_{i\ge1} C_i$.
 
\medskip
 
Since $f_t$ is a diffeomorphism for every $t\in T$,
thus a continuous open  map,
we deduce that each $f_t(C_i)$ is a connected
open set contained in $\supp\mu$, by the c-invariance.
Hence there is some $j=j(i,t)$ such that
$f_t(C_i)\subset C_j$ by openness and
connectedness.
 
In particular, by the same reasoning, we see
that every point in $C_i$ is sent by $f_t$
in the interior of $\supp\mu$ for all $t\in T$ and
$i\ge1$.
 
\medskip
 
We show that $j=j(i,t)$ does not
depend on $t\in T$.
 
\medskip
 
By contradiction, let us suppose there are
$i\ge1$, $t_0$ and $t_1$ in $B$ such that
$j_0=j(i,t_0)\neq j(i,t_1)=j_1$
and let us fix $x\in C_i$. 
We take a continuous curve
$\ga: [0,1] \mapto T$ with endpoints
$t_0$ and $t_1$ in $B$:
$\ga(0)=t_0$ and $\ga(1)=t_1$.
We know that
\[
f(x,\ga(s)) \in \inter{\supp\mu} = \cup_{i\ge1} C_i
\;\; \text{for all} \;\; s\in[0,1],
\]
but since $f(x,\ga(0))=f(x,t_0)\in C_{j_0}$
and $f(x,\ga(1))=f(x,t_1)\in C_{j_1}$ with
$C_{j_0}, C_{j_1}$ distinct connected components
of $\inter{\supp\mu}$, we conclude there is
$\tra{s}\in]0,1[$ such that
\[
f(x,\ga(\tra{s}))\in \pa C_{j_0}
\subset \pa (\supp\mu) = \supp\mu\setminus(\inter{\supp\mu}),
\]
a contradiction. So every $C_i$ is sent into
some $C_{j(i)}$ by any $f_t$ and the permutation
$i\mapsto j(i)$ does not depend on $t\in T$.

We remark, in particular, that
if for $x\in C_i$
we have $f^k(x,\subl{t})\in C_j$ for some $j$,
$k\ge1$ and $\subl{t}\in\De$, then
$f^k(x,\De)\subset C_j$.
 
\medskip
 
Since $\mu\times\nuinf (C_i\times\De) >0$
($i\ge1$) Poincar\'e's recurrence theorem
guarantees that $\mu\times\nuinf$-a.e.
pair $(x,\subl{t})\in C_i\times\De$ is
$\om$-recurrent with regard to the action of $S$.
By last remark, we see that $f^k(C_i,\De)$
returns to $C_i$ infinitely often, for every
fixed $i$.
Hence, again by hypothesis $A)$, each $C_i$
contains a $\xi_0$-ball.
Thus, because $M$ is compact,
the pairwise disjoint family
$C_1,C_2,\ldots$
must be finite and so
$\inter{\supp\mu}= C_1 \dot{\cup} \ldots
\dot{\cup} C_l$ (a disjoint union).

\medskip

The open sets $C_1,\ldots,C_l$ may not
be pairwise separated. However, the following
reflexive and symmetric relation
$C_i\sim C_j \Leftrightarrow
\tra{C_i}\cap\tra{C_j}\neq\emptyset,
(1\le i,j \le l)$
{\em generates\/} a unique
equivalence relation $\simeq$ such that,
if $\til{C}_1,\ldots,\til{C}_q$ are the
$\simeq$-equivalence classes,
then
$W_1=\cup\til{C}_1,\ldots,W_q=\cup\til{C}_q$
are pairwise separated open sets.
Moreover, these sets are interchanged by
any $f_t$ ($t\in T$) in the same way the
$C_1,\ldots,C_l$ were, that is,
the permutation of their indexes by the
action of $f_t$ does not depend on $t$.

\medskip

The permutation of the indexes of the
$W_1,\ldots,W_q$ has a finite number of
cycles which are a finite collection
of pairwise separated open sets
satisfying definition \ref{defdomsimetr}.
We have proved
\begin{proposition}\label{estac-dominv}
Every stationary measure $\mu$ is such
that the interior of its support is
made of a finite number of s-invariant domains.
\end{proposition}

\begin{remark}
If $\mu$ were ergodic,
then
$
\lim_{n\to\infty}
\frac1n \sum_{j=0}^{n-1}
1_{W_i} ( f_{\subl{t}}^j (x) )
= \mu ( W_i ) > 0
$
for $\mu\times\nuinf$-a.e.
$(x,\subl{t})\in M\times\De$
and $1\le i \le q$.
So almost every point
of $W_1\cup\ldots\cup W_q$
returns to $W_i$ infinitely many times.
In this case the {\bf interior of $\supp\mu$
is made of a single s-invariant domain\/}.
\end{remark}

\medskip

Let now $\SM_1,\ldots,\SM_h$ be {\em all\/}
the minimal domains inside the s-invariant
domains given by proposition \ref{estac-dominv}
(recall section \ref{minimalexist}).
Provisionally we assume the following

\begin{lemma}\label{estac-restric}
The normalized restriction of a
stationary measure to
a c-invariant set is a
stationary probability.
\end{lemma}

\medskip

Let the normalized restrictions be
$
\mu_{\SM_i}(A) = \frac1{\mu(\SM_i)} \cdot
\mu( A\cap \SM_i), \; i=1,\ldots,h,
$
where $A$ is any Borel set and $\mu(\SM_i)>0$
(because $\SM_i$ is a collection of open
sets inside $\inter{\supp\mu}$).
By proposition \ref{physical},
$\mu_{\SM_i}$ must be the characteristic
probability of $\SM_i$, $i=1,\ldots,h$.

\begin{remark}\label{bordo-0}
This means the characteristic probability
of each $\SM_i$ must give zero mass to
the border $\pa\SM_i$, since it
coincides with its normalized restriction
to the interior of $\SM_i$.
\end{remark}

To see that these probabilities are enough
to define $\mu$, we consider
$
\la=\mu - \mu(\SM_1)\cdot \mu_{\SM_1} -
\ldots - \mu(\SM_h)\cdot \mu_{\SM_h}.
$
If $\la\not\equiv 0$, then $\la$ is
a stationary measure (of course,
being stationary is an additive property)
whose support is nonempty.
By proposition \ref{estac-dominv}
and by section \ref{minimalexist}
we have some minimal domain $\SM$ in
$\supp\la$ with $\la(\SM)>0$.
But
$\supp\la\subset \supp\mu
\setminus
( \SM_1\cup\ldots\SM_h )$
and the $\SM_1,\ldots,\SM_h$ are
the only minimals in $\supp\mu$.
We have reached a contradiction,
so $\la\equiv 0$ and we have
proved (apart lemma \ref{estac-restric})

\begin{proposition}\label{linearcomb}
Every stationary probability is a
linear finite and convex combination
of characteristic probabilities.
\end{proposition}

\medskip

Let us note that these arguments show
that
$
\supp\mu=\supp \mu_{\SM_1} \dot{\cup}
\ldots \dot{\cup}\, \supp \mu_{\SM_h}
$
and consequently
$\mu(\SM_1)+\ldots+\mu(\SM_h) = 1$,
that is, the linear combination above
is indeed convex.



\medskip

To end this section we prove the remaining lemma.

\begin{rproof}{\ref{estac-restric}}
Let $\mu$ be a stationary measure and
$C$ a c-invariant set.

We remark that we know every point of $C$
stays in $C$, but we do not know whether
points in the complement $\supp\mu\setminus C$
enter in $C$ by the action of $f_t$.

First, we show
$D=\supp\mu\setminus C$ to be
{\em almost\/} completely invariant.
 
In fact, we may assume $\mu(D)>0$
(otherwise $C=\supp\mu$, $\mu$-mod $0$)
and write
\[
0< \mu(D) = \int 1_D(x) \, d\mu(x)
= \int \int 1_D( f(x,t) ) \, d\mu(x) \, d\nu(t)
\]
because $\mu$ is $S$-invariant. By the invariance of $C$,
$
x\in C \impl f(x,t)\in C \impl 1_D( f(x,t) ) = 0
$
for every $t\in T$
and so
$
\int \int 1_D( f(x,t) ) \, d\mu(x) \, d\nu(t)
= \int \int_D 1_D( f(x,t) ) \, d\mu(x) \, d\nu(t) .
$
 
Defining
$D_1(t)=\{ x\in D: f(x,t)\in D \}$ and
$D_2(t)=\{ x\in D: f(x,t)\not\in D \}$
for $t\in T$, we have
$
\mu(D) =
\int \int_{D_1(t)\cup D_2(t)}
1_D( f(x,t) ) \, d\mu(x) \, d\nu(t)
= \int \mu \left( D_1(t) \right) \, d\nu(t) > 0,
$
where
$\mu \left( D_1(t) \right) \le \mu (D)$
for every $t\in T$.
Thus $\mu \left( D_1(t) \right) = \mu (D)$
for $\nu$-a.e. $t$, that is,
$f(x,t)\in D$ for $\mu\times\nu$-a.e.
$(x,t)\in D\times T$.
In other words, points outside $C$ almost
never enter in $C$.
 
\medskip

Now we know that
$1_C(x) = 1_C( f(x,t))$
for $\mu\times\nu$-a.e. pair $(x,t)$.
Hence,
\[
\int \vfi\cdot 1_C \, d\mu =
\int \int \vfi( f_t(x) ) \cdot 1_C( f_t(x) )
\, d\mu(x) \, d\nu(t) =
\int \int \vfi( f_t(x) ) \cdot 1_C( x )
\, d\mu(x) \, d\nu(t)
\]
for any $\vfi\in C^0(M,\real)$, that is,
the restriction of $\mu$ to $C$ is stationary.

\end{rproof}

\section{Time Averages and Minimal Domains}
\label{conclusion}

What remains to be done is essentially to
fit together previous results.
Indeed, sections \ref{numerofinitodeminimais} and
\ref{ergodica} prove items $1$ and $2$ in the
statement of Theorem \ref{teorema1}.
To achieve the decomposition of item $3$ we are
going to show that every point $z\in M$ is
sent into some minimal domain by $\nuinf$-a.e.
perturbation of $\De$ and the $\nuinf$-mod $0$
partition of $\De$ obtained by this property
satisfies $3a$, $3b$ and $3c$, since we
already know that $m$-a.e. point inside a
minimal belongs to the respective ergodic
basin.

\medskip

Let $z\in M$ and let $\mu$ be a stationary probability
given by some accumulation point of the averages
\refer{average}.
By proposition \ref{linearcomb}
we know $\mu$ decomposes
in the following way
\begin{equation}\label{decomposes}
\mu = \al_1\cdot \mu_1 + \ldots + \al_h \cdot \mu_h
\end{equation}
where
$0 < \al_1,\ldots,\al_h < 1$,
$\al_1+\ldots+\al_h=1$ and
$\mu_1,\ldots,\mu_h$ are the characteristic
probabilities of the minimals
$\SM_1,\ldots,\SM_h$, respectively.

\medskip

Decomposition \refer{decomposes} and the
construction of $\mu$ ensure there is,
for every $i=1,\ldots,h$, a set
$V_i\subset\De$ with $\nuinf(V_i)>0$
such that
there is $k\in\natur$ satisfying
$f^k (z,\subl{s}) \in \SM_i$
for every $\subl{s}\in V_i$.

Indeed, $\mu( \SM_i ) >0$ implies there exist
open sets $U\subset\tra{U}\subset V \subseteq \SM_i$
such that $\mu(U)>0$ and so $\vfi\in C(M,\real)$
with $0\le \vfi \le 1$, $\supp\vfi\subseteq V$
and $\vfi_{|U}\equiv 1$ satisfies
$
\mu(\vfi)=
\lim\limits_{i\to\infty}
\frac{1}{n_i} \sum\limits_{j=1}^{n_i}
\int \vfi
\left(
f^j(z,\subl{t})
\right)  \,
d\nuinf (\subl{t})
> 0.
$
Then we have, for some $j\in\natur$:
\[
\nuinf
\left\{ \subl{t}\in\De: f^j (z,\subl{t})
\in \SM_i \right\}
\ge
\int \vfi    
\left(
f^j(z,\subl{t})
\right)  \,  
d\nuinf (\subl{t})
> 0.
\]
 
\medskip
 
Now we claim the sets $V_i$ occupy the entire space $\De$
or equivalently (cf. definition \ref{defGH})

\begin{proposition}\label{enche}
For every $z\in M$ we have
$
G_{\SM_1}(z) \cup \ldots \cup G_{\SM_l}(z) = \De,
\;\; \nuinf-\text{mod} \; 0$
and
$G_{\SM_i}(z) \cap G_{\SM_j}(z) = \emptyset$
for every pair $1\le i < j \le l$
where $\SM_1,\ldots,\SM_l$ are {\em all\/}
the minimal invariant domain of $\SD$.
\end{proposition}

\begin{proof}
By contradiction,
let us suppose there is $V\subset\De$ with
$\nuinf(V)>0$ such that
$\nuinf(V\cap G_{\SM_i}(z) ) =0$,
$i=1,\ldots,l$
(or
$V\subset\cap_{i=1}^l H_{\SM_i}(z), \; \nuinf-\text{mod}\; 0$).

\medskip

Let $\subl{t}$ be a $V$-generic vector and let
$w\in\om(z,\subl{t})$.
By lemma \ref{H=Delta} we have
$\cap_{i=1}^l H_{\SM_i}(w) = \De, \; \nuinf-\text{mod}\; 0$,
that is, the orbit of $w$
under almost every perturbation never falls in
$\SM_1\cup\ldots\cup\SM_l$.
Consequently any stationary probability
obtained from the orbits of $w$ as in section
\ref{abscont} will admit a (nontrivial) decomposition
(according to proposition \ref{linearcomb})
$
\mu = \be_1\cdot \til{\mu}_1 + \ldots
+ \be_{\til{h}} \cdot \til{\mu}_{h}
$
such that $0\le \be_1,\ldots,\be_{h} \le 1$,
$\be_1+\ldots+\be_{h}=1$ and
each $\til{\mu}_i$ is the characteristic
probability of $\til{\SM}_i$, $i=1,\ldots,h$, where
each of the
$\til{\SM}_1,\ldots,\til{\SM}_{h}$ is
distinct from $\SM_1,\ldots,\SM_l$.

\medskip

This contradict the supposition that the
$\SM_1,\ldots,\SM_l$ are {\em all\/} the
minimal invariant domains of $\SD$ and
so such a set $V$ cannot exist.
\end{proof}

\medskip

We now easily derive the continuous dependence
of the sets $V_i(x)$ from $x\in M$ with respect
to the distance between $\nuinf-\bmod 0$ sets
$A,B\subset\De$ given by
$d_{\nu}(A,B)=\nuinf(A\difsim B)$.

\medskip

We fix $x\in M$ and
note that each $V_i(x)$ can be
written as
\begin{equation}\label{decomp=kinM}
V_i(x) = \bigcup_{k=1}^\infty V_{i,k}(x)
\;\; \text{where}\;\;
V_{i,k}(x) = \{
\subl{t}\in\De : f^k(x,\subl{t})\in \SM_i
\},\;\; k\ge1,
\end{equation}
are open and
$V_{i,k}(x)\subseteq V_{i,k+1}(x)$
for all $k\ge1$
by the complete invariance of $\SM_i$,
$i=1,\ldots,l$.
This implies that for some $\de>0$ we can
find $k_0\in\natur$ such that
$
\nuinf( V_i(x)\setminus  V_{i,k_0}(x) )
\le \de
$
for all  $1\le i \le l$.

\medskip

On the one hand, by the finiteness of $k_0$,
property \ref{prop0} and the openness of
the domains that form $\SM_i$,
we get the existence
of $\ga>0$ with the property
$V_{i,k_0} (y)\supseteq V_{i,k_0}(x)$
for all $y\in B(x,\ga)$.
Hence
$
\nuinf(V_i(y)) \ge \nuinf( V_{i,k_0}(y) )
\ge \nuinf( V_{i,k_0}(x) ) \ge
\nuinf(V_i(x)) - \de
$
whenever $d_M(y,x) < \ga$ and
for every $i=1,\ldots, l$.

On the other hand
\begin{eqnarray*}
\nuinf( V_i(y) )
&=&
1 - \nuinf(V_1(y)) -\ldots - \nuinf(V_{i-1}(y))
- \nuinf(V_{i+1}(y)) - \ldots - \nuinf(V_h (y))
\\
&\le&
1 - \nuinf(V_1(x)) -\ldots - \nuinf(V_{i-1}(x))
- \nuinf(V_{i+1}(x)) - \ldots - \nuinf(V_h (x))
+
\\
& & + (h-1)\cdot \de
\\
&=&
\nuinf(V_i(x)) + (h-1)\cdot\de
\end{eqnarray*}
for all $1\le i \le l$ and continuity follows.

\medskip

We are left to show item $3c$ of theorem \ref{teorema1}
holds with respect to this decomposition.

Let us fix $1\le i \le l$ such that $\nuinf(V_i)>0$

We note that \refer{decomp=kinM}, the openness of
the $\SM_i$ and the continuity property \ref{prop0}(1)
imply the $V_i(z)$ to be open subsets of $\De$, that
is, for every $\subl{t}\in V_i(z)$ there are
$k\in\natur$ and $\ro>0$ such that
$
f^k( z , B(\subl{t},\ro) ) \subset \SM_i
\;\;\text{and so}\;\;
V_i(z)\supset B(\subl{t},\ro).
$
According to section~\ref{ergodicidade}
we have
\begin{equation}\label{mednul}
\SM_i \subset E(\mu_i) 
\;\;\text{and thus}\;\;
\left\{
\subl{s}\in V_i(z) : f^k( z,\subl{s} ) \in
E(\mu_i) 
\right\} \supset B(\subl{t},\rho).
\end{equation}
This means that every $\subl{s}$
in $B(\subl{t},\ro)\subset V_i=V_i(z)$
is such that $w=f^k(z,\subl{s}) \in E(\mu_i)$,
that is,
\[
\lim_{n\to\infty} \frac1n \sum_{j=0}^{n-1}
\vfi( f^j(w,\subl{u}) ) = \int \vfi \, d\mu_i
\;\; \text{for every} \;\;
\vfi\in C^0(M,\real)
\;\;\text{and}\;\;
\nuinf-\text{a.e.} \; \subl{u}\in\De.
\]


Since time averages do not depend on any
finite number of iterates,
item $3c$ of Theorem \ref{teorema1} follows
and the proof of Theorem \ref{teorema1} is complete.

\begin{remark}
\label{contopen}
We note that {\em diffeomorphism\/} in the arguments
and definitions of sections~\ref{notacoes}
through~\ref{conclusion} may be replaced
throughout by {\em continuous open map\/}.
This means Theorem~\ref{teorema1} is a
result of {\em continuous Ergodic Theory\/}
and not specific of {\em differentiable
Ergodic Theory\/}: a $C^0$-continuous and regular
family of continuous open maps
$f_t: M\mapto M$, $t\in B$, would suffice, i.e., for fixed
$x\in M$, $t\mapsto f(x,t)$ sends Lebesgue measure zero sets
into sets of m-measure zero.
\end{remark}

\begin{remark}
Hypothesis $B)$ of theorem \ref{teorema1}
was utilized in very specific points of
the proof, whereas hypothesis $A)$ was frequently
used throughout the arguments.

Naturally enough, condition $B)$ of
theorem \ref{teorema1} was used in the derivation
of the absolute continuity of stationary probability
measures for the perturbed system and in the proof
of lemma \ref{ergodiclaim}, that is, in the
proof of ergodicity for a stationary measure
$\mu$ supported in some minimal invariant
domain $\SM$, $\supp\mu\subset\tra{\SM}$.

This was the sole role of
hypothesis $B)$ in the proof of theorem \ref{teorema1}.
\end{remark}

\begin{remark}
The probability $\nuinf$ defined on $\De$ was
chosen for simplicity.
The property of $\nuinf$ used in the proof of
theorem \ref{teorema1} was, besides condition $B)$,
that every open set of $\De$ has positive
$\nuinf$-measure.
This implies that
\begin{enumerate}
\item any set $Y\subset\De$ with $\nuinf(Y)=1$
is a dense set, $\tra{Y}=\De$.
\end{enumerate}
Moreover, because $\nuinf$ is an uniform measure
\begin{enumerate}
\item[2.] for every $\de>0$ there is $\ga>0$
such that each set $Y\subset\De$ with
$\nuinf(Y)\ge 1-\de$ is $\ga$-dense,
$\bigcup\limits_{y\in Y} B(y,\ga) \supset \De$.
\end{enumerate}

These properties were of use in the proof
of the fundamental lemmas ( number 2 above)
in connection with the general property
(0-1 type law) given by lemma \ref{lema0-1}
and again (now number 1 above) in the arguments
of the proof of lemma \ref{ergodiclaim}.

Any other probability satisfying numbers 1 and
2 above would do for the statement and proof
of theorem \ref{teorema1}.
\end{remark}

\begin{remark}\label{weakeningAB}
The conclusions of Theorem \ref{teorema1}
can be obtained with weaker hypothesis
instead of the stated $A)$ and $B)$.

Indeed, it is very easy to see that the integer
$N$ may depend on $x$ in the statement of $A)$.
Thus it can be replaced by
\begin{enumerate}
\item[A')] There is $\xi_0>0$ such that for
all $x\in M$ there exists $N=N(x)\in\natur$
satisfying
$f^k(x,\De)\supset B(f^k(x),\xi_0)$
for all $k\ge N$.
\end{enumerate}
Moreover, $B)$ can be weakened so that the absolute
continuity of a stationary probability $\mu$
still holds by allowing $f^k(x,\nuinf) \ll m$
for {\em some\/} $k\ge1$.
If this $k$ does not depend on $x\in M$, then
we can still prove Proposition~\ref{physical}
in the same way.

Other weakenings of $B)$ are possible,
one such will be of use following
section \ref{bifurcations} dealing with
random parametric perturbations near
homoclinic bifurcations.
\end{remark}

\section{Bowen's Example}
\label{Bowen}

This is the answer to a question raised by C. Bonatti.
This example captures the meaning of Theorem
\ref{teorema1}:
even if a given {\em deterministic\/}
({\em noiseless\/}) system is devoid
of physical measures (its Birkhoff averages do
not exist almost everywhere) we may nevertheless
get a finite number of physical probabilities
describing the asymptotics of almost every orbit
just by {\em adding\/} a small amount of
random noise.

\begin{example}\label{exemploBowen}
{\em Bowen's example\/}
(see \cite{T2} for the not very clear
reason for the name) is a folklore example
showing that Birkhoff averages need not exist
almost everywhere.
Indeed, in the system pictured in figure \ref{Bowensys}
Birkhoff averages for the flow {\em do not exist
almost everywhere\/}, they only exist for the sources
$s_3$, $s_4$ and for the set of separatrixes and saddle
equilibria
$W=W_1\cup W_2 \cup W_3 \cup W_4 \cup \{s_1,s_2\}$.

\begin{figure}[h]
 
\centerline{\psfig{figure=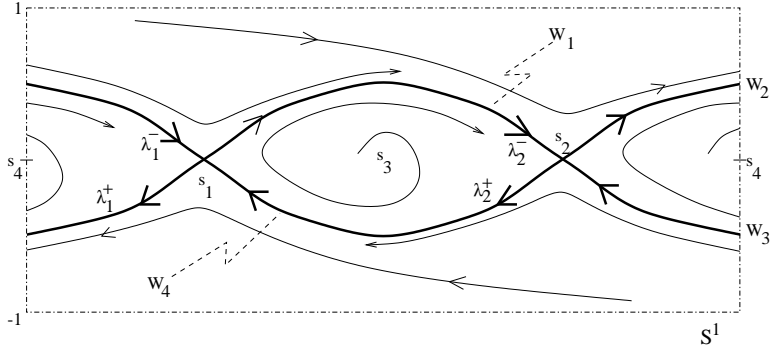,width=4.0in}}
 
\caption{\label{Bowensys} A sketch of Bowen's example flow.}
\end{figure}

The orbit under this flow $\phi_t$ of every point
$z\in S^1\times [-1,1] = M$ not in $W$
accumulates on either
side of the separatrixes, as suggested in the figure,
if we impose the condition
$\la_1^{-} \la_2^{-} > \la_1^{+} \la_2^{+}$
on the eigenvalues of the saddle fixed points
$s_1$ and $s_2$ (for more specifics on this see
\cite{T2} and references therein).

\medskip

We apply Theorem \ref{teorema1} to this case.
We remark that $M$ is not a boundaryless manifold,
but its border
$S^1\times\{ \pm 1 \}$
is sent by $\phi_1$ into $S^1\times [-1,1]$.
Moreover, Theorem \ref{teorema1} refers not
to perturbations of flows, so we will
consider the time one map $\phi_1$ as our
diffeomorphism $f:M\mapto M$ and, since
$M$ is parallelizable, we can make an
absolutely continuous random perturbation,
as in example~1 of section~\ref{teoremaprincipal}.
In this circumstances the proof of Theorem
\ref{teorema1} equally applies.

For everything to be properly defined, though,
we must restrict the noise level $\ep>0$ to
a small interval $]0,\ep_0[$ such that the
perturbed orbits stay in $S^1\times ]-1,1[$.
After this minor technicalities we proceed to
prove

\begin{proposition}\label{Bowenfisico}
The system above, under random absolutely
continuous noise of level $\ep\in ]0,\ep_0[$,
admits a single physical absolutely continuous
probability measure $\mu$ whose support is a
neighborhood of the separatrixes:
$\inter{\supp\mu} \supset W$.
Moreover the ergodic basin of $\mu$ is the
entire manifold:
$E(\mu)= M , \; \mu\bmod 0$.
\end{proposition}

\begin{proof}
Let $\ep\in ]0,\ep_0[$ be the fixed noise level
from now on and let $U$ be the ball of radius
$\ep/4$ around $s_1$. We will build fundamental
domains for the action of $f=\phi_1$ over
$M\setminus W$ in $U$, as explained below.

\begin{figure}[h]
 
\centerline{\psfig{figure=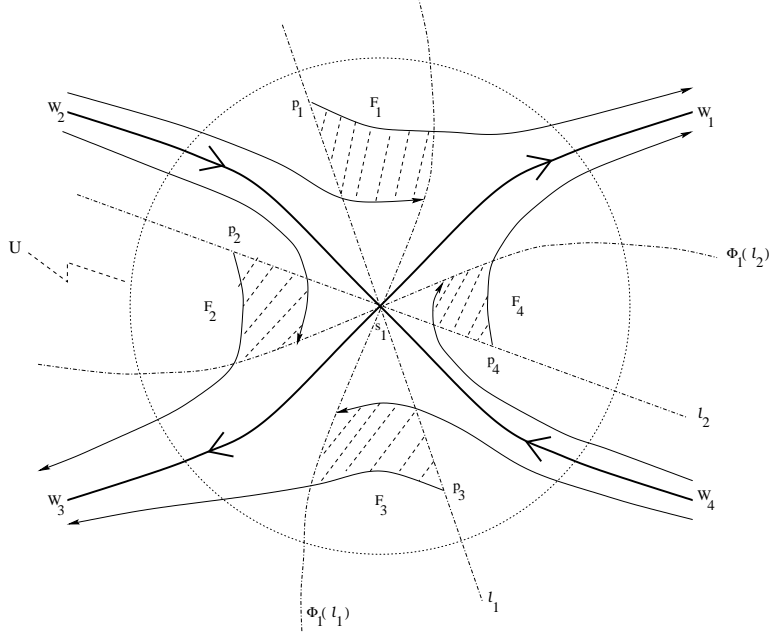,width=4.0in}}
 
\caption{\label{grandeplanosela}
How the fundamental domains are obtained.}
\end{figure}

We choose two strait lines $l_1,l_2$ through
$s_1$ crossing $U$ and let $l_1',l_2'$ be their
images under $\phi_1$ as sketched in figure
\ref{grandeplanosela}.
Now we choose two points in each line $l_1,l_2$
on either side of $s_1$: $p_1, p_2, p_3$ and $p_4$,
and consider their orbits under the flow $\phi$
for positive time, until they return to $U$ and
cut $l_1',l_2'$, as depicted in the abovementioned
figure.

\medskip

The four intersections of the orbit of $p_i$ with the
proper $l_j,l_j'$, together with portions of the
orbit and of $l_j,l_j'$ define a "square"
$F_i$ (shadowed in figure \ref{grandeplanosela})
which is a fundamental domain for the dynamics of
$f=\phi_1$ on the connected components of
$M\setminus W$, $i=1,2,3,4$ and $j=1$ or $2$.

\medskip

This means that every $z\in M\setminus W$ is
such that there is a $k\ge1$ with
$z_k=f^k(z)\in F = F_1\cup F_2 \cup F_3 \cup F_4$.
Moreover, by the choice of $U$, $z_k$ may
be sent into any $F_i$, $i=1,2,3,4$, by adding to
a vector of length smaller than $\ep$.
Thus we deduce that
$f^k(z,\De) \supset F_1\cup F_2 \cup F_3 \cup F_4$
and even more: $f^k(z,\De) \supset U$.

Keeping in mind that for $m\ge1$ we have
$
f^{k+m}(z,\De) = f^m \left( f^k(z.\De) , \De \right)
= \{ f^m(w,\De) : w\in f^k(z,\De) \},
$
we see that $f^{k+m}(z,\De)$ will contain all the
$f$-images of each $F_1,F_2,F_3$ and $F_4$,
which will return to $U$ infinitely many times.
Furthermore, at each return the points may again
be sent into any $F_1,F_2,F_3$ or $F_4$ by
an $\ep$-perturbation.
Hence the sets of the sequence
$\{ f^n(z,\De) \}_{n=1}^\infty$
contain $F_1,F_2,F_3$ or $F_4$ for
infinitely many $n's$ and also all their
$f$-images.

\medskip

We conclude that $\om(z,\De)$ contains a
neighborhood of $W$.

The same holds for $w\in W$, since $f^1(w,\De)$
is an open set and so contains some $z\in M\setminus W$.
That is,
{\em every $z\in M$ is such that $\om(z,\De)$ contains
a neighborhood of $W$\/}.

\medskip

Therefore, there can be only one minimal $\SM$ in
the perturbed system, such that $\SM\supset W$ and
into which every point $z\in M$
finally falls by almost every perturbed orbit
(this is a consequence of sections
\ref{minimais-disjuntos}, \ref{ergodica},
\ref{decomposition} and \ref{conclusion}).
We have further that the
characteristic probability $\mu_{\SM}$ is
{\em the physical probability\/} $\mu$
of the system, with $E(\mu)=M,\; \mu\bmod 0$,
and $\supp\mu\supset\SM\supset W$,
as stated.
\end{proof}
\end{example}

\section{Homoclinic Bifurcations \\
and Random Parametric Perturbations}
\label{bifurcations}

We consider arcs (one-parameter families) of
diffeomorphisms exhibiting a 
quadratic homoclinic tangency and
derive similar properties for their random parametric
perturbations to those stated in Theorem \ref{teorema1}.

\subsection{One-Parameter Families}
\label{1parfam}

The arcs we will be considering are given by
a $C^\infty$ function
$
f: M^2 \times ]-1,1[ \mapto M^2
$
such that for every $-1<t<1$,
$f_t:M^2\mapto M^2, \; x\mapsto f(x,t)$
is a diffeomorphism of the boundaryless
surface $M^2$.
The family of diffeomorphisms
$\SF=(f_t)_{-1<t<1}$ satisfies the following
conditions.

\begin{figure}[h]
 
\centerline{\psfig{figure=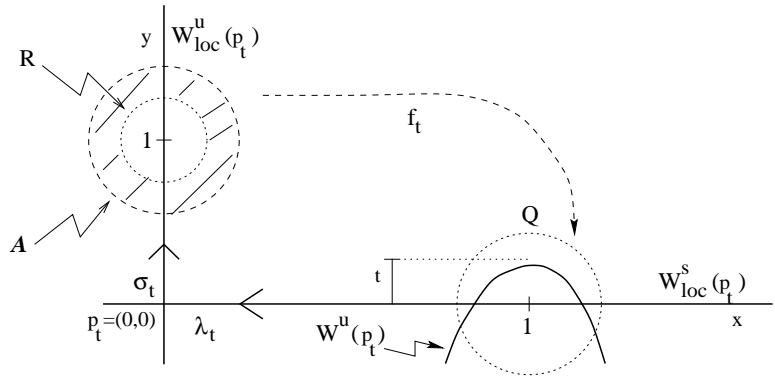,width=4.0in}}

\caption{\label{tangencia} A sketch of the situation to be considered}
\end{figure}

\begin{enumerate}
\item $\SF$ has a {\em first tangency\/} at $t=0$,
that is (v. \cite[Appendix 5]{PT})     \label{1tang}

\begin{enumerate}
\item for $t<0$, $f_t$ is persistently hyperbolic; \label{pershyp}

\item for $t=0$ the nonwandering set $\Om(f_0)$ consists of
a closed hyperbolic set
$\til\Om(f_0)=\lim_{t\nearrow 0} \Om(f_t)$
together with a homoclinic orbit of tangency
$\SO$ associated with a hyperbolic fixed saddle point $p_0$,
so that $\Om(f_0)=\til\Om(f_0)\cup\SO$;       \label{cnerrante}

\item the branches $W_+^s(p_0)$, $W_+^u(p_0)$ of the invariant
manifolds $W^s(p_0)$, $W^u(p_0)$ have a quadratic tangency along
$\SO$ unfolding generically as pictured in figure~\ref{tangencia}
(v. \cite[Capt. 3]{PT}):
$\SO$ is the only orbit of tangency between stable and unstable
separatrixes of periodic orbits of $f_0$;       \label{quadtang}
\end{enumerate}



\item The saddle $p_0$ has eigenvalues $0<\la_0<1<\si_0$
satisfying the conditions for
the existence of $C^2$ linearizing coordinates in
a neighborhood of $(p_0,0)$ in $M^2\times]-1,1[$
(v. \cite{T1}).   			\label{nressona}
\end{enumerate}


Condition~\ref{1tang} imposes bounds on the region
where new accumulation points can appear for $t>0$
(small) --- section~\ref{adaptinglinear} will specify
this (cf.~\cite[Appendix 5]{PT}).


We note that condition~\ref{nressona} above is generic
in the space of all $C^\infty$ one-parameter families
satisfying \ref{1tang}.
Moreover, those families that satisfy \ref{1tang}
are open
(cf. \cite[Capt. 3, Appendix 5]{PT} and references therein).

\subsection{Statement of The Results}
\label{resultados2}

For some small $t^\star>0$, to be explained in the
following sections,
we fix $t_0\in]0,t^\star[$,
$\ep_0 = \min\{ | t_0 |, |t^\star - t_0 | \}$
and the {\em noise level\/} $\ep\in]0,\ep_0[$.
We consider the system $f_{t_0}$ under a random
parametric perturbation of noise level $\ep$,
$\SF_{t_0,\ep}$, as defined in section~\ref{pertaroundpar}.
We let $\De=\De_{\ep}(t_0)$ be the perturbation space
$[t_0-\ep,t_0+\ep]^\natur$.

\medskip

We will be interested in studying what happens in
$\SQ$, a closed neighborhood of $q$ to be constructed.
We need an effective definition of {\em interesting points\/}.

\begin{definition}[First Return Times]
\label{deft1ret}
Given some $z\in M^2$ and $\subl{t}\in\De$ we let
\[
r(z,\subl{t},1)=\min\{ k\ge 0 : f^k(z,\subl{t})\in \SQ \}
\]
and inductively define
$
r(z,\subl{t},n+1)=\min\{ k \ge 1 :
f^{R(z,\subl{t},n)+k}(z,\subl{t})\in\SQ \}
$
for every $n\ge1$,
where $R(z,\subl{t},n) = \sum_{i=1}^n r(z,\subl{t},i)$,
with the convention $\min\emptyset=+\infty$.
\end{definition}

\begin{definition}\label{def:Vrecorrente}
A {\bf $V$-recurrent point\/} is a $z\in\SQ$ for which
there exists a $V\subset\De$ satisfying
\begin{enumerate}
\item $\nuinf(V)>0$;
\item $r(z,\subl{t},n)<\infty$ for every $n\ge1$
and $\nuinf$-a.e. $\subl{t}\in V$.
\end{enumerate}
\end{definition}

In other words, $z\in\SQ$ is {\em interesting\/} if its perturbed
orbits pass through $\SQ$ infinitely often under a
positive measure set of perturbations.

\medskip

We can now state

\begin{maintheorem}\label{teorema2}
For every $C^\infty$ arc of diffeomorphisms as described
in subsection~\ref{1parfam} and any given homoclinic tangency
point $q$ associated to the saddle $p_0$, there are a
closed neighborhood $\SQ$ of $q$ and $t^\star>0$
such that, for each $t_0,\ep>0$ satisfying
$0< t_0 <t^\star$ and
$0<\ep<\ep_0=\min\{|t_0|,|t^\star - t_0|\}$,
the random parametric perturbation $\SF_{t_0,\ep}$ of $f_{t_0}$
with noise level $\ep$ admits a finite number of
probabilities $\mu_1,\ldots,\mu_l$ whose support
intersects $\SQ$ and that
\begin{enumerate}
\item $\mu_1,\ldots,\mu_l$ are physical absolutely continuous
probability measures;
\item $\supp{\mu_i}\cap\supp{\mu_j}=\emptyset$
for all $1\le i < j \le l$;       
\item for all $z\in \SQ$ and $V\subset\De$ such that
$z$ is $V$-recurrent
there are open sets $V_1=V_1(z),\ldots,V_l=V_l(z)\subset V$
such that
\begin{enumerate}
\item $V_i \cap V_j = \emptyset,\;\; 1\le i < j \le l$;
\item $\nuinf( V \setminus (V_1\cup\ldots\cup V_l)) =0$;
\item for all $1\le i\le l$ and $\nuinf$-a.e. $\subl{t}\in V_i$
we have
\[
\lim_{n\to\infty} \frac1n \sum_{j=0}^{n-1}
\vfi( f^j(z,\subl{t}) ) = \int \vfi \, d\mu_i, \;\;
\text{for every} \; \vfi\in C(M,\real)  .
\]

\end{enumerate}
\end{enumerate}
\end{maintheorem}

\subsection{Adapting the Linearization}
\label{adaptinglinear}

As preparation for the proof of Theorem~\ref{teorema2}
by using Theorem~\ref{teorema1}
we study the adaptation of the linearizing coordinates
to our setting.

Condition~\ref{nressona} enables us to consider a change of
coordinates $\vfi_t:L\subset\real^2\mapto M^2$
in a neighborhood $L$ of every $p_t$,
where $|t|<t^\star$ for some small  $t^\star>0$
and
\begin{equation}\label{linearconjugacy}
f_t(\vfi_t(x,y))=\vfi_t( \la_t\cdot x, \si_t\cdot y)
\end{equation}
with $0<\la_t < 1 < \si_t$ the
eigenvalues of the hyperbolic saddle fixed point $p_t$.
These coordinates will be adapted much like
\cite[p.49 and Appendix 5]{PT}.
Specifically, after choosing
a homoclinic point $q$ associated to $p_0$:

\begin{enumerate}
\item[I)] we suppose $q\in W^u(p_0)\cap W^s(p_0)$
to be in $L$ --- to achieve this we may extend $L$
along $W^s(p_0)$ as explained in~\cite[Capt.2]{PT};

\item[II)] we extend $L$ along $W^u(p_0)$ in order that
$r=f_0^{-1}(q)$ be in $L$;

\item[III)] we use the implicit function theorem and
two independent rescalings of the x- and y-axis to get,
because of condition~\ref{quadtang}:

\begin{enumerate}
\item $q=(1,0)$, $r=(0,1)$, $p_t=(0,0)$ and
$W^s_{\rm loc}(p_t), W^u_{\rm loc}(p_t)$
are the x- and y-axis, respectively;

\item $f_t(0,1)$ is a local maximum of the y-coordinate
restricted to $W^u(p_t)$;

\item $\vfi_t^{-1}\circ f_t \circ\vfi_t(0,1)=(1,t)$;
\end{enumerate}
for every $|t|<t^\star$ in the coordinates defined
by $\vfi_t$;

\item[IV)] writing $\La_0$ the basic
set to which $p_0$ belongs
(possibly $\La_0=\{p_0\}$ trivially)
by condition~\ref{cnerrante} we have
$W^s(\La_0)=W^s(\La_0\cup\SO)$ and
$W^u(\La_0)=W^u(\La_0\cup\SO)$ and there exists
a filtration
$\emptyset\neq M_1 \subset M_2 \subset M$ such that
(v.~\cite[Appendix 5, pp. 212-214]{PT}
and cf.~\cite[Capt. 1]{S})

\begin{enumerate}
\item $M_i$ is closed and $f_0(M_i)\subset \inter{M_i}$
for $i=1,2$;
\item $M_1\subset\inter{M_2}$, and
\item $\La_0\cup\SO=\left(
\cap_{j\ge0} f_0^j(M_2) \right) \cap
\left( \cap_{j\ge0} f_0^{-j}(M_1^c)  \right)$;
\end{enumerate}

\item[V)] since $\La_0$ is a basic set (of saddle type)
there is a small compact
neighborhood $U$ of $\La_0$ where extensions
$\SH^s$, $\SH^u$ of the stable and unstable foliations
$W^s(\La_0)$, $W^u(\La_0)$ are defined
(v.~\cite[Appendix 1]{PT} and references therein),
and by IVc
there is $N^\star\in\natur$ such that

\begin{enumerate}
\item $\left( \cap_{j=0}^{N^\star-1} f_0^j(M_2) \right) \cap
\left( \cap_{j=0}^{N^\star-1} f_0^{-j}( M_1^c ) \right)
\subset U \cup \SQ^\star$
where $\SQ^\star\subset L$ is a neighborhood of
the portion of $\SO$ outside $U$
with finitely many components
$\SQ_1,\SQ_2,\ldots,\SQ_l$ and
$\SQ^\star\cap U=\emptyset$.
Moreover we can assume they satisfy
$f_0(\SQ_1)\subset\SQ_2, \ldots, f_0(\SQ_{l-1})\subset\SQ_l$
with $q\in\til{\SQ}=\SQ_i$, $i\in\{1,\ldots,l\}$;

\item making $t^\star>0$ smaller if need be
and $\SQ^\star$ and $U$ a little bigger, we get also for all
$t_1,\ldots,t_N, t_1',\ldots,t_N'\in]-t^\star,t^\star[$
\[
\left( \bigcap_{j=0}^{N^\star-1} f_{t_j}\circ\ldots\circ f_{t_1} (M_2)
\right) \cap \left(
\bigcap_{j=0}^{N^\star-1} f^{-1}_{t_1'}\circ\ldots\circ f^{-1}_{t_N'}
( M_1^c )\right)  = U \cup \SQ^\star
\]
and also
$f_{t_1}(\SQ_1)\subset\SQ_2,\ldots,f_{t_{l-1}}(\SQ_{l-1})\subset\SQ_l$
for all $t_1,\ldots,t_{l-1} \in ]-t^\star,t^\star[$;
\item $f_t(M_i) \subset \inter{M_i}$ for all $|t|<t^\star$
and $i=1,2$;
\item $\La_t = \cap_{n\in\relativ} f^n_t( U )$
is the analytic continuation of $\La_0$ for all $|t|<t^\star$;
\end{enumerate}

\item[VI)] for every closed neighborhood
$\SQ\subset\til{\SQ}\subset L$ of $q$ and $t^\star>0$
small we have that

\begin{enumerate}
\item there is $N_Q \in \natur$ such that
$f_{t_i}\circ\ldots\circ f_{t_1}( \SQ ) \subset L$
for all $t_1,\ldots,t_i \in ]-t^\star,t^\star[$ and
$i=1,\ldots,N_Q$;

\item in the neighborhood
$\SR =\tra{ \cup_{|t|<t^\star} f^{-1}_t(\SQ) }$
of $r=(0,1)$ --- we may suppose $\SR \subset L$
by  making $\SQ$ and $t^\star$ smaller,
keeping $(a)$ by increasing $N_Q$
--- the map
$\til{f}_t=\vfi_t^{-1}\circ f_t \circ\vfi_t$
has the form
\begin{equation}\label{formaemR}
(x,1+y) \mapsto
(1+\al y+\eta x + H_1(t,x,y)\; ;\; \be y^2 + \ga x + t +H_2(t,x,y) )
\end{equation}
where $\al \cdot \be \cdot \ga \neq 0$,
$H_1$ is of order 2 or higher and $H_2$ is of
order 3 or higher in $y$ and order 2 or higher in
$x,t$ and $y\cdot t$;


\item for all $|t|<t^\star$ we make $f_t(0,1)\in \inter{\SQ}$
by taking $t^\star$ smaller if needed and keeping $\SQ$
and $N_Q$ unchanged satisfying $(a)$ and
(re)defining $\SR$ as in $(b)$.

\item for any given $\de_0>0$ and all sufficiently small
$\SQ$ and $t^\star$, we may keep
everything up until now increasing $N_Q$ and imposing
$
| D_2 H_i | , | D_3 H_i | < \de_0,\; i=1,2;
$
\end{enumerate}


\item[VII)] since all of the above holds for every
small (compact) neighborhood $\SQ\subset\til{\SQ}$
of $q$ and $t^\star>0$, except that $N_Q$ increases,
we may suppose $\SQ$ is so small that $N_Q> N^\star$
and then make $t^\star$
so small that item V) holds with $\SQ$ in the place
of $\til{\SQ}$ for some integer $N>N^\star$.
Furthermore writing $\SQ'$ for this new
neighborhood we may suppose that
$\La_t$ still is the maximal invariant set
inside $B( U,\rho ) = \cup_{z\in U} B(z,\rho)$
for $|t|<t^\star$ and
$\tra{ B( U,\rho ) } \cap \tra{ B(\SQ',\rho) } = \emptyset$
for some small $\rho>0$;

\item[VIII)] we may suppose the extended foliations
$\SH^s,\SH^u$, which are defined in a neighborhood of $p_0$
(since $p_o\in\La_0$), were extended by positive and
negative iterations of $f_0$ to cover all of $L$.
Moreover we may assume also that there are extended
foliations $\SH^s_t,\SH^u_t$ defined all over $L$
with respect to $f_t$ for every $|t|<t^\star$;

\item[IX)] in a small neighborhood $\SA$ of $\SR$ given
by
$\SA= \left( \cup_{z\in\SR} B(z,\xi) \right) \setminus \SR$,
$\xi>0$ small
(we may think of it as a small
{\em annulus\/} around $\SR$),
every point is sent by $f_t$ outside of $U\cup\SQ'$,
for every $t\in T$,
because $U$ and $\SQ'$ are separated according to
item VII. $\SA$ is open and will be called
the {\em nonreturn\/} annulus.
\end{enumerate}

We note that figure~\ref{tangencia} was made having these
items already in mind.

\subsection{Another Tour of Another Proof}
\label{tour2}

To begin with, pick a $V$-recurrent point $z\in\SQ$
and deal with its generic $\om$-limit points $w$,
which are always {\em regular\/} by the following

\begin{proposition}\label{regular}
There exists $J\in\natur$ such that if $z\in\SQ$ is $V$-recurrent
for some $V\subset\De=\De_\ep(t_0)$
with $\nuinf(V)>0$, then
the first return times of $w\in\om(z,\subl{t})$,
for all $V$-generic $\subl{t}$, do not depend on
$\subl{s}\in\De$ and are bounded by $J$:
\[
r(w,\subl{s},n)\equiv r(w,n) \le J \;\;
\text{for every}\;\; n\ge1.
\]
\end{proposition}

\begin{definition}\label{def:regular}
The points $w\in M^2$ which satisfy the conclusion of
the proposition
above will be called {\bf regular points}
(with respect to $\SF_{t_0,\ep}$).
\end{definition}

Taking advantage of the regularity of $w$,
the expression~\refer{formaemR} for $f_t|_{\SR}$ and
condition~\ref{1tang},
we will derive versions of hypothesis $A)$ and $B)$
of Theorem \ref{teorema1}:

\begin{proposition}\label{fisica}
Let $w\in M^2$ be a regular point.
Writing $r_n=r(w,n)$, $n\ge1$, the following holds.
\begin{enumerate}
\item For every $s\in ]t_0-\ep, t_0+\ep[$
there is a $\xi_0=\xi_0(s)>0$ such that
for all $n\ge2$
\[
f^{R_n}(w,\De) \supset B( f^{R_n}_s (w), \xi_0)
\;\;\text{where}\;\;
R_n=\sum_{i=1}^n r_i ;
\]

\item For all $n\ge2$ it holds that $f^{R_n}(w,\nuinf)\ll m$.
\end{enumerate}
\end{proposition}

In other words, we get conditions $A)$ and $B)$ of Theorem
\ref{teorema1} for the return times of $w$, which do not
depend on the perturbation chosen, since $w$ is regular.
Behind proposition~\ref{fisica} is the geometrically
intuitive idea of mixing expanding and contracting
directions near $q$ due to the
homoclinic tangency, together with condition~\ref{1tang}
that keeps the orbits of regular points confined in a
neighborhood of $\La_0\cup\SO$ (v. section~\ref{fisica:prova}).

\medskip

This is enough to prove Theorem~\ref{teorema2}.

Indeed, setting $K=2(J+1)$ then $R_2=R_2(w)\le K$
for every regular point $w$ and for $k\ge K$ there are
$n\ge2$ and $0\le i \le r_{n+1}-1\le J$
(by proposition~\ref{regular}) such that
$k=R_n+i$.
After item 1 of proposition~\ref{fisica} we have
$f^k(w,\De)= f^i( f^{R_n}(w,\De), \De) \supset
f^i_{t_0}( B( f^{R_n}_{t_0}(w), \xi_0) )$
and since $0\le i \le J$ there is some $\xi_0'>0$
such that
$f^i_{t_0}( B( f^{R_n}_{t_0}(w), \xi_0) ) \supset
B( f^{R_n+i}_{t_0} (w), \xi_0') = B( f^k_{t_0}(w), \xi_0')$
because $f_t$ is a diffeomorphism. We have hypothesis $A)$.

For hypothesis $B)$ we let $w$ and $k\ge K$ be as above.
Then $k=R+i$ with $i\ge0$ and $R=R_2=R_2(w)$.
We suppose $i\ge1$ for otherwise item 2 of
proposition~\ref{fisica} does the job.
We take a measurable set $E\subset M^2$ such that
$m(E)=0$ and observe that
$f^{R+i}(w,\nuinf)E = \nu^k(F)$
where
$F=\{ (t_1,\ldots,t_k)\in T^k:
f^{R+i}(w,t_1,\ldots,t_k) \in E \}$.

Defining for every $(t_{R+1},\ldots,t_k)\in T^i$
the section
$
F(t_{R+1},\ldots,t_k)=
\{ (s_1,\ldots,s_R)\in T^R:
(s_1,\ldots,s_R,t_{R+1},\ldots,t_k)\in F \}
$
we have by Fubini's theorem
\begin{equation}\label{mednulseccao}
\nu^k(F) = \nu^{R+i}(F) =
\int \nu^R( F(t_{R+1},\ldots,t_k) ) \,
d\nu^i(t_{R+1},\ldots,t_k).
\end{equation}
However
$
F(t_{R+1},\ldots,t_k) = \{
(s_1,\ldots,s_R)\in T^R:
f^i_{t_{R+1},\ldots,t_k} \circ
f^R_{s_1,\ldots,s_R} (w) \in E
\}
=
\{
(s_1,\ldots,s_R)
: f^R_{s_1,\ldots,s_R} (w) \in
\left( f^i_{t_{R+1},\ldots,t_k}  \right)^{-1} (E)
\}
$
and each $f_t$ is a diffeomorphism, so the
inverse image of a set of measure zero is a
set of measure zero.
Hence
$
\nu^R( F(t_{R+1},\ldots,t_k) )
$
is given by
$
f^R(w,\nuinf)
\left[
\left( f^i_{t_{R+1},\ldots,t_k}  \right)^{-1} (E)
\right] = 0
$
since $f^{R_2}(w,\nuinf) \ll m$ by proposition~\ref{fisica}(2).
We deduce from~\refer{mednulseccao} that
$
f^{R+i}(w,\nuinf)E = f^k(w,\nuinf)(E) = \nu^k(F) = 0
$
whenever $m(E)=0$, i.e.,
$f^k(w,\nuinf)\ll m$ for every $k\ge K$.

\medskip

It is clear that Theorem~\ref{teorema2} holds
by considering $(\SD,\preceq)$ as the set of
s-invariant domains $D=(\SU_0,\ldots,\SU_{r-1})$
with respect to $\SF_{t_0,\ep}$ whose points
$\SU_0 \cup \ldots \cup \SU_{r-1}$ are
regular points, with the same relation $\preceq$
as before, and using Theorem~\ref{teorema1}.

We should explain how to get the decomposition of
item 3 of Theorem~\ref{teorema2} for $V$-recurrent
point $z\in\SQ$. We use two previous ideas:
\begin{nliste}{\Parentalgarismo}
\item Going back to section~\ref{conclusion},
taking a generic $w\in\om(z,\subl{t})$
(i.e., $\subl{t}$ is $V$-generic) provides a
stationary probability $\mu$, as in sections~\ref{ergodica}
and~\ref{decomposition}, which decomposes as
in~\refer{decomposes} and we get the sets
$
V_i=\{ \subl{s}\in\De: \exists k\ge1
\;\;\text{s.t.} \;\;
f^k(w,\subl{s})\in\SM_i \}
$
as in item 3 of Theorem~\ref{teorema1}.

\item The previous item
together with proposition~\ref{regular}
just says that a $V$-recurrent
point $z\in\SQ$ satisfies lemma~\ref{finalmentecai},
i.e.,
there are $W\subset V$ with $\nuinf(W)>0$ and $m\in\natur$
such that $f^m_{\subl{\te}}(z)\in\SM$ for every
$\subl{\te}\in W$, where $\SM$ is some minimal
of $(\SD,\preceq)$.
We know there is just a finite number
$\SM_1,\ldots,\SM_l$ of minimals in $(\SD,\preceq)$
and define
$
V_i=V_i(z)=\{ \subl{s}\in V :
\exists k\ge1 \;\;\text{s.t.}\;\; f^k(z,\subl{s})\in \SM_i \},
$
$1\le i \le l$.
Repeating the arguments of proposition~\ref{enche}
with $\De$ replaced by $V$ throughout
gives item 3 of Theorem~\ref{teorema2}
and completes the proof.
\end{nliste}

\section{{\em Physical\/} Parametric Noise
with a Single Parameter}
\label{fisica:prova}

We start the proof of proposition \ref{fisica}
deducing the following consequence of condition~\ref{1tang}
in section~\ref{1parfam} and items V and VI.

\begin{lemma}\label{confina}
For every small $t^\star>0$ and $\SQ$
and every $z\in\SQ$ recurrent under
some vector $\subl{t}=(t_j)_{j=1}^\infty$ with $|t_j|<t^\star$,
$j\ge1$, i.e., such that $\om(z,\subl{t})\cap \SQ \neq\emptyset$,
the following holds
\begin{equation}\label{confinado}
f^j(z,\subl{t}) \in L \;\;
\text{for} \;\; 0\le j \le N_Q \;\;
\text{and} \;\; f^j(z,\subl{t}) \in U\cup \SQ^\star \;\;
\text{for} \;\; j \ge N_Q.
\end{equation}
\end{lemma}

\begin{proof}
We let $z\in\SQ\subset U\cup\SQ^\star \subset M_2\cap M_1^c$
be a recurrent
point under $\subl{t}$ as stated, and suppose that
$f^j(z,\subl{t})\not\in U\cup \SQ^\star$ for some
$j\ge N^\star$.
Then by item Vb it must hold
\[
f^j(z,\subl{t})\in\bigcup_{i=0}^{N^\star-1} f_{t_j}\circ\ldots\circ
f_{t_{j-i}}(M_2^c) \;\; \text{or} \;\;
f^j(z,\subl{t})\in\bigcup_{i=0}^{N^\star-1}
f^{-1}_{t_{j+1}}\circ\ldots\circ f^{-1}_{t_{j+i}}(M_1).
\]

Since $z\in U\cup\SQ^\star \subset M_2$ we have by
item Vc that $f^i(z,\subl{t}) \in M_2$ for every $i\ge 0$.
Hence only the right hand side alternative above can hold,
otherwise we would have for some $0\le i \le {N^\star-1}$ that
$
f_{t_j}\circ\ldots\circ f_{t_1}(z) \in
f_{t_j}\circ\ldots\circ f_{t_{j-i}} (M_2^c)
$
and so
$
f_{t_{j-i-1}}\circ\ldots\circ f_{t_1}(z) \in M_2^c
$
with $j-i-1\ge0$ because we took $j\ge N^\star$, a contradiction.
But then we get
$
f_{t_j}\circ\ldots\circ f_{t_1}(z) \in
f^{-1}_{t_{j+1}}\circ\ldots\circ f^{-1}_{t_{j+i}}(M_1),
$
i.e.,
$
f^{j+i}(z,\subl{t})\in M_1,
$
and item Vc says $f^{j+i+k}(z,\subl{t})\in M_1$
for all $k\ge0$ with
$\SQ\subset U\cup\SQ^\star \subset M_1^c$.
That is, $\om(z,\subl{t})\cap\SQ=\emptyset$,
contradicting the choice of $z$ and $\subl{t}$.

We have show \refer{confinado} to hold for $j\ge N^\star$,
since $\SQ^\star \subset L$.
However, by item VIa, we know
$f^j(z,\subl{t})\in L$ for $1\le j\le N_Q$,
where $N_Q > N^\star$
by item VII.
\end{proof}

\begin{remark}\label{Qlinha}
The arguments above show that if we replace $N^\star$
by $N$ and assume $\til{\SQ}=\SQ$
as in item VII, then writing $\SQ'$ for this new
neighborhood of the portion of $\SO$ outside $U$, we
may ensure under the same conditions of lemma~\ref{confina}
that $f^j(z,\subl{t})\in U\cup \SQ'$ for all $j\ge N$.
\end{remark}

This {\em confinement\/} property in turn implies

\begin{lemma}\label{empeh}
For every given $b_0>0$, $c_0>0$ and $\si>1$ there are
\begin{itemize}
\item a sufficiently small compact neighborhood
$\SQ\subset\SQ^\star\subset L$ of $q$, and
\item a small enough $t^\star>0$
\end{itemize}
such that $N_Q$ of item VIa be  big enough
in order that whenever
\begin{itemize}
\item $v_0\in T_{z_0} M^2$ with $z_0\in\SQ$;
\item $\subl{t}=(t_j)_{j=1}^\infty$ is a sequence
satisfying $|t_j|<t^\star$, $j\ge1$, and
\item there is $k\in\natur$ such that $N_Q\le k <\infty$
is the first integer satisfying $f^k_{\subl{t}}(z)\in \SR$;
\end{itemize}
then we have
\begin{enumerate}
\item $\slope{v_0}\ge c_0 \impl$
$\slope{ Df^k_{\subl{t}} (z_0)v_0 }\ge b_0 $ and

\item $\| Df^k_{\subl{t}} (z_0)v_0 \| \ge \si \| v_0 \|$,
\end{enumerate}
where $\|\cdot\|$, the maximum norm on $L\subset\real^2$,
and the slope are to be measured in the
linearizing coordinates
given by $\vfi_0: L\mapto M^2$.
\end{lemma}

In other words, every vector sufficiently away from
the tangent directions of $\SH^s$ at $\SQ$ will
keep pointing away from $\SH^s$ when it first
arrives at $\SR$, i.e., there are no folds
in between by the action of $f_t$.

\begin{proof}
By items I through VII of section~\ref{adaptinglinear}
there is an expanding cone field $\SC^u$ defined over
$U\cup \SQ^\star \cup L$ respected by all $f_t$ with
$|t| < t^\star$ outside of $\SR$.
It may be seen as a cone field centered around
the tangent vectors to $\SH^u$, and we may assume
that vectors in $\SC^u$ at points of $L$ have slope
$\ge b_0$, since $\SH^u$ is given by $x=cont.$
in the domain $L$ of the coordinate chart $\vfi_0$.

\medskip

We let $v_0\in T_{z_0} M^2$, $z_0\in\SQ$,
$\subl{t}$ and $N_Q \le k < \infty$ be as in
the statement of the lemma.
If $\slope{v_0)\ge c_0$, then by VIa it holds that
$z_{N_Q}=f_{\subl{t}}^{N_Q} (z_0) \in L$ and
$v_{N_Q}= Df_{\subl{t}}^{N_Q} (z_0) v_0 \in \SC^u(z_{N_Q}}$.
Indeed by~\refer{linearconjugacy} we have
$\slope{v_{N_Q}}\ge C^{N_Q} \cdot \slope{v_0}$,
where $C\approx \si_0\la_0^{-1} >1$,
and $N_Q$ may be taken sufficiently big
according to item VI, by shrinking $Q$
and $t^\star$.
Likewise we may arrange for
$\| v_{N_Q} \| \ge \si \| v_0 \|$
to hold.

If $k=N_Q$, then the lemma is proved.
Otherwise we can write
$z_k=f_{\subl{t}}^k(z_0)=f_{\subl{s}}^{k-N_Q}(z_{N_Q})
\in \SR$
where $\subl{s}=\si^{k-N_Q}\subl{t}$
and
$v_k=Df_{\subl{t}}^k(z_0)v_0 =
Df_{\subl{s}}^{k-N_Q}(z_{N_Q})v_{N_Q}$.
Moreover, lemma~\ref{confina}, the construction
of $\SC^u$ and the definition of $k\ge1$
as the first iterate to arrive at $\SR$ together
imply that the iterates
$v_{N_Q},\ldots,v_{k-1},v_k$
are all in the respective cones of $\SC^u$,
and therefore
$\slope{v_k}\ge b_0$ and
$\| v_k \| \ge \| v_{N_Q} \| \ge \si \| v_0 \|$.
\end{proof}

\medskip

Now for the effect of the tangency in $\SQ$,
recalling that the slope and norm are measured
in the $\vfi_0$ coordinates.

\begin{lemma}\label{deita}
Given $\ze>0$ there is $b_0>0$ such that
for all sufficiently small compact neighborhoods
$\SQ$ of $q$ and small $t^\star>0$ it holds for
every $|t|<t^\star$ that
\[
\left.
\begin{array}{l}
z\in\SR, v\in T_z M^2
\\
\text{and}\;\;\slope{v}\ge b_0
\end{array}
\right\} \impl
\left\{
\begin{array}{l}
\slope{Df_t(z)v}\le \ze
\;\;\text{and}
\\
\| Df_t(z)v \| \ge \frac{\al}{100} \cdot \| v \| .
\end{array}
\right.
\]
\end{lemma}

\begin{proof}
We take $z\in\SR$, $v\in T_z M^2$ and $\ze>0$.
By the differentiability of $\vfi_t$ with respect to
$t$ we know that
$\tra{f}_t=\vfi_0^{-1}\circ f_t \circ\vfi_0$
has the same local expression~\refer{formaemR}
as $\til{f}_t$.
We may suppose $\vfi_0(z)=(z,y+1)$
and $D\vfi_0(v)=(v_1,v_2)$ and
derive from~\refer{formaemR} that
\begin{eqnarray*}
\slope{ D\vfi_0^{-1}( f_t(x,y+1) )
Df_t(x,y+1) (v_1,v_2)
} =
\\
= \left|
\frac{
[2\be y + D_3 H_2(t,x,y)]\cdot v_2 +
[\ga + D_2 H_2(t,x,y)]\cdot v_1
}{
[\al + D_3 H_1(t,x,y)]\cdot v_2 +
[\rho+ D_2 H_1(t,x,y)]\cdot v_1
}
\right|
\\
\le
\frac{
|2\be y + D_3 H_2(t,x,y)| +
|\ga + D_2 H_2(t,x,y)|\cdot |v_1/v_2|
}{
|\,
|\al + D_3 H_1(t,x,y)| -
|\rho+ D_2 H_1(t,x,y)|\cdot |v_1/v_2|
\,|
}.
\end{eqnarray*}

If $\slope{v_1,v_2}\ge b_0$ then we can write
\[
\le
\frac{
|2\be y + D_3 H_2(t,x,y)| +
|\ga + D_2 H_2(t,x,y)|\cdot b_0^{-1}
}{
|\,
|\al + D_3 H_1(t,x,y)| -
|\rho+ D_2 H_1(t,x,y)|\cdot b_0^{-1}
\,|
}. 
\]

We easily see that if $b_0$ is big enough
and $\de_0>0$ in item VI is small enough,
then since $\al\cdot\be\cdot\ga\neq 0$
the last quotient approximates
$|2\be y|/|\al| = |2\be\al^{-1}|\cdot|y|$,
which can be made smaller then any positive
$\ze>0$ by shrinking $\SR$ via
taking $\SQ$ and $t^\star>0$ smaller.
Moreover  making the compact neighborhood
$\SQ$ of $q$ and $t^\star>0$ smaller just
enables $\de_0$ to be smaller,
so we are safe.

The denominator in the last quotient has a
modulus bigger than
\begin{eqnarray*}
& &
|\,
|\al + D_3 H_1(t,x,y)|\cdot |v_2| -
|\rho+ D_2 H_1(t,x,y)|\cdot |v_1|
\,|
\ge
\\
&\ge&
|v_2|\cdot |\,
|\al + D_3 H_1(t,x,y)| -
|\rho+ D_2 H_1(t,x,y)|\cdot |v_1/v_2|
\,|
\\
&\ge&
\| (v_1,v_2) \| \cdot |\,
|\al + D_3 H_1(t,x,y)| -
|\rho+ D_2 H_1(t,x,y)|\cdot b_0^{-1}
\,|
\\
&\ge&
\frac{\al}{100} \cdot \| (v_1,v_2) \|
\end{eqnarray*}
since $\al\neq 0$ and
$|D_3 H_1|, |D_2 H_1|$ and $b_0^{-1}$
may be made very small.
Also
$|v_2|=\max \{ |v_1|,|v_2| \}$
because we may take $|v_2/v_1|\ge b_0 >1$.
This provides the result on the norm.
\end{proof}

\medskip

We let $t_0,\ep>0$ be such that
$|t|<t^\star$ and
$\ep < \min\{ |t| , |t^\star - t_0| \}$
as in the statement of Theorems~\ref{teorema1}
and~\ref{teorema2} and observe the following.

\medskip

\begin{remark} \label{slope+speed}
Expression~\refer{formaemR} for
$\til{f}_{t|\SR}$
implies there are $l_0, \eta>0$ such that
the smooth curve
$c_z: T=[t_0-\ep,t_0+\ep] \mapto M^2$,
$t \mapsto f(z,t)$
has slope $\ge\eta$ and velocity $\ge l_0$
at every point $c_z(t)$ independently of
$z\in\SR$ and $t\in T$.
\end{remark}

If we make $\ze=\eta/3$ we get,
by lemma~\ref{deita},
a $b_0>0$ such that this lemma 
holds for all sufficiently small $\SQ$
and $t^\star$.
Setting $c_0=\eta$ and using the $b_0$ just
obtained, lemma~\ref{empeh} holds for
every sufficiently small $t^\star$ and $\SQ$.
We note that~\refer{confinado} of lemma~\ref{confina},
on which both lemmas~\ref{empeh} and~\ref{deita}
rest, still holds if we shrink $\SQ$ and $t^\star$ and,
moreover, lemmas~\ref{empeh} and~\ref{deita}
are independent of each other.

\medskip

{\em Hence there are a compact neighborhood $\SQ$
of $q$ and $t^\star>0$ such that both
lemmas~\ref{empeh} and~\ref{deita} hold
with some $b_0>0$ and $c_0=\eta, \ze=\eta/3>0$.}

\medskip

We are now ready for the

\begin{rproof}{\ref{fisica}}
We let $w\in\SQ$ be a regular point with respect to
$\SF_{t_0,\ep}$ according to definition~\ref{def:regular}
and pick some $\subl{t}\in\De=\De_\ep(t_0)$ and $n\ge1$.
Then $w_n=f^{R_n}(w,\subl{t})\in\SQ$ and
$z=f^{R_n-1}(w, \subl{t}) \in \SR$.
Moreover since $w$ is regular, its perturbed orbits
$\SO(w,\subl{s})$ have the same return times
to $\SQ$ independently of $\subl{s}\in\De$,
and so $c_z$ is a smooth curve in $\SQ$
with slope $\ge c_0=\eta$ and speed $\ge l_0$.

Setting $\subl{s}=\si^{r_n} \subl{t}$ then
$c=f^{r_{n+1}-1}_{\subl{s}}\circ c_z:
t\in T \mapsto f^{r_{n+1}-1}
( c_z(t) , t_{r_n+1} , \cdots , t_{r_{n+1}-1} )$
is a curve in $\SR$ with slope $\ge b_0$
and speed $\ge \si_0 l_0$ by lemma~\ref{empeh},
whereas, by lemma~\ref{deita},
$f_u\circ c$ is a curve in $\SQ$
with slope $\le \ze = \eta/3$ and speed
$\ge \frac{\al}{100}\si_0 l_0$
for all $u\in T=[t_0-\ep,t_0+\ep]$.

\begin{figure}[h]
 
\centerline{\psfig{figure=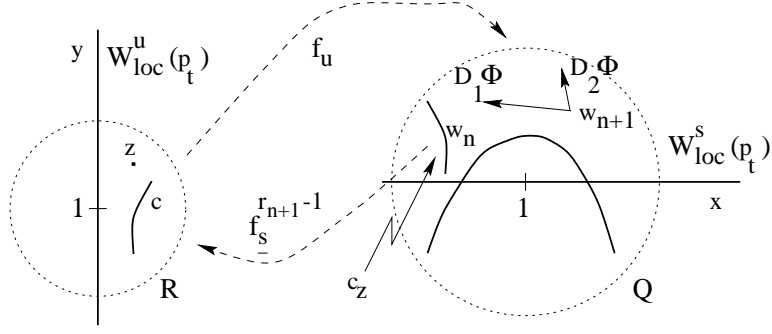,width=4.0in}}
 
\caption{\label{dobra} The iterations in the proof of
proposition~\ref{fisica}.}
\end{figure}

The regularity of $w$ implies
$\Phi(t,u)=f( c(t) , u)$
to be such that
$\Phi(t,u)\in f^{R_{n+1}}(w, \De) \subset \SQ$
for every $(t,u)\in T\times T$.
In short we have
\begin{equation}\label{fisica-geom}
\left\{
\begin{array}{l}
\slope{ D_1 \Phi } \le \eta/3
\\
\| D_1 \Phi \| \ge \frac{\al}{100}\cdot\si_0 l_0
\end{array}
\right.
\;\;\text{and}\;\;
\left\{
\begin{array}{l}
\slope{ D_2 \Phi } \ge \eta
\\
\| D_2 \Phi \| \ge l_0
\end{array}
\right.   .
\end{equation}

Noting that $D\Phi$ is the derivative of
$f^{R_{n+1}} (w,.)$ with respect to the
$R_nth$ and $R_{n+1}th$ coordinates at $\subl{t}$,
we have
$D\Phi=D_{R_n,R_{n+1}} f^{R_{n+1}} (w,\subl{t}):
\real^2 \mapto T_{w_{n+1}} M^2$
is a surjection for every $\subl{t}\in\De$.
We conclude that
$f_w^{R_{n+1}}:\De \mapto M^2$,
$\subl{t}\mapsto f^{R_{n+1}} (w,\subl{t})$
is a submersion.
This immediately gives~\ref{fisica}(2)
by definition of $f^{R_{n+1}} (w,\nuinf)$,
because the inverse image by a submersion
preserves sets of measure zero.

Making $\subl{t}=(s,s,s,\ldots)\in\De$
for some $s\in T$, since the bounds in~\refer{fisica-geom}
do not depend on $\subl{t}$, we deduce from
$f^{R_{n+1}}(w, \De) \supset \Phi(T\times T)$
that there is $\xi_0=\xi_0(s)>0$
such that $f^{R_{n+1}}(w, \De)$ contains a ball
of radius $\xi_0$ around
$\Phi(s,s)= f_s^{R_{n+1}}(w)$
as stated in~\ref{fisica}(1).
\end{rproof}

\section{Regularity of Limit Points}
\label{regular:prova}

Let $z\in\SQ$ be $V$-recurrent with $\nuinf(V)>0$
and let $\subl{t}$ be a $V$-generic vector and
$w\in\om(z,\subl{t})$.

\begin{claim}\label{cl-A}
If some $\subl{\te}\in\De$ takes $w$ to $\SQ$
after $k\ge1$ iterates,
then every other $\subl{\vfi}\in\De$ must
do the same.
\end{claim}

Indeed, if $k\ge1$ and $\subl{\te}\in\De$ are
such that $f^k(w,\subl{\te})\in\SQ$ and there is
$\subl{\vfi}\in\De$ such that $f^k(w,\subl{\vfi})\not\in\SQ$,
then we must have
$f^{k-1}(w,\subl{\te})\in\SR$
and $f^{k-1}( w,\subl{\vfi})\not\in\SR$.

By connectedness of $T^{k-1}$ and continuity of
$f^k$ (v. property~\ref{prop0}) there must be
$\subl{\psi}\in\De$ such that
$f^{k-1}(w,\subl{\psi})\in\SA$.
Since $\SA$ is open and $w\in\om(z,\subl{t})$
with $\subl{t}$ $V$-generic,
we may find for small $\de>0$ a $n\in\natur$
(according to lemma~\ref{lema0-1>0})
such that for every $\subl{s}\in V$ satisfying
$d(\subl{s},\subl{t})<\de$,
$d_M( f^n(z,\subl{s}) , w ) < \de$
and $d(\si^n\subl{s}, \subl{\psi}) < \de$
it holds that
$f^{n+k-1}(z,\subl{s})\in\SA$, and so
$f^{n+k}(z,\subl{s}) \in ( U\cup \SQ')^c$.
Moreover, these points form a set of positive
$\nuinf$-measure.

According to remark~\ref{Qlinha}
(the $n$ above can be made arbitrarily big,
bigger than $N$ in particular),
those $\subl{s}$ cannot define a perturbed orbit
$\SO(z,\subl{s})$ with infinitely many
returns to $\SQ$, which contradicts the assumptions
on $z$ and $V$.

The previous arguments readily prove

\begin{claim}\label{cl-B}
The orbit of $w$ under any $\subl{\te}\in\De$
cannot fall outside of $U\cup\SQ'$.
\end{claim}

\begin{claim}\label{cl-C}
If some $\subl{\te}\in\De$ keeps the orbit
$\SO(w,\subl{\te})$ inside $U$ for all $k$th
iterates with $k\ge k_0$,
then every other $\subl{\vfi}\in\De$ must
do likewise.
\end{claim}

In fact, if $\subl{\te}\in\De$ is such that
$f^k(w,\subl{\te}) \in U$ for all $k\ge k_0$
for some $k_0\in\natur$ and there are
$k_1\ge k_0$ and $\subl{\vfi}\in\De$
such that $f^{k_1}(w,\subl{\vfi})\not\in U$,
then by the connectedness of $T^{k_1}$,
property~\ref{prop0} and the separation between
$U$ and $\SQ'$ given by item VII,
there is $\subl{\psi}\in\De$
satisfying $f^{k_1}(w,\subl{\psi}) \in (U\cup\SQ')^c$.
We may now repeat the arguments proving the
preceding claim.

\medskip

For $w\in\om(z,\subl{t})$ with $\subl{t}$ a $V$-generic vector
we have the following alternatives:
\begin{enumerate}
\item $w$ returns to $\SQ$ a finite number of times only
under every $\subl{\te}\in\De$;
\item $w$ never passes through $\SQ$ under every
$\subl{\te}\in\De$;
\item $w$ returns to $\SQ$ infinitely often
and $r(w,\subl{s},n)=r(w,n)$, $\subl{s}\in\De$,
$n\ge1$.
\end{enumerate}

Since $w$ cannot get out from $U\cap\SQ'$ by claim~\ref{cl-B},
alternatives 1 and 2 imply that the orbits of $w$ stay forever
in $U$ after some finite number of iterates or never leave $U$,
respectively.
For our purposes it is enough to suppose
$w\in\om(z,\subl{t})\cap\SQ$.

\subsection{Finite Number of Returns}
\label{alt1}

First we eliminate alternative 1.
By claims~\ref{cl-A} and~\ref{cl-C}
the return times to $\SQ$ and the iterate after which the
orbits remain forever in $U$ do not depend on the
perturbation vector.

\medskip

Let $r_0\in\natur$ be the last return iterate of $w$ to
$\SQ$ under every $\subl{\te}\in\De$.
The point $w$ is like a
{\em regular point up to iterate $r_0$\/}
and so the arguments in section~\ref{fisica:prova}
show that $f^{r_0}(w,\De)$ contains a curve $c$
with slope $\ge\eta$ and speed $\ge l_0$ at every point.
So its length is
$\ge 2\ep\cdot l_0 = a_0 >0$
and since $w\in\om(z,\subl{t})$, no orbit
is allowed to leave $U\cup\SQ'$.
Hence $f^k(c,\De)\subset U$ for all $k\ge1$.
In particular,
$c_k=f^k_{t_0}(c)=f^k(c,\subl{t_0}) \subset U$,
$k\ge1$.

According to the previous section, after $N_Q$
iterates curve $c$ will have all its tangent
vectors in $\SC^u$ and keep them this way for all
iterates onward, because $c_k\subset U$ for all $k\ge1$.
Since $\SC^u$ is a field of unstable cones,
the length of $c_k$ will grow without bound
with $c_k$ being an {\em unstable curve\/} always
inside $U$.

This is a contradiction, since $U$ is a small neighborhood
of a hyperbolic set $\La_{t_0}$ of saddle type which is
the maximal invariant set inside $U$.

\subsection{No Returns}
\label{alt2}

\begin{figure}[h]
 
\centerline{\psfig{figure=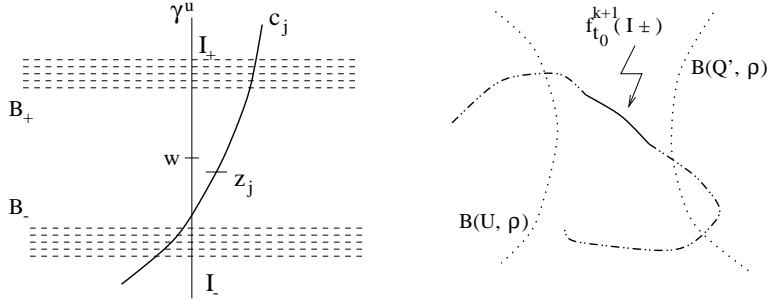,width=4.0in}}
 
\caption{\label{faixaslimite} The situation near $w$ and the
image of $I_{\pm}$}
\end{figure}

Let $w$ be as in alternative 2.
Consequently $f^k_{t_0}(w)\in U$
for all $k\ge1$.
Since $\La_{t_0}$ is the maximal invariant set
inside $U$, we deduce that if $\ga^u$ is a small
segment of $\SH^u_{t_0}(w)$ centered at $w$,
then it is not possible
that $f^k_{t_0}(\ga^u)\subset U$ for all $k\ge1$.
Likewise if we replace $U$ by $B(U,\rho)$,
by item VII.
Hence, writing $\ga^u_{+}, \; \ga^u_{-}$ the two
segments such that
$\ga^u_{+} \cup \ga^u_{-} = \ga^u$ and
$\ga^u_{+} \cap \ga^u_{-} = \{ w \}$,
there are $k_{\pm}\ge1$ and nonempty intervals
$I_{+}\subset \ga^u_{+}$,
$I_{-}\subset \ga^u_{-}$ satisfying
$f^i_{t_0}(I_{\pm})\subset B(U,\rho)$
for $1\le i \le k$ and
$f^{k+1}_{t_0}(I_{\pm}) \subset
( \tra{B(U,\rho)} \cup \tra{B( {\SQ}' ,\rho)} )^c$
--- because
$\tra{B(U,\rho)} \cap \tra{ B(\SQ' ,\rho) } = \emptyset$
and by connectedness of
$\ga^u_{\pm}$ (v. figure~\ref{faixaslimite}).

Let $x\in I_{\pm}$ and $y\in \SH^s_{t_0}(x)$.
Then we have
$d_M (f^k_{t_0}(x),f^k_{t_0}(y)) \le C\la^k d_M(x,y)$
where $1>\la \ge |\la_t|$ for $|t|<t^\star$.
So every $y\in \SH^s_{t_0}(x)$ with
$d_M(x,y) \le C^{-1}\la^{-k}\cdot \rho/2$
satisfies $f^k_{t_0}(y) \in (U\cup\SQ')^c$.


Geometrically this means that near $w$ there are two strips
$B_{\pm}$ made of $\SH^s_{t_0}$-leaves with length
$C^{-1}\la^{-k}\cdot \rho/2$ and whose intersection
with $\ga^u$ is $I_{\pm}$ (cf. figure~\ref{faixaslimite}).

Making $\ga^u$ small and $k$ big we can make the length of
$B_{\pm}$ big and the distance to $w$ small. The angle between
leaves in $B_{\pm}$ and $\ga^u$ is near a straight angle
in the $\vfi_{t_0}$-coordinates of $L\supset\SQ\ni w$,
since the slope of $\SH^s_{t_0}$ is near $0$.

\medskip

Let $n_1<n_2<n_3<\ldots$ be such that
$z_j=f^{n_j}(z,\subl{t}) \mapto w$ when $j\to\infty$.
We define $c_j: T \mapto L$,
$u\mapsto f_u( f^{n_j-1}(z,\subl{t}) )$,
the {\em perturbation curve\/} through $z_j$ and
observe that either $B_{+}$ or $B_{-}$ intersects
$c_j(T)$ in a segment of positive length $\ge a_1>0$,
since $\slope{c_j'(u)}\ge\eta$ and the length of
$c_j$ is $\ge a_0>0$, for all $u\in T$ and $j\ge1$.

This means there is a segment $S_j$ of length $\ge a_2>0$
in $T$ such that $c_j(S_j)\subset B_{\pm}$
and thus
$f^{k+1}_{t_0}(c_j( S_j )) \subset ( U\cup \SQ' )^c$.

According to lemma~\ref{lema0-1-2},
for every $0<\ga,\de<1$ we can find $k_0\in\natur$
such that for all $j\ge k_0$ we have
$\nu ( p_{n_j} V_{n_j-1} ( \subl{t}, n_j-1, \subl{s} ) ) \ge 1-\de$
for a positive measure set $V_{n_j-1}\subset V$
and a set of $\subl{s}\in\De$ with $\nuinf$-measure
$\ge 1-\ga$.
Hence, since $k$ is fixed,
we may find for $j$ big a $\subl{s}\in\De$
very close to $\subl{t_0}=(t_0,t_0,\ldots)$
(taking $\ga>0$ small)
such that
$\nu( S_j \cap p_{n_j} V_{n_j-1}( \subl{t}, n_j-1, \subl{s}) )>0$
and
$f^{k+1}_{\subl{s}} ( c_j(S_j) ) \subset (U\cup\SQ' )^c$.

\medskip

We have shown that inside $V$ there is a positive measure
set whose perturbation vectors send $z$ into
$(U\cup\SQ' )^c$ after $n_j+k+1$ iterates,
where $j$ (and $n_j$) may be made arbitrarily big.
This contradicts the assumption of $V$-recurrence
on $z$, since those perturbed orbits will never again
return to $\SQ$. Alternative 2 is thus impossible.

\subsection{Bounded First Return Times}

The points $w\in\om(z,\subl{t})\cap\SQ$ with $\subl{t}$
a $V$-generic vector satisfy alternative 3.
Going back to the arguments in subsection~\ref{alt1},
we have an unstable curve $c$ in $f^k(w,\De)$ whose
length cannot grow unbounded.
Therefore it must leave $U$ and go to $\SQ$
(since $w$ no orbit may leave $U\cup\SQ' $)
after a finite number or
iterates bounded by some $J\in\natur$.
We observe that since the length of $c$ is
$\ge 2\ep \cdot l_0$
and the diameter of $\SR$ is finite,
we must have $(2\ep l_0) \cdot \si^J \approx \diame{\SR}$.

\medskip

This proves proposition~\ref{regular}
and Theorem~\ref{teorema2}.

\begin{remark}\label{semilocal}
We may drop the first tangency condition of
subsection~\ref{1parfam} if we strengthen
definition~\ref{def:Vrecorrente} of $V$-recurrent point
by adding the following item
\begin{enumerate}
\item[3.] for $\nuinf$-a.e $\subl{t}\in V$
there is $n=n(\subl{t})\ge1$ s.t.
$f^k(z,\subl{t})\in U\cup\SQ$ for all $k\ge n$;
\end{enumerate}
where $U$ is a fixed neighborhood of the basic
set $p_0$ belongs to and $\SQ$ a neighborhood
of the piece of the orbit of tangency outside $U$.

Lemma~\ref{confina} is now needless and the rest
of the proof is unchanged. The scope of the theorem
is enlarged and next section shows how this extra
condition on $V$-recurrence is not too restrictive.
\end{remark}

\section{Infinitely many attractors}
\label{sinks}

We start with the particular case 
of perturbations of sinks.

\begin{definition}\label{defpertsink}
We say $f\in \text{Diff }^l(M)$, $l\ge1$, has a
{\em perturbation of a sink\/}
in a finite collection $\domcomplinv$
of pairwise disjoint open sets of $M$
if there exists a neighborhood $\SV$ of
$f$ in $\text{Diff }^l(M)$ such that,
for every continuous arc
$\SG=\{ g_t \}_{t\in B} \subset \SV$
with $g_0\equiv f$,
the following holds:
\begin{enumerate}
\item $\tra{ g_{\subl{t}}^n(\SU_i) } \subset \SU_{(i+n) \bmod r}$
for every $n\ge1$, $\subl{t}\in B^\natur$ and $0\le i \le r-1$ ;
\item there is a constant $\be>0$ such that for every point
$x\in\SU_i$, $0\le i \le r-1$,
every $v\in T_x M \setminus \{0\}$ and
every $\subl{t}\in B^\natur$ it holds that
\[
\limsuppn \log \| Dg_{\subl{t}}^n (x) \cdot v \| \le -\be ;
\]
\item with the notation introduced at definition
\ref{defswlimite} we have
\[
\diame{\om \left( \SU_j,\De_\ep (0)  \right) \cap \SU_i}
\longrightarrow 0 \;\; \text{when} \;\; \ep\to 0^{+}
\]
for every $0\le i \le r-1$.
\end{enumerate}
[Where $B=\tra{B}^j (0,1)$ and $\De_\ep (0) =( \tra{B}^j(0,\ep) )^\natur$
as in subsection \ref{pertaroundpar}].
\end{definition}

Next proposition characterizes this kind of invariant domains.

\begin{proposition}\label{pococomruido}
Let $f$ be a $C^l$ diffeomorphism of $M$, $l\ge1$.
Then $f$ has a hyperbolic sink $s_0$ with period $k\ge1$
if, and only if,
$f$ has a perturbation of a sink in a neighborhood
$\domcomplinv$ of the orbit
$s_0, s_1=f(s_0),\ldots,s_{r-1}=f^{r-1}(s_0)$
of $s_0$.
\end{proposition}

\begin{proof}
First some results that locate the
limit points near a perturbed sink.

If $s_0\in M$ is a hyperbolic sink
for $f$ with period $r$, then
for some $0 < \la_1 <1$ every eigenvalue $\la\in\complex$
of $Df^r(s_0)$ satisfies $| \la | \le \la_1$.
Moreover, given some $\la_1 < \ups <1$ there are
$\de>0$ and a neighborhood $\SV$ of $f$ in
$\text{Diff }^l(M)$ -- both may be made arbitrarily small --
such that
each eigenvalue $\la\in\complex$ of
$Dg^r(x)$ satisfies $|\la| \le \ups$ for every
$g\in\SV$ and $x\in B(s_0,\de)$.
Consequently
\begin{equation}\label{contraccao}
d_M ( g^r(x), g^r(y) ) \le \ups\cdot d(x,y)
\;\; \text{for every} \;\;
x,y\in B(s_0,\de)
\;\;\text{ and }\;\;
g\in\SV.
\end{equation}

So,
writing $s_i=f^i(s_0)$, we see that
$
\left( \SU_0=B(s_0,\de) , \ldots ,
\SU_{r-1}=B(s_{r-1},\de)  \right)
$
is a finite collection of pairwise disjoint
(we may take $\de<\frac12
\min\{ d_M(s_i,s_j) : 0\le i<j \le r-1 \}$)
open sets of $M$ that satisfies conditions
1 and 2 of definition \ref{defpertsink}.

To get condition 3 we have the next

\begin{lemma}
Let $\SG=\{g_t\}_{t\in I} \subset \SV$ be some
continuous arc in $\text{Diff }^l(M)$ with $g_0\equiv f$.
Let
$P_i=\{ s_i(t) : t\in B \}$
be the set of analytic continuations
of the orbit $\SO(s_0)$ of the sink $s_0$
with respect to $g_t$, $t\in B$.

If we fix $x\in\SU_i$, $0\le i \le r-1$,
and $\subl{t}\in \De$, then we have
\[
d_M(y, P_j ) \le \frac{\ups}{1-\ups} \cdot
\max \{ \diame{ P_k  }: 0\le k \le r-1 \},
\]
for every $y\in\om(x,\subl{t})\cap\SU_j$,
$j=0,\ldots,r-1$.
\end{lemma}

\begin{proof}
This is an easy consequence of~\refer{contraccao}.
\end{proof}

\medskip

We now know that
$\om(x,\subl{t})\subset B( P, \ga )$
where
$\ga=\frac{\ups}{1-\ups}\cdot\max \{ \diame{P_h}: 0\le h \le r-1\}$
and, since $s_0$ is an hyperbolic sink for $f\equiv g_0$,
we have
\[
\diame{ \{ s_h(t) : t \in \tra{B}^j(0,\ep) \}  } \longrightarrow 0
\;\;\text{when}\;\;
\ep\longrightarrow 0^{+}
\]
by the structural stability results for such attractor.

Therefore item 3 holds for $\domcomplinv$
constructed above and we have shown
that in a neighborhood of the orbit of every
hyperbolic sink there is a perturbation of sink.

\medskip

Conversely, let us suppose $f$ has a perturbation
of a sink in some collection $\domcomplinv$
of pairwise disjoint open sets and take
$\SG=\{g_t\}_{t\in I}$ as in definition~\ref{defpertsink}.
Then we will have by definition
$
\om( \SU_i, \De_{\ep_1}(0) ) \subseteq
\om( \SU_i, \De_{\ep_2}(0) )
$
for every small $0<\ep_1<\ep_2$ and every $0\le i \le r-1$.
Property 3 of definition~\ref{defpertsink}
now ensures there is a point $s_0$ such that
$
\{ s_0 \} = \cap_{\ep>0}
[ \om( \SU_0,\De_\ep(0) ) \cap \SU_0  ]
$
since $\om( \SU_0,\De_\ep(0) )$ is a closed set.

\medskip

Writing $\subl{0}=(0,0,\ldots)$ then $\subl{0}\in\De_\ep(0)$
and $\om(s_0,\subl{0})\subset\om(\SU_0,\De_\ep(0))$
for every $\ep>0$.
Thus $\{ s_0 \} = \om(s_0,\subl{0}) \cap \SU_0 =
\om_f(s_0)\cap\SU_0$.
Considering the dynamics induced in
$\domcomplinv$ by the arc $\SG$ we see that
$
\om_f(s_0)=\{ s_0,\dots,s_{r-1} \}$
where
$s_i=f^i(s_0)$, $i=0,\ldots,r-1$.

Since the limit is $f$-invariant, we have
$f^r(s_0)=s_0$ and found a $r$-periodic orbit of $f$.
In addition, property 2 of definition~\ref{defpertsink}
guarantees that for each
$v\in T_{s_0}^1 M=\{u\in T_{s_0}M: \|u\|=1 \}$
such that $v$ is an eigenvector of $Df^r(s_0)$
corresponding to the eigenvalue $\la\in\complex$
(using the complexification of
$Df^r(s_0):T_{s_0}M \to T_{s_0}M$ if need be)
the following holds
\[
0> -\be \ge \limsup_{n\to +\infty}
\frac1{rn} \log
\left\| Df^{n\cdot r} (s_0) \cdot v \right\|
= \frac1r \log |\la| \impl
|\la| \le \exp ( -r\be ) < 1
\]
and so ${\rm sp}( Df^r(s_0) ) \subset
\{ z\in\complex : |z|<1 \}$.
Hence $s_0,\ldots,s_{r-1}$ is the orbit of an hyperbolic
sink for $f$.

\end{proof}

\subsection{Newhouse's and Colli's Phenomena}

Let us suppose the family $f$ satisfying the conditions
specified in subsection \ref{1parfam}
is also in the conditions of Newhouse's
theorem (cf.~\cite{N1,N2} and ~\cite{PT})
on the coexistence of infinitely many sinks, that is,
$p_0$ is a dissipative ($|\deter{Df_0(p_0)}|<1$)
saddle point.

\medskip

We may now choose a parameter $a>0$ such that $f_a$
has infinitely many hyperbolic sinks in $\SQ$.
Moreover $a>0$ may be taken arbitrarily close to zero
(see \cite[Capt. 6]{PT}) and thus all the results
of previous sections apply to the present setting.

\medskip

Let $N$ be some positive integer and let us pick
$N$ distinct orbits of hyperbolic sinks for $f_a$ in
$\SQ$: $\SO(s^{(i)}), \; i=1,\ldots,N$.
Since they are hyperbolic attractors, they are isolated:
there exist pairwise disjoint -- even separated --
open neighborhoods $V_i$ of $\SO(s^{(i)}), \; i=1,\ldots,N$.
Moreover, by the previous subsection, we may construct
a perturbation of a sink inside each $V_i$
associated to $\SO(s^{(i)})$ with respect to an arc
$\SF_{a,\ep_i}$, for some $\ep_i>0$, and every $1\le i \le N$.

\medskip

We now observe that a perturbation of a sink obviously is,
in particular, a completely and symmetrically invariant
domain. Specifically, each perturbation of a sink
constructed in $V_i$ is a completely and symmetrically
invariant domain with respect to the arc $\SF_{a,\ep_i}$,
$i=1,\ldots, N$.

Hence, setting $\ep_0=\min\{\ep_1,\ldots,\ep_N\}$, we
have $\ep_0>0$ and the former invariant domains are also
completely and symmetrically invariant with respect to the arc
$\SF_{a,\ep}$ for every $0<\ep<\ep_0$.
Then, by subsection \ref{minimalexist}, there is a minimal
domain $\SM^{\ep}_i$ inside each perturbation of a sink $V_i$,
for every $1\le i \le N$ and noise level $0<\ep<\ep_0$.

We have thus constructed $N$ distinct minimal invariant
domains in $\SQ$ for the arc $\SF_{a,\ep}$
for every $0<\ep<\ep_0$ and proved

\begin{proposition}\label{muitosminimais}
Given an arc $\SF$ as in subsection \ref{1parfam}
where $p_0$ is a dissipative saddle point, for every
parameter $a>0$ sufficiently close to zero such that
$f_a$ has infinitely many sinks in $\SQ$, we have
the following.

For every $N\in\natur$ there exists $\ep_0>0$ such that,
for every $0<\ep<\ep_0$, the number of minimal invariant
domains in $\SQ$ for the arc $\SF_{a,\ep}$
is no less than $N$.
\end{proposition}

We now remark that what enables us to build an
invariant domain in a neighborhood of a sink is
the fact that it is attractive: given any
neighborhood $U$ of the orbit of a sink
$s_0,\ldots,s_{r-1}$ there is another neighborhood
$V\subset\tra{V}\subset U$ of the same orbit
such that $f(\tra{V})\subset V$
(a {\em trapping region\/}). By continuity,
this persists for any diffeomorphism $g$ close
to $f$ and hence we get an invariant domain.

In~\cite{C} E. Colli shows how to have infinitely
many H\'enon-like attractors when generically
unfolding an homoclinic tangency under the same
conditions of Newhouse's theorem. These attractors
are separated like the infinity of sinks in the
Newhouse phenomenon and each one admits a
{\em trapping region\/} according to~\cite{BM}
and~\cite{V}. Specifically, the constructions
described in~\cite{C} can be carried out verbatim
within a restricted set of parameter values
having this property, without altering the
statements of any theorem in that paper.

Consequently we may state and prove a proposition
analogous to~\ref{muitosminimais}
replacing {\em sink\/} by {\em H\'enon-like attractor\/}
in the paragraphs above.

\section{Some Conjectures}
\label{conjecturas}

The methods used in this paper are prone
to generalization. We propose some here.

\begin{nliste}{\Parentalgarismo}
\item Is there some {\em similar\/} result
to Theorem~\ref{teorema1} for flows?
The kind of perturbation to perform is part
of the question.

\item In section~\ref{sinks} a characterization
is given for invariant domains originating from
a perturbation of a sink.
Is there some {\em similar\/} characterization
of an invariant domain obtained by a perturbation
of an H\'enon-like strange attractor?

\item The same question regarding perturbations
of elliptic islands. This is more subtle: we
may ask whether there is some invariant domain
near an elliptic island.

\item We did not look at what happens to the
physical probabilities when the noise level $\ep>0$
tends to zero. Does the limit exist?
If it does then it must be an $f$-invariant
probability measure. Is it an SRB-measure?

\item Globally what can we say about the
stochastic stability of the infinitely many
Dirac (in Newhouse's phenomemon) or SRB
(in Colli's phenomenon) measures in a
neighborhood of a homoclinic tangency point?
Here a global notion of stochastic stability
is required, see e.g.~\cite{V2}:
if $\mu_i$ are the SRB measures of $f$
($i=1,2,\ldots$), time averages of each continuous
$\vfi$ along almost all random orbits should
be closed to the convex hull of the $\int \vfi\,d\mu_i$
for small $\ep>0$.
\end{nliste}

{\em
Vítor Araújo

Centro de Matemática da Universidade do Porto

Rua do Campo Alegre  687, 4169-007 Porto, Portugal.

{\rm and}

Instituto de Matemática,  Universidade Federal do Rio de Janeiro,

C. P. 68.530,  21.945-970 Rio de Janeiro, RJ-Brazil


Email: vitor.araujo@im.ufrj.br {\rm and} vdaraujo@fc.up.pt

}

\end{document}